\documentclass{amsart}

\usepackage{graphicx}
\usepackage{amsfonts}
\usepackage{amsmath}
\usepackage{amsthm}

\newtheorem{Thm}{Theorem}

\newtheorem{Rule}{Rule}
\newtheorem{Lem}{Lemma}
\newtheorem{Prop}{Proposition}

\newtheorem{Def}{Definition}
\newtheorem{Rem}{Remark}

\newcommand{\If}{\mbox{~~if~~}}

\newcommand{\z}{\mathbb{Z}}

\newcommand{\TT}{\tilde{T}}
\newcommand{\Galph}{{G_\alpha}^{\Phi}}

\begin{document}
\title[Cluster Algebras of Classical Type]{A Graph Theoretic Expansion Formula for  Cluster Algebras of Classical Type}
\author{Gregg Musiker}
\date{November 5, 2007}
\thanks{This work was done with support of an NSF Mathematical Sciences Postdoctoral Fellowship.}
\maketitle

\begin{abstract}
In this paper we give a graph theoretic combinatorial interpretation for the cluster variables that arise in most cluster algebras of finite type.  In particular, we provide a family of graphs such that a weighted enumeration of their perfect matchings encodes the numerator of the associated Laurent polynomial while decompositions of the graphs correspond to the denominator.  This complements recent work by Schiffler and Carroll-Price for a cluster expansion formula for the $A_n$ case while providing a novel interpretation for the $B_n$, $C_n$, and $D_n$ cases.  
\end{abstract}

\tableofcontents

\section{Introduction}

Several years ago, Sergey Fomin and Andrei Zelevinsky introduced a new mathematical object known as a cluster algebra which is related to a host of
other combinatorial and geometric topics.  Some of these include canonical bases of semisimple algebraic groups, generalized associahedra, quiver representations, tilting theory, and Teichm\"{u}ller theory. In
the proceeding we will use the definitions and conventions used in Fomin and Zelevinsky's inital papers,
\cite{ClustI, ClustII}.  Starting with a subset $\{x_1,x_2,\dots, x_n\}$ of cluster algebra $\mathcal{A}$, one applies binomial exchange relations to obtain additional generators of $\mathcal{A}$, called \emph{cluster variables}.  The (possibly infinite) set of cluster variables obtained this way generate $\mathcal{A}$ as an algebra.
It was proven in \cite{ClustI} and
\cite{Laurent} that any cluster variable is a Laurent polynomial in $\{x_1,x_2,\dots, x_n\}$, i.e. of the form
$$\frac{P(x_1,\dots, x_n)  }{  x_1^{a_1}x_2^{a_2}\cdots x_n^{a_n}} \hspace{3em} (\mathrm{Note~that~}x_i = \frac{1}{x_i^{-1}}\mathrm{~is~also~allowed})$$
where $P(x_1,\dots, x_n)$ is a polynomial with integer coefficients (not divisible by any monomial) and the exponents $a_i$ are (possibly negative) integers.  
It is further conjectured that the polynomials $P(x_1,\dots, x_n)$ have \emph{nonnegative} integer coefficients for any cluster algebra.  However, this conjecture has been proved in a limited number of cases, including the finite type case, as proved in \cite{ClustII}, the case of rank two affine cluster algebras as demonstrated in \cite{Caldero}, \cite{MusPropp}, \cite{SherZel}, and \cite{Zel} as well as cluster algebras arising from acyclic quivers \cite{CalReit}.

The finite type case is defined as the case where the cluster variable generation procedure only yields a finite set of cluster variables associated to $\mathcal{A}$.  By a combinatorial and geometric miracle, one which has sparked much interest in these algebras, the cluster algebras of finite type exactly correspond to the Lie algebras of finite type.  Furthermore, in these cases, the cluster variables (except for the $x_i$'s) have denominators with \emph{nonnegative} exponents and can be put in a $1$-to-$1$ correspondence with the positive roots of the associated root system. 

Study of the particular finite type cluster algebra of type $A_n$, 
also known as the Ptolemy algebra has been especially fruitful as it can be 
realized in terms of the Grassmannian and Pl\"{u}cker embedding.  In 2003, as part of the REACH research group under Jim Propp's direction, Gabriel Carroll and Gregory Price 
\cite{CPpre} described two combinatorial interpretations of the 
associated cluster variables, one in terms of paths and one in terms of 
perfect matchings.  Further, Ralf Schiffler recently independently discovered and extended the paths interpretation \cite{Schiffler}.

In the present paper, we go beyond $A_n$, and describe a combinatorial interpretation for the cluster variables in all four families of finite type, namely $A_n$, $B_n$, $C_n$, and $D_n$, for the coefficient-free case.  Our combinatorial model will involve perfect matchings, in the spirit of 
\cite{MusPropp}, and agrees with Carroll and Price's interpretation in the 
$A_n$ case.  Unlike the aforementioned work we do not attempt to give the Laurent expansion of cluster variables in terms of any seed but only in terms of the initial bipartite seed, whose definition we remind the reader of below.  By restricting ourselves to expansions in this initial seed, we are able to 
explicitly write down families of graphs which encode the cluster algebra using 
weighted perfect matchings.

We shall use the following notation throughout this paper.  Let $G=(V,E)$ be a finite graph with vertex set $V = \{v_1,\dots, v_m\}$ and edge set $E \subseteq \{ \{u,v\}:u,v\in V\}$.  For each edge $e\in E$, we set $w_e$ to be the weight of $e$, where $w_e$ is allowed to be $1$ or $x_i$ for $i\in \{1,2,\dots, n\}$.  A \emph{perfect matching} $M$ of graph $G$ is a subset of $E$ such that for every vertex $v\in V$, there is exactly one edge $e\in M$ containing $v$.  The weight of a perfect matching is defined to be the product $w(M) = \prod_{e\in M} w_e$, and we let $P(G)$ denote the matching polynomial, or matching enumerator, of graph $G$, defined as $$P(G) = \sum_{M \mathrm{~is~a~perfect~matching~of~}G} w(M).$$  
The result of this paper is the following theorem.

\begin{Thm} \label{vargraph} 
Let $\Phi$ be a root system of classical type and denote its positive roots as $\Phi_+$.  For each such $\Phi$, we explicitly construct a family of graphs, $\mathcal{G}_{\Phi}$, with the following three properties.

\begin{enumerate}
\item $|\mathcal{G}_{\Phi}| = |\Phi_{+}|$.

\item For each $\alpha = (\alpha_1, \alpha_2, \dots, \alpha_n)$, there exists a unique $G_\alpha^{\Phi} \in \mathcal{G}_{\Phi}$ that can be effeciently identified.

\item We have the cluster expansion formula $$x[\alpha]^{\Phi} = \frac{P({G_\alpha}^{\Phi})}{x_1^{\alpha_1}\cdots x_n^{\alpha_n}},$$ where $x[\alpha]^{\Phi}$ denotes the cluster variable corresponding to positive root $\alpha$ (in type $\Phi$) under Fomin and Zelevinsky's bijection.
\end{enumerate}
\end{Thm}

Given graph $G \in \mathcal{G}_\Phi$, we are able to determine for which $\alpha \in \Phi_+$ we have $G = G_\alpha^{\Phi}$ by breaking down $G$ into tiles.  More precisely, we let a family of tiles $\mathcal{T}=\{T_1,\dots, T_n\}$ be a finite set of graphs, with weighted edges, such that each $T_i$ is isomorphic to a cycle graph.  Given the faces and edge weighting of graph $G$, we decompose $G$ into a union of such tiles by gluing together certain edges.

\vspace{1em}\begin{center}
$\begin{array}{cccc}
\includegraphics[width = 0.2in , height = 0.2in]{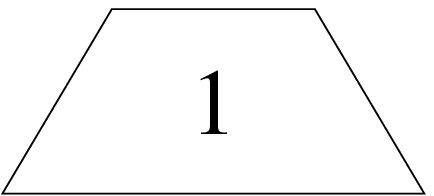}  &~~~& \includegraphics[width = 0.2in , height = 0.2in]{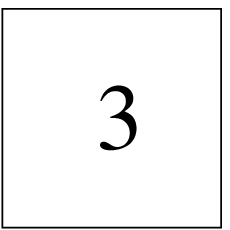} &~~~\\
~~~& \includegraphics[width = 0.6in , height = 0.45in]{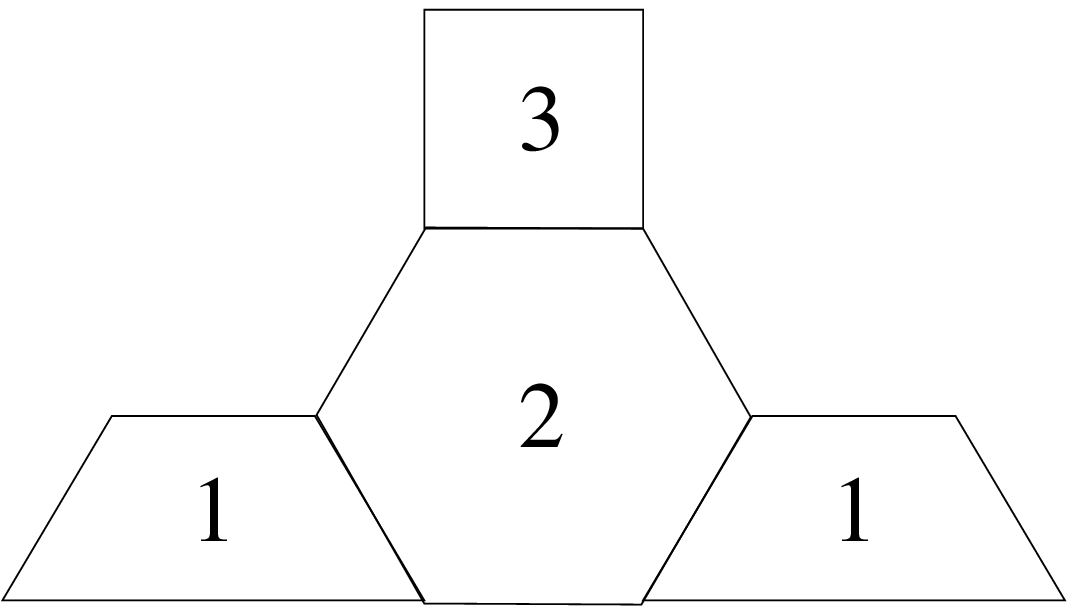} &~~~& \includegraphics[width = 0.2in , height = 0.4in]{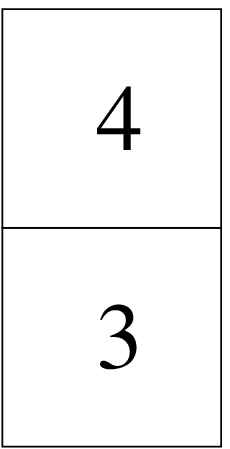}  \\
\includegraphics[width = 0.4in , height = 0.6in]{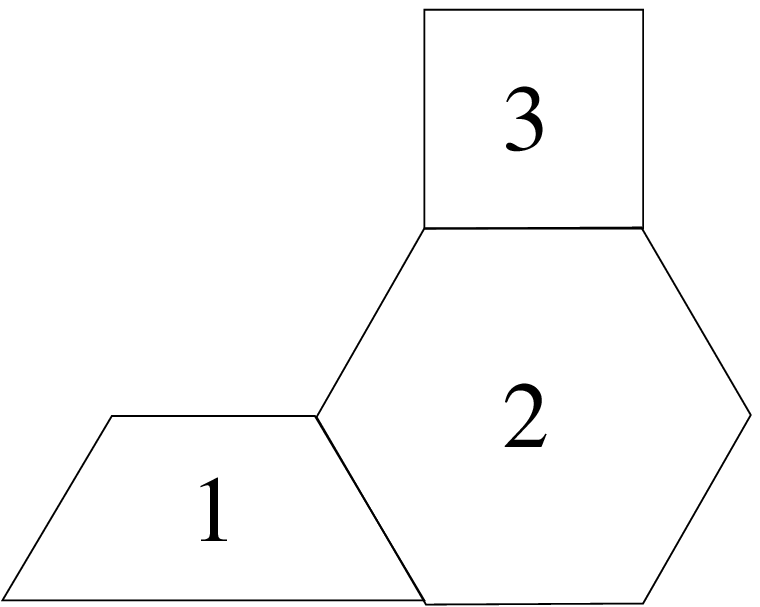}  &~~~& \includegraphics[width = 0.6in , height = 0.7in]{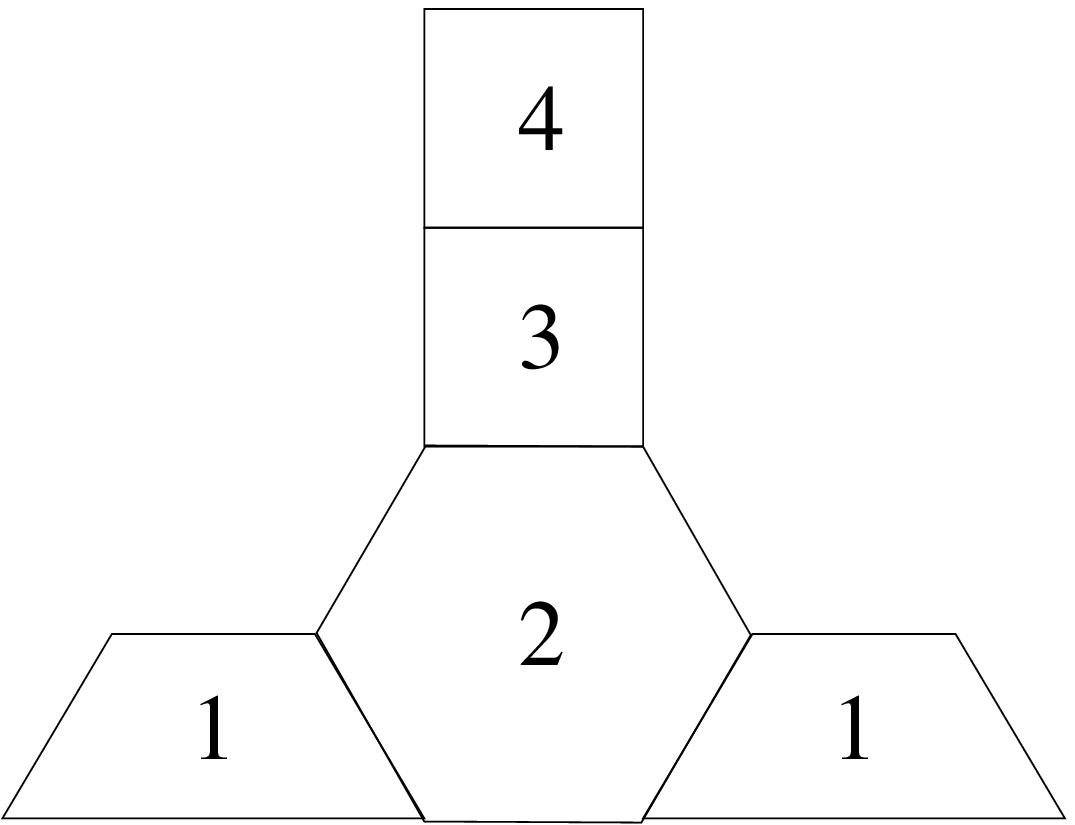} &~~~ \\
~~~& \includegraphics[width = 0.8in , height = 0.6in]{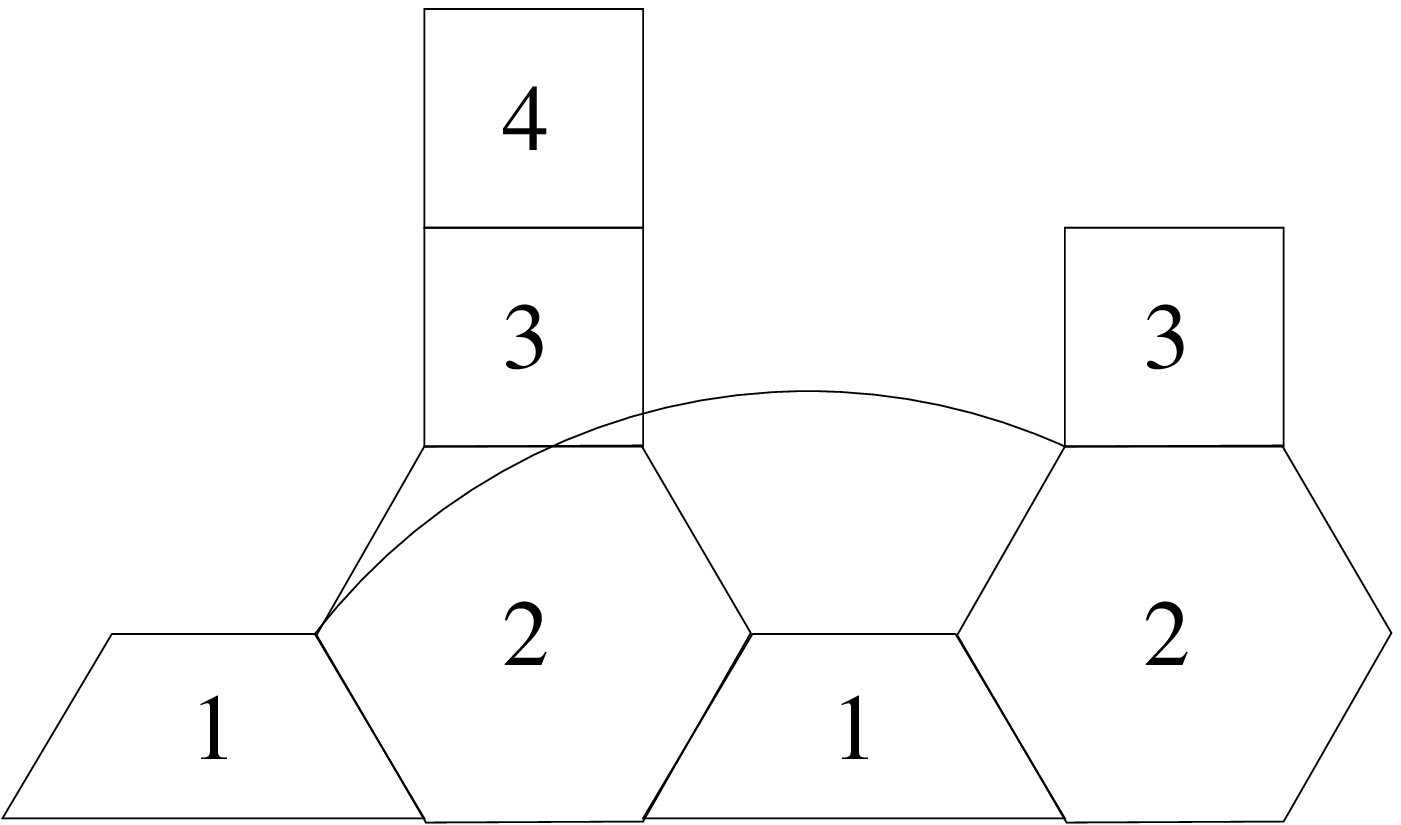} &~~~& 
\includegraphics[width = 0.6in , height = 0.4in]{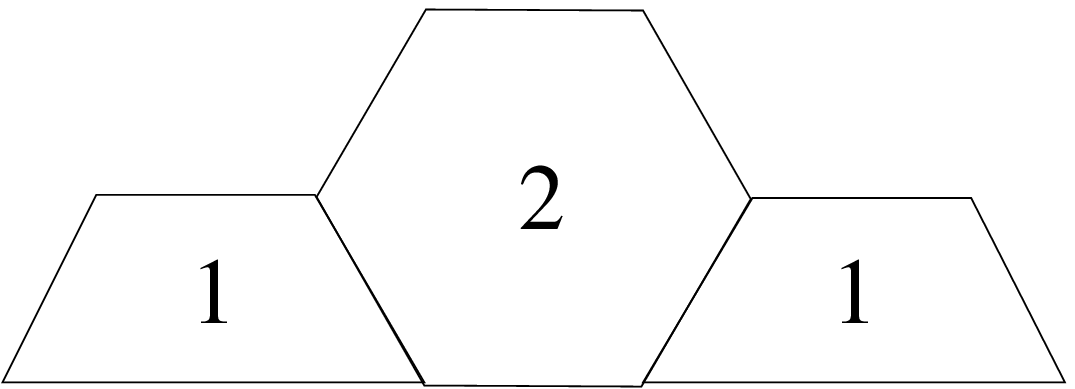} \\
\includegraphics[width = 0.4in , height = 0.6in]{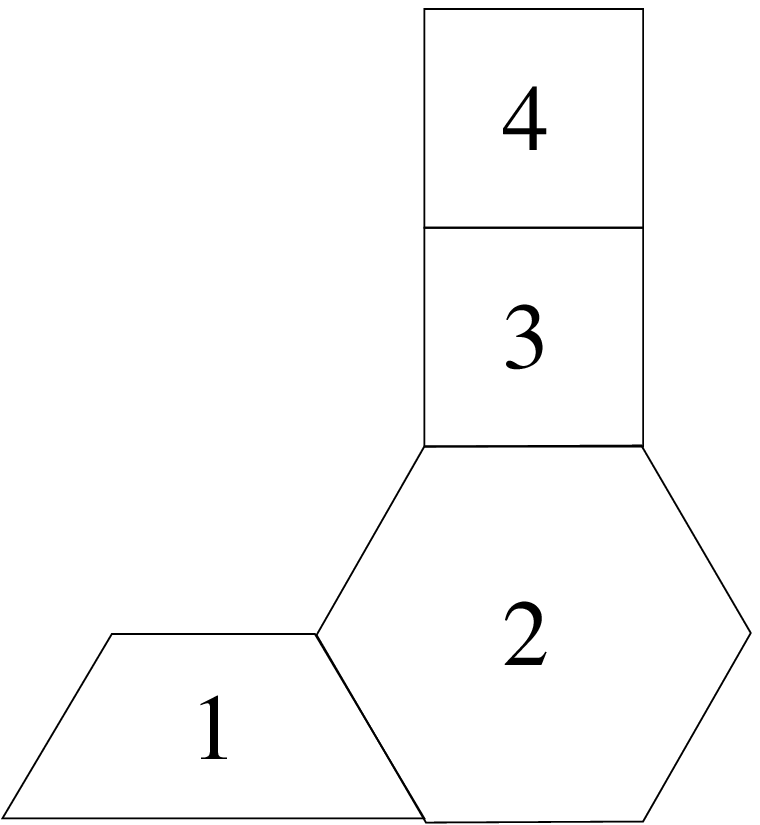}  &~~~& \includegraphics[width = 0.8in , height = 0.5in]{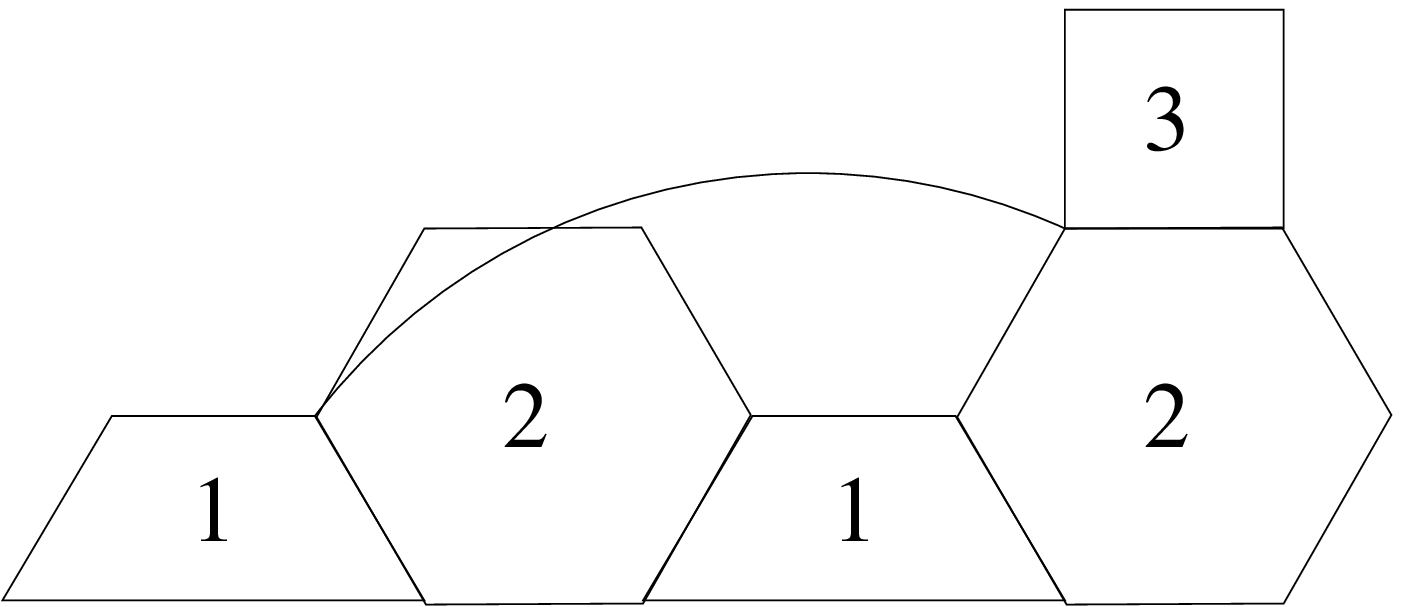} &~~~ \\
~~~& \includegraphics[width = 0.8in , height = 0.6in]{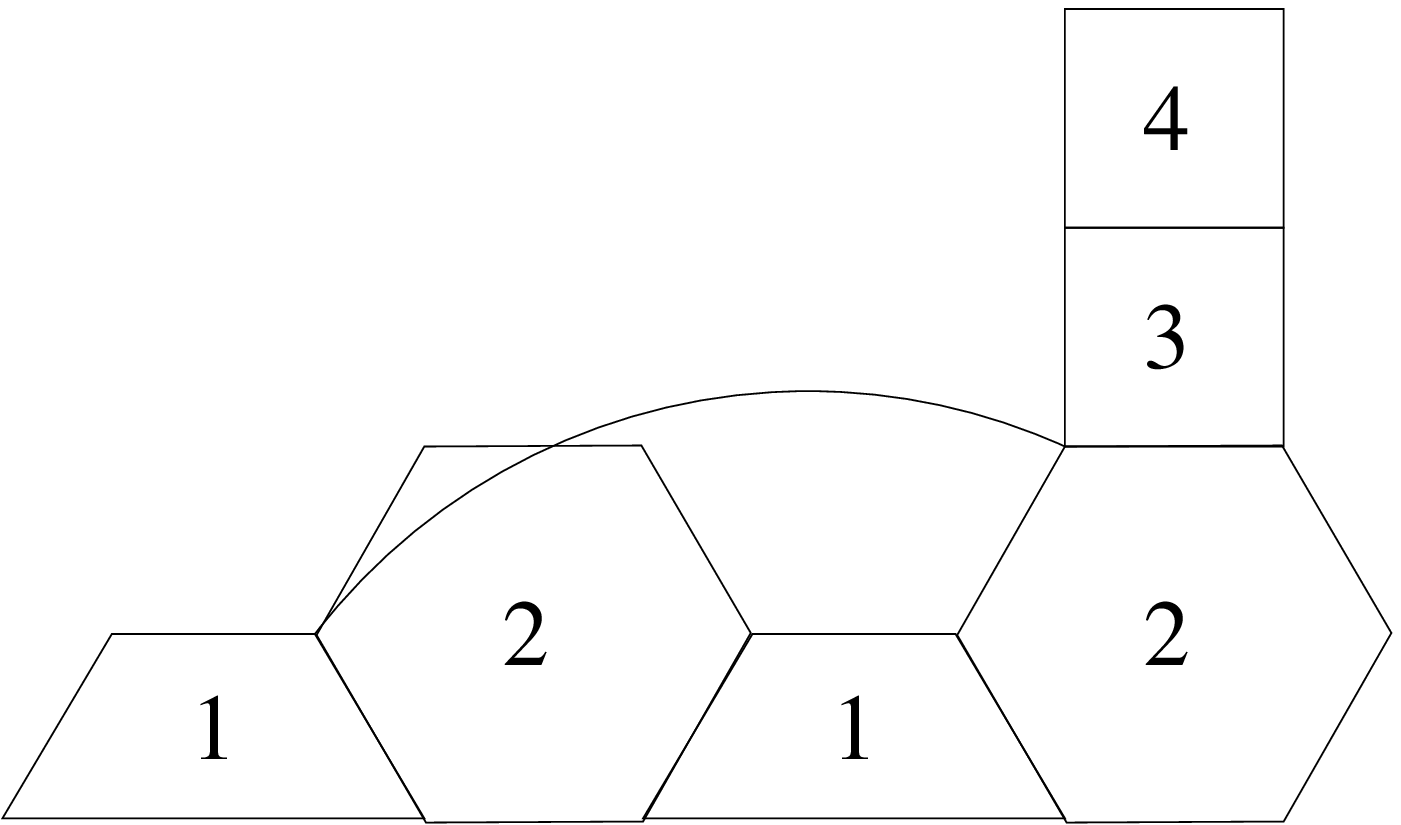} &~~~& \includegraphics[width = 0.2in , height = 0.5in]{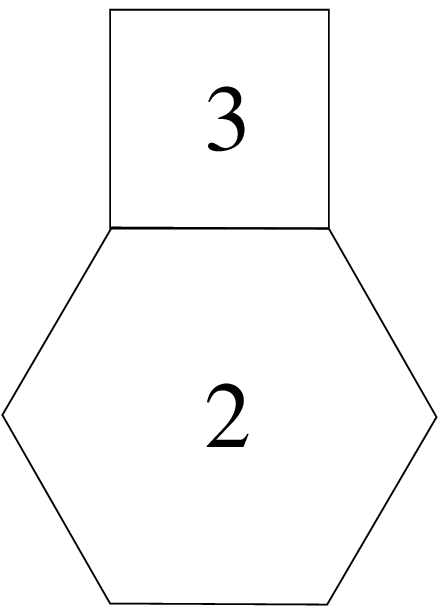}  \\
\includegraphics[width = 0.4in , height = 0.4in]{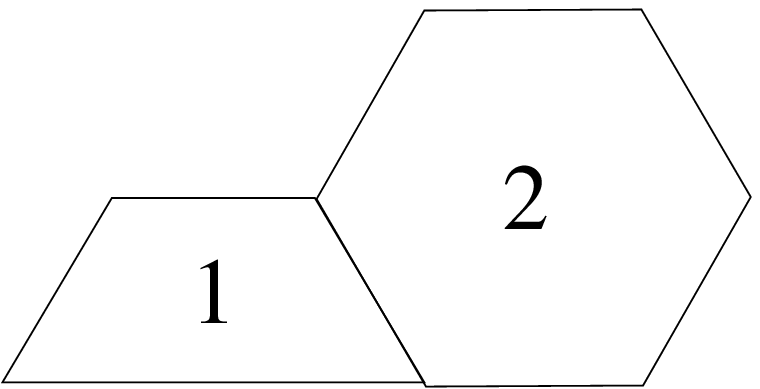}  &~~~& \includegraphics[width = 0.2in , height = 0.6in]{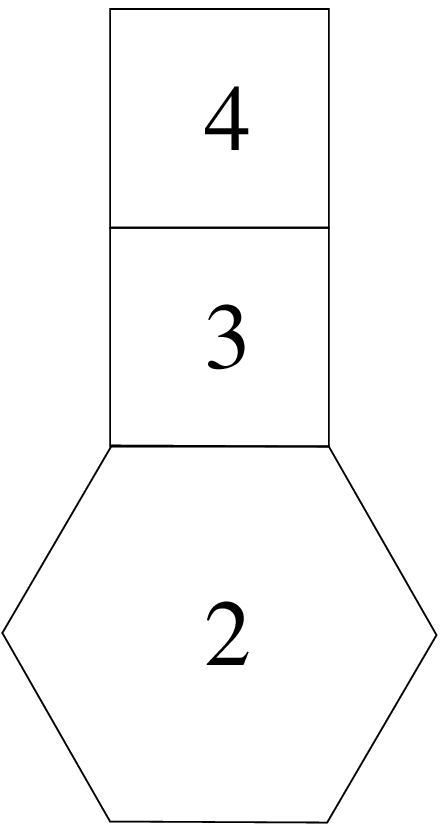} &~~~ \\
~~~& \includegraphics[width = 0.25in , height = 0.4in]{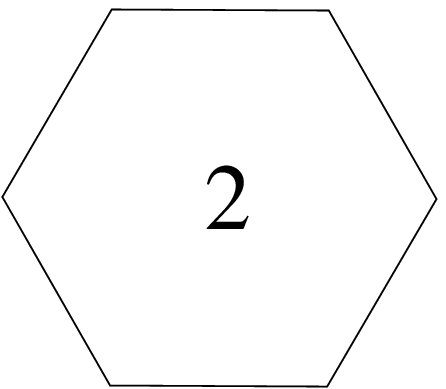} &~~~& \includegraphics[width = 0.2in , height = 0.2in]{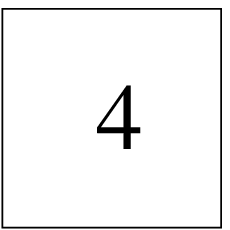} \\
\end{array}$ \\
The collection $\mathcal{G}_{B_4}$ (edge weights described in Section $4$). 
\end{center}\vspace{1em}

We shall use the convention from \cite{ClustII}, so that the initial exchange
matrix $B =||b_{ij}||_{i,j=1}^n$ contains rows of like sign.  Any rank $n$ cluster algebra of finite type has such a seed consisting
of a cluster of initial variables $\{x_1,\dots, x_n\}$ and a set
of $n$ binomial exchange relations of the form
$$x_ix_i^\prime = \prod_{j =1}^nx_j^{|b_{ij}|} + 1.$$
After mutating in the $k$th direction, i.e. applying an exchange relation of the form $x_k x_k^\prime =$ binomial, we obtain a new seed with cluster $\{x_1,x_2,\dots, x_n\}\cup \{x_k^\prime\}\setminus \{x_k\}$ and exchange matrix $B^\prime  = ||b_{ij}^\prime||_{i,j=1}^n$ such that the $b_{ij}^\prime$'s satisfy
$$b_{ij} = \begin{cases} -b_{ij} &\If i=k \mathrm{~or~}j=k, \\ 
b_{ij} + \max(-b_{ik},0)\cdot b_{kj} + b_{ik}\cdot \max(b_{kj},0) &\mathrm{~~otherwise}.\end{cases}$$
As we mention below in Remark \ref{sameseed}, we shall use an ordering of mutations in this paper so that we need only work with binomial exchanges of the form $x_k x_k^\prime =~($Monomial $+1)$.  Note that we shall use the notation $P_\alpha(x_1,x_2,\dots, x_n)$ to denote the numerator of the cluster variable with denominator $x_1^{\alpha_1}x_2^{\alpha_2}\cdots x_n^{\alpha_n}$ despite its similarity with the notation of $P(G)$ for the matching polynomial of graph $G$.

\vspace{1em}
The outline of the paper is as follows.  We proceed to prove Theorem 
\ref{vargraph} separately for the four families of non-exceptional type, 
starting with the well-studied case of $A_n$.  We will use different language than in \cite{CPpre}, \cite{ClustII}, or \cite{Schiffler}, and we include our own proof of this case to familiarize the reader with the techniques which we will utilize later in the paper.  Since the type of the cluster algebra will frequently be clear from context, we will simply denote tiles as $T_i$ or graphs as $G_\alpha$ (instead of $\Galph$).  We end with some comments and directions for further research.

\begin{Rem}
In \cite{YSys}, Fomin and Zelevinsky explicitly constructed Fibonacci polynomials for types $A_n$ and $D_n$, which provide an alternate combinatorial expansion formula for cluster varibles.  Generalizations of these polynomials, for other types, are defined in \cite{ClusIV}, where they are referred to as F-polynomials.
\end{Rem}

\section{$A_n$} \label{an}

The work in this section was done independently of the work of Carroll-Price \cite{CPpre} and the work of Schiffler \cite{Schiffler} mentioned in the introduction.  We will use the notation and the techniques of this section later in the paper for the $B_n$ and $D_n$ cases.  Thus we include this section even though the combinatorial interpretation given by Proposition \ref{CaseAn} is not new in this case, although we believe our proof via excision, as described by Lemma \ref{exciseGr}, is new.  This excision technique will also be utilized for the $B_n$ and $D_n$ cases in Section $4$.

We begin by reviewing the necessary characteristics of the cluster algebra 
of type $A_n$.  Recall that Lie algebra 
$A_n$ has a Dynkin diagram consisting of a line of $n$ vertices
connected by edges of weight one.

$$\bullet\line(1,0){3}\bullet\line(1,0){3}\bullet\line(1,0){3}\bullet\line(1,0){3}\bullet\line(1,0){3}
 \bullet\line(1,0){3} \dots \dots \line(1,0){3} \bullet$$

\noindent Thus the associated Cartan matrix has the form

$$\begin{bmatrix}
 2  & -1 &  0  &  0  & \dots & 0  & 0 \\
-1  &  2 & -1  &  0  & \dots & 0  & 0 \\
 0  & -1 &  2  & -1  & \dots & 0  & 0 \\
 0  &  0 & -1  &  2  & \dots & 0  & 0 \\
\dots & \dots & \dots & \dots & \dots & \dots & \dots \\
0  &  0 & 0  & 0  & \dots & -1 & 2 \\
\end{bmatrix},$$

\noindent and thus using the convention given in \cite{ClustII}
the associated exchange matrix is

$$B^{A_n} = ||b_{ij}|| = \begin{bmatrix}
 0  & 1 &  0  &  0  & \dots & 0  & 0 \\
-1  &  0 & -1  &  0  & \dots & 0  & 0 \\
 0  & 1 &  0  & 1  & \dots & 0  & 0 \\
 0  &  0 & -1  &  0  & \dots & 0  & 0 \\
\dots & \dots & \dots & \dots & \dots & \dots & \dots \\
0  &  0 & 0  & 0  & \dots & (-1)^{n+1} & 0 \\
\end{bmatrix}.$$

\vspace{1em} \noindent Notice that every row has like sign and that the matrix
is skew-symmetrizable (and in fact skew-symmetric in this case).
 The bipartite seed for a cluster algebra of type $A_n$ therefore
consists of a initial cluster of variables $\{x_1,x_2,\dots,
x_n\}$ and exchange matrix $B^{A_n}$ which encodes the following
exchange binomials

\begin{eqnarray*}
x_1x_1^\prime &=& x_2 + 1 \\
x_2x_2^\prime &=& x_1x_2 + 1 \\
x_3x_3^\prime &=& x_2x_4 + 1 \\
\dots \\
x_{n-1}x_{n-1}^\prime &=& x_{n-2}x_n + 1 \\
x_nx_n^\prime &=& x_{n-1} + 1.
\end{eqnarray*}

We describe a set of tiles from which we will build our family of graphs.  In the case of $A_n$, let tiles $T_1,\dots, T_n$ be squares defined as follows: 

\begin{Def} Tile
$T_1$'s northern edge is given weight $x_2$ while the other three
are given weight $1$.  Tile $T_n$'s southern edge is weighted
with value $x_{n-1}$ and the rest are weighted with value $1$.
Finally all other $T_i$ have a weight of $x_{i+1}$ given to their
northern edge, $x_{i-1}$ for their southern edge while the eastern
and western edges are given weight $1$.
\end{Def}

\begin{center} \includegraphics[width = 3in , height
= 0.6in]{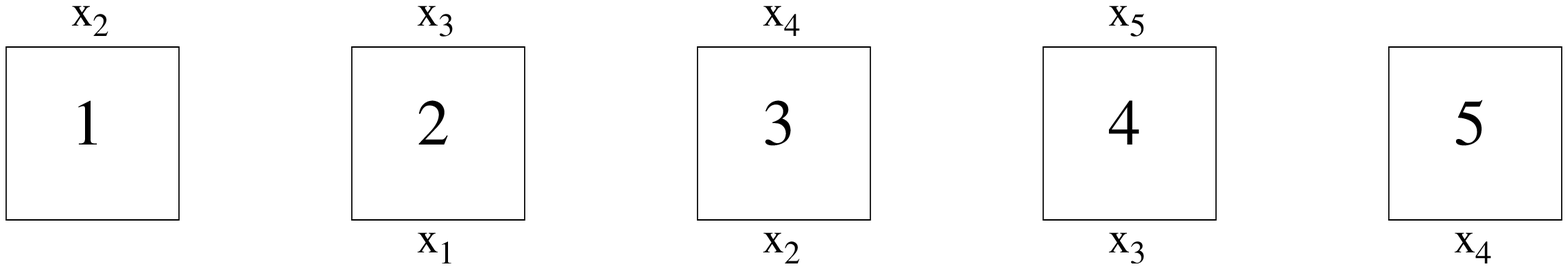}\\ The tiles for cluster algebra of type $A_5$.
\end{center}\vspace{1em} 

\noindent Let $\mathcal{G}_{A_n}$ be the set of graphs that can be
built from these $n$ tiles given the following gluing rule.

\begin{Rule} \label{Gluing} Without allowing reflections or rotations of the tiles, tile $T_i$ can be glued to tile $T_j$ if and only if the identified edge (as
an edge of $T_i$) lies clockwise from an edge weighted $x_j$ and
clockwise from an edge weighted $x_i$ (as an edge of $T_j$).
\end{Rule}

\noindent Since tile $T_i$ only contains edges of weight $x_{i+1}$
and $x_{i-1}$, and these weights appear across from each other,
this rule uniquely describes how the blocks can connect.

\begin{Lem}
Given the above tiles, $\mathcal{T}_{A_n}$, and the above gluing rule, the collection of possible graphs is enumerated by the set of subsets
$$\{T_i,T_{i+1},\dots,
T_{j-1},T_j\}$$ for $1 \leq i < j \leq n$. \end{Lem}

\noindent  This collection $\mathcal{G}_{A_n}$ has the same
cardinality as the set of positive roots of the Lie algebra of
type $A_n$ using the bijection $$
T_i \cup T_{i+1} \cup T_{i+2} \cup \dots \cup  T_{j-1} \cup T_j
\rightarrow \alpha_i + \dots + \alpha_j.$$  As shown
in \cite{ClustII}, this implies that the cardinality is also the
same as the number of non-initial cluster variables for the
bipartite cluster algebra of type $A_n$.

\begin{Prop} \label{CaseAn} The set of graphs $\mathcal{G}_{A_n}$ is in bijection
with the set of non-initial cluster variables for a coefficient-free cluster algebra of type $A_n$ and satisfy the statement of Theorem \ref{vargraph}.
\end{Prop}

\vspace{1em} \begin{center} \includegraphics[width = 3in ,
height = 3in]{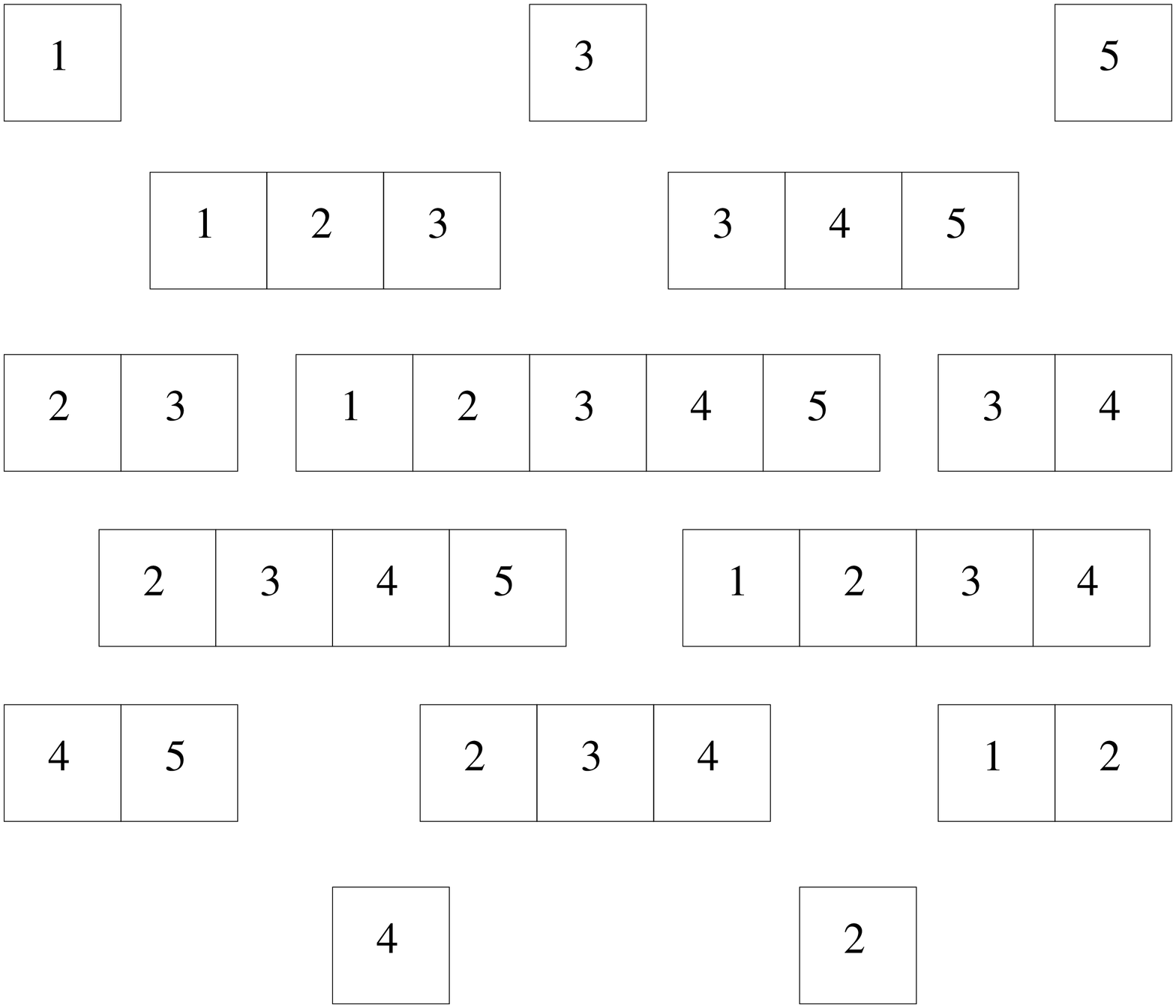} \\  The collection $\mathcal{G}_{A_5}$. \end{center} \vspace{1em}

To prove this proposition we will take a detour through a case
we refer to as $A_{\infty}$. In this case, our set of tiles is in bijection with the integers, and we define $T_i$ to have $x_{i+1}$ on its northern edge, $x_{i-1}$ on its southern edge for all $i \in \z$.  Without the issue of boundaries, it is easier to show that a certain lattice of graphs corresponds to the non-initial cluster
variables.  After doing so, we choose a periodic specialization for the initial variables to recover a corresponding region of this lattice for any specific $A_n$.  We start with the following observation.

\begin{Rem} \label{sameseed}For all $n$, if we start with the above exchange matrix $B^{A_n}$ and apply the binomial exchanges corresponding to relations $1$, $3$, $5, \dots n$ (resp. $n-1$) if $n$ is odd (resp. even) the resulting
exchange matrix is $-B$.  Afterward, applying the relations $2$, $4$,
$6, \dots n-1$ (resp. $n$) if $n$ is odd (resp. even) to exchange
matrix $-B^{A_n}$ results in the intial exchange matrix $B^{A_n}$.  In fact, in both of these cases, the order of the exchanges does not matter, and the intermediate exchange matrices will have rows of like sign for all relevant $x_k$ not already exchanged.  By the definition of matrix mutation, this procedure will in
fact work for any cluster algebra where the seed has an exchange
matrix that is tri-diagonal ($b_{ij} = 0$ if $|i-j| \not = 1$).
 Thus we can calculate a row of cluster variables at a time by
applying the exchange relations relative to the two previous rows.
 The tri-diagonal condition includes the cases $A_n$, $B_n$, $C_n$,
and $G_2$ and minor modifications to the procedure will allow it
to work for $D_n$.   
\end{Rem}

Returning to the $A_n$ case, after applying exchange $1,3,5,\dots, n$ (resp. $n-1$) we
have cluster

$$\{x_1^{(1)},~x_2,~x_3^{(1)},~x_4,~x_5^{(1)},\dots, ~x_{n-1},~x_n^{(1)} \} $$
$$(\mathrm{resp.~~} \{x_1^{(1)},~x_2,~x_3^{(1)},~x_4,~x_5^{(1)},\dots, ~x_{n-2},~x_{n-1}^{(1)},~x_n
\}~)$$\vspace{1em} 

\noindent where $x_i^{(1)} = \frac{x_{i-1}x_{i+1}+1  }{  x_i}$ using
the convention $x_0 = x_{n+1} = 1$.  Analogously, applying
exchanges $2,4,6, \dots, n-1$ (resp. $n$) we obtain the cluster

$$\{x_1^{(1)},~x_2^{(2)},~x_3^{(1)},~x_4^{(2)},~x_5^{(1)},\dots, ~x_{n-1}^{(2)},~x_n^{(1)} \} $$
$$(\mathrm{resp.~~~~} \{x_1^{(1)},~x_2^{(2)},~x_3^{(1)},~x_4^{(2)},~x_5^{(1)},\dots, ~x_{n-2}^{(2)},~x_{n-1}^{(1)},~x_n^{(2)}
\}~)$$\vspace{1em} 

\noindent where $x_i^{(2)} =
\frac{x_{i-2}x_{i+2}+x_{i-2}+x_{i+2}+x_{i-1}x_{i+1}+1  }{ 
x_{i-1}x_ix_{i+1}}$ for $1\leq i \leq n$, if we set $x_{-1}=x_{n+2}=0$.  By Remark
\ref{sameseed} we can make a lattice of cluster variables by
applying exchanges iteratively (row-by-row) in this order.

\vspace{1em} $\begin{array}{cccccccccccc}
x_1 &~~~& x_3 &~~~& x_5 &~~~& \dots &~~~& x_{n-2} &~~~& x_{n} \\
 ~~~& x_2 &~~~& x_4 &~~~& x_6   &~~~& \dots  &~~~& x_{n-1} & ~~~ \\
x_1^{(1)} &~~~& x_3^{(1)} &~~~& x_5^{(1)} &~~~& \dots &~~~& x_{n-2}^{(1)} &~~~& x_{n}^{(1)} \\
 ~~~& x_2^{(2)} &~~~& x_4^{(2)} &~~~& x_6^{(2)}   &~~~& \dots  &~~~& x_{n-1}^{(2)} & ~~~ \\
x_1^{(3)} &~~~& x_3^{(3)} &~~~& x_5^{(3)} &~~~& \dots &~~~& x_{n-2}^{(3)} &~~~& x_{n}^{(3)} \\
 ~~~& x_2^{(4)} &~~~& x_4^{(4)} &~~~& x_6^{(4)}   &~~~& \dots  &~~~& x_{n-1}^{(4)} & ~~~
\end{array}$ \\ \begin{center} Top six rows of this lattice (assuming $n$ odd).
\end{center}\vspace{1em} 

By the binomial exchange relations, this lattice satisfies the
diamond condition which states that the relation $ad = bc + 1$
holds for any four elements arranged as a diamond.

\begin{center}$\begin{array}{ccc}
&a& \\
b&&c \\
&d&
\end{array}$\end{center}

\noindent For example $x_ix_i^{(1)} = x_{i-1}x_{i+1}+1$ for $i \in
\{2,3,\dots, n-1\}$.  Furthermore, this lattice can be extended
periodically by using the conventions $x_{-1}=x_{n+2}=0$, $x_0 =
x_{n+1} = 1$, and extending further using negatively weighted
variables.  Note that the negatives are necessary since we wish the configurations 
\begin{center}$\begin{array}{ccccccc}
&0&  &~& &0& \\
b&&1 &~&1 && c \\
&0&  &~& &0&
\end{array}$\end{center}
to satisfy the diamond condtion.

\vspace{1em} $\begin{array}{ccccccccccccccc}
 -x_3^{(2)} &~~~& -x_1^{(2)} &~~~& 0 &~~~& x_1^{(2)} & \dots & x_{n}^{(2)} &~~~& 0 &~~~& -x_{n}^{(2)} \\
~& -x_2^{(1)} &~~~& -1 &~~~& 1 & ~~~& x_2^{(1)} & \dots & 1 &~~~ & -1 &~~~  \\
 -x_3 &~~~& -x_1 &~~~& 0 &~~~& x_1 & \dots & x_{n} &~~~& 0 &~~~& -x_{n}  \\
~& -x_2 &~~~& -1 &~~~& 1 & ~~~& x_2 & \dots  &  1 &~~~ & -1 &~~~ \\
 -x_3^{(1)} &~~~& -x_1^{(1)} &~~~& 0 &~~~& x_1^{(1)} & \dots & x_{n}^{(1)} &~~~& 0 &~~~& -x_{n}^{(1)} \\
~& -x_2^{(2)} &~~~& -1 &~~~& 1 & ~~~& x_2^{(2)} & \dots  & 1 &~~~
& -1 &~~~
\end{array}$ \\ \begin{center} Six rows of extended lattice (assuming $n$ odd).
\end{center}\vspace{1em} 

This lattice can also continue infinitely in the vertical direction, as well as
horizontally, extending vertically in the unique way
that preserves the diamond condition throughout the entire
lattice. Consequently all $A_n$ can be treated simultaneously by
considering the infinite diamond pattern ($A_{\infty}$) which
starts with sequence $\{\dots, -y_2,-y_1,y_0,y_1,y_2,\dots\}$
zig-zagging to create the initial two rows. To obtain the extended
$A_n$-lattice for a specific $n$ we let
\begin{eqnarray*}
y_1 &=& 1 \\
y_i &=& x_{i-1} \mathrm{~for~} i \in \{2,\dots, n+1\} \\
y_{n+2} &=& 1 \\
y_{n+3} &=& y_0 \\
y_{n+3+k} &=& -y_{n+3-k} \mathrm{~for~} k \in \{0,\dots, n+3\} \\
y_{2n+6+k} &=& y_{k}
\end{eqnarray*}
\noindent and then take the limit as $y_0$ goes to zero.  We do not set $y_0$ and $y_{n+3}$ to be zero directly since this would sometimes result in indeterminate expressions of the form ``$0/0$''.  Also we use the shifted indices here since it will make the ensuing arguments more symmetrical.

As a consequence of these substitutions, it suffices to start the proof of Proposition \ref{CaseAn} by proving the combinatorial
interpretation for the infinite diamond pattern corresponding to
$A_{\infty}$, which we write in terms of $y_i$'s for $i\in \z$.  
Even though we now have a boundary-less lattice, every given $y_i^{(j)}$ can be computed locally by considering the necessary mutations stemming from a finite half diamond extending back to the initial two rows of $y_i$'s.

\begin{center} \includegraphics[width = 5in ,
height = 0.5in]{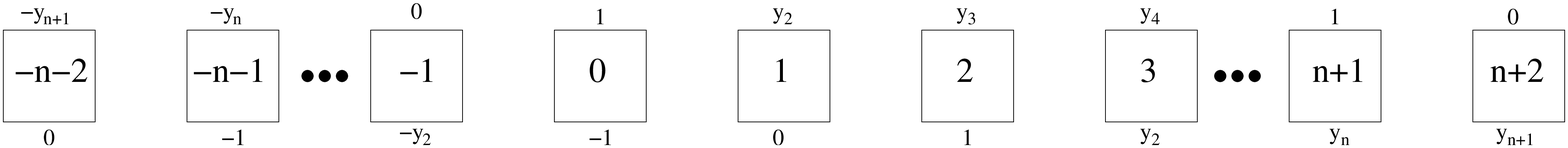}\\ The tiles for the extended
$A_n$-lattice. \end{center}\vspace{1em} 

For the purposes of the $A_\infty$ case, for all $i\in\z$, we let $\tilde{T_i}$ denote the tile with $y_{i+1}$ on its northern edge and $y_{i-1}$ on its southern edge.  We will utilize variables $y_i$'s and tiles $\tilde{T_i}$'s until we do the final substitution at the end of the proof of Proposition \ref{CaseAn}.  We now proceed to prove the combinatorial interpretation in this boundary-less version.  We start with the base case where we can easily see that the combinatorial interpretation works for cluster variables with denominator $y_i$.  To see this, we observe that $y_iy_i^{(1)} = y_{i-1}y_{i+1}+1$ corresponds to the two perfect matchings of the graph consisting of tile $\tilde{T}_i$ by itself.  Similarly, we observe the
bijection works for the second row of non-initial cluster
variables by the definition of $y_i^{(2)}$.  We see this by verifying that the
graph containing tiles $\tilde{T}_{i-1}$, $\tilde{T}_{i}$, $\tilde{T}_{i+1}$ connected in that order bijects to cluster variable $y_i^{(2)}$, i.e. $y_i^{(2)}y_{i-1}y_iy_{i+1} = P(\tilde{T}_{i-1}\cup \tilde{T}_i \cup \tilde{T}_{i+1})$.  By a technique of graphical condensation developed by Eric Kuo \cite{Kuo}, we obtain the following combinatorial interpretation for the rest of the rows.

\begin{Lem} \label{diag} Cluster variable $y_i^{(j)}$ bijects to graph 
$$\tilde{T}_{i-j+1}\cup \dots \cup \tilde{T}_{i+j-1},$$
i.e. the grid graph containing exactly $2j-1$ tiles, namely tiles $\tilde{T}_{i-j+1}$ through $\tilde{T}_{i+j-1}$ connected in order. \end{Lem}

\begin{proof} The proof follows from a slight variant of the argument given in \cite{MusPropp}.  Here we
need to be more careful with the labeling scheme, but the same pairings will yield the
desired result.  In fact if one lets $y_i = x$ if $i$ even and $y_i =y$ if $i$ odd, one recovers
the $A(2,2)$ case analyzed in \cite{MusPropp}.
In particular, we inductively assume for all $i\in\z$ that $y_i^{(j-1)}$ bijects to graph $G_1^i = \tilde{T}_{i-j+2}\cup\dots\cup\tilde{T}_{i+j-2}$ and 
$y_i^{(j-2)}$ bijects to graph $G_2^i = \tilde{T}_{i-j+3}\cup\dots\cup\tilde{T}_{i+j-3}$, in the sense that 
$y_i^{(j-1)} = \frac{P(G_1^i)}{y_{i-j+2}\cdots y_{i+j-2}}$ and 
$y_i^{(j-2)} = \frac{P(G_2^i)}{y_{i-j+3}\cdots y_{i+j-3}}$.  It thus suffices to show, for all $i\in\z$, that Laurent polynomials $y_i^{(j)}$, defined as 
$\frac{y_{i-1}^{(j-1)}y_{i+1}^{(j-1)}+1}{y_i^{(j-2)}}$, equal $\frac{P(G_0^i)}{y_{i-j+1}\cdots y_{i+j-1}}$ where $G_0^i = \tilde{T}_{i-j+1}\cup\dots\cup\tilde{T}_{i+j-1}.$  
We use our induction hypothesis and normalize to rewrite our desired equation as 
{\small
\begin{eqnarray} \label{norma} 
P(G_0^i)P(G_2^i) &=& P(G_1^{i-1})P(G_1^{i+1})+y_{i-j+1}y_{i-j+2}y_{i-j+3}^2\cdots y_{i+j-3}^2y_{i+j-2}y_{i+j-1}.
\end{eqnarray}
}

One can decompose graph $G_0^i$ into a superposition of graphs $G_1^{i-1}\cup G_1^{i+1}$ so that $G_2^i$ is the intersection of overlap.  Out of the two subgraphs, only $G_1^{i-1}$ contains tiles $\tilde{T}_{i-j+1},~\tilde{T}_{i-j+2}$ and only $G_1^{i+1}$ contains $\tilde{T}_{i+j-1}, ~\tilde{T}_{i+j-2}$.  Let $M(G)$ denote the set of perfect matchings of graph $G$, ${m_0}^\prime$ denote the matching of $G_0^i$ using the horizontal edges of $\tilde{T}_{i-j+2},~\tilde{T}_{i-j+4},\dots,~ \tilde{T}_{i+j-4}, ~\tilde{T}_{i+j-2}$, and ${m_2}^\prime$ denote the matching of $G_2^i$ using the horizontal edges of $\tilde{T}_{i-j+3},~\tilde{T}_{i-j+5},\dots,~ \tilde{T}_{i+j-5}, ~\tilde{T}_{i+j-3}$.  The pair of matchings $(m_0^\prime, m_2^\prime)$ has exactly the weight of the excess monomial 
$$y_{i-j+1}y_{i-j+2}y_{i-j+3}^2\cdots y_{i+j-3}^2y_{i+j-2}y_{i+j-1}.$$  We finish the proof of Lemma \ref{diag} by exhibiting a weight-preserving bijection between
$M(G_0^i)\times M(G_2^i)\setminus\{ ( {m_0}^\prime, {m_2}^\prime )$ and $M(G_1^{i-1})\times M(G_1^{i+1})$, thus showing (\ref{norma}).

We define our bijection piece-meal on $M(G_0^i) \times M(G_2^i) \setminus \{ ( {m_0}^\prime, {m_2}^\prime )$, first considering the case where the horizontal edges of penultimate tile $\TT_{i+j-2}$ in $G_0^i$ are not used.  In this case, the pair of matchings from 
$M(G_0^i)\times M(G_2^i)$ reduces to a pair from $M(G_1^{i-1}\cup \TT_{i+j-1})\times M(G_2^i)  $.  We define $\phi(m_0,m_2)=(m_{-1},m_1)$ for such matchings by letting $m_{-1}$ be the corresponding matching of $G_1^{i-1}$ and build matching $m_1$ by adjoining the matching of $G_2^i$ to the matching of $\TT_{i+j-1}$.  In other words, map $\phi$ takes tile $\TT_{i+j-1}$ and slides it down from $G_0^i$ onto $G_2^i$ to obtain $G_1^{i+1}$ with the matching included.  This also leaves $G_1^{i-1}$ in place of $G_0^i$.  

If on the other hand, the horizontals of $\TT_{i+j-2}$ \emph{are used}, then the situation is more complicated.  If we restrict further to the case where the rightmost vertical edge is used in $G_2^i$, we can define $\phi$ analogously by sliding down tiles $\TT_{i+j-2}$ and $\TT_{i+j-1}$.  We are also forced to use the rightmost vertical edge of $G_1^{i-1}$ in this case.

We can continue defining $\phi$ iteratively, defining it for classes characterized by the length of the pattern of horizontals on the right-hand sides of $G_0^i$ and $G_2^i$.  If the horizontals of $\TT_{i+j-2},~  \TT_{i+j-4}, \dots, \TT_{i+j-2\ell}$ in $G_0^i$, the horizontals of $\TT_{i+j-3},~\TT_{i+j-5},\dots,$ $\TT_{i+j-2\ell - 1}$ (resp. $\TT_{i+j-2\ell + 1}$)  in $G_2^i$ are used, accompanied by a vertical edge between tiles $\TT_{i+j-2\ell-2}$ and $\TT_{i+j-2\ell-1}$ of $G_0^i$ (resp. $\TT_{i+j-2\ell-1}$ and $\TT_{i+j-2\ell}$ of $G_2^i$), then $\phi$ swaps the right-hand sides of these two graphs, leaving the left-hand sides alone up until tile $\TT_{i+j-2\ell-3}$ (resp. $\TT_{i+j-2\ell-2}$).  This construction makes sense as long as we eventually encounter a vertical edge as we move leftward via the use of these horizontal edges since given such patterns, neither the matching of $G_0^i$ nor of $G_2^i$, will use the horizontal edges of tile $\TT_{i+j-2\ell-2}$ (resp. $\TT_{i+j-2\ell-1}$).

Map $\phi$ is injective since the inverse map just swaps back the right-hand sides as dictated by the alternating pattern of horizontals.  Since we have exhaustively enumerated the pairs of matchings $(m_{-1},m_1)$ by splitting into classes according to the longest alternating pattern of horizontals, we also have surjectivity.  Lastly, it is easy to verify that $(m_0^\prime,m_2^\prime)$ is the unique matching which cannot be decomposed into a pair $(m_{-1},m_1)$.
\end{proof}

Analogous pairings will also appear below in the
arguments for the case of $B_n$.  Notice that by Lemma \ref{diag}, the diagonals of the
lattice satisfy the following two properties:

\vspace{1em} $\bullet$ On any of the diagonals travelling from \emph{SW to NE}, all graphs
\emph{end} with the same tile.

$\bullet$ On any of the diagonals travelling from \emph{NW to SE}, all graphs \emph{start}
with the same tile. \vspace{1em}

\noindent We now wish to show how to specialize to the case of a specific $A_n$ by imposing periodicty and boundary conditions on the intial two rows of variables.  To accomplish this goal, we must (a) verify that the $y_i^{(j)}$'s satisfy the correct horizontal periodicity of the extended-$A_n$ lattice once we apply the proper substitutions of variables, (b) verify that we have a vertical periodicity as well and really only need to worry about a finite collection of graphs which we readily identify as the set $\mathcal{G}_{A_n}$.
To show (a), it suffices to show $y_{-i}^{(j)}=-y_{i}^{(j)}$, $y_0^{(j)}=0=y_{n+3}^{(j)}$, and $y_1^{(j)}=1=y_{n+2}^{(j)}$ for all $j$.  The diamond condition will then induce the horizontal periodicity for non-initial cluster variables.

Notice that $y_{n+3+k} = -y_{n+3-k}$ and the periodicity $y_{2n+6+k}=y_k$ imply the relation $y_k = y_{-k}$ for all $k\in \z$.
Thus Lemma \ref{diag}, together with $y_{-i}=-y_i$, imply the relations 
\begin{eqnarray*}
y_i^{(j)} &=& \frac{P(\TT_{i-j+1}\cup \dots \cup \TT_{i+j-1})}{y_{i-j+1}\cdots y_{i+j-1}} \mathrm{~~~and}\\
y_{-i}^{(j)} &=& \frac{P(\TT_{-i-j+1}\cup \dots \cup \TT_{-i+j-1})}{y_{-i-j+1}\cdots y_{-i+j-1}} = 
\frac{P(\TT_{-(i-j+1)}\cup \dots \cup \TT_{-(i+j-1)})}{(-1)^{2j-1}y_{i-j+1}\cdots y_{i+j-1}},
\end{eqnarray*}
where the second equality comes from reversing the order of the tiles.

Furthermore, tile $\TT_{-k}$ has $y_{-k+1}$ on its northern edge and $y_{-k-1}$ on its southern edge, while $\TT_{k}$ has $y_{k+1}$ on its northern edge and $y_{k-1}$ on its southern edge.  Thus the relation $y_{-i}=-y_i$ for all $i\in\z$ induces $y_{-i}^{(j)}=-y_i^{(j)}$ for all $i\in \z$ and $j\geq 1$.

As a corollary, we obtain $y_0^{(j)}=0$ for all $j$.  Under Lemma \ref{diag}, the graph corresponding to $y_0^{(j)}$ is centered around tile $\TT_0$, i.e. $$\TT_{-j+1}\cup \dots \cup \TT_0 \cup \dots \cup \TT_{j-1}.$$  Laurent polynomial $y_{n+3}^{(j)}$ analogously corresponds to a graph centered around tile $\TT_{n+3}$, which is equivalent to tile $\TT_0$ after applying the above periodicity and substution in the extended $A_n$-lattice.  Consequently, we induce that for all $j\geq 1$, $y_{n+3}^{(j)}= 0$ similarly.  
To complete our proof of (a), we prove the following Lemma.

\begin{Lem} \label{CenterOne}
If $H_j$ signifies the graph $\TT_{-j}\cup \dots \cup \TT_{j+1}$, i.e. a graph with an even number of tiles, with leftmost tile $\TT_{-j}$, and centered around subgraph $\TT_0\cup \TT_1$, then (using $y_{-i}=-y_i$), after dividing the matching polynomial $P(H_j)$ by the proper monomial, we find that graph $H_j$ bijects to Laurent polynomial $y_{j+2}$.

Moreover, any graph $\hat{H}_j$ with an odd number of tiles centered around tile $\TT_1$ (resp. $\TT_{n+2}$) bijects to $1$ as a Laurent polynomial, regardless of the number of tiles.
\end{Lem}

\begin{proof}
We start by proving the result for graphs with an even number of tiles, beginning with two base cases.  The graph $\TT_0\cup \TT_1$ has three perfect matchings, but due to signs, two of them cancel with one another.  We thus obtain $P(\TT_0\cup \TT_1)= y_0y_2$.  After dividing through by $y_0y_1=y_0(1)$ we get Laurent polynomial $y_2$.

Secondly, the graph $\TT_{-1}\cup\TT_0 \cup \TT_1 \cup \TT_2$ has eight perfect matchings, but one of them has a weight containing submonomial $y_0^2$, and six of the other seven cancel with each other after letting $y_{-i}=-y_i$ and $y_1=1$.  We are left with one perfect matching, which has weight $-y_0y_2y_3$.  This time the denominator is $y_{-1}y_0y_1y_2=(-1)y_0(1)y_2$ and we get $y_3$ after division.

We now assume that $j\geq 2$ and wish to show $$P(H_j)= y_{-j}y_{1-j}y_{2-j}\cdots y_{-2}(-1)y_0(1)y_2\cdots y_{j}y_{j+1}y_{j+2}$$ for all such $j$.  Notice that if the rightmost vertical edge of graph $H_{j}$ is used in a matching, our computation of $P(H_j)$ reduces to the computation of $P(H_j^{\prime})$, where $H_j^{\prime} =\TT_{-j}\cup \dots \TT_0\cup \dots \cup \TT_{j}$, a graph centered around tile $\TT_0$.  However, we know from previous arguments that such a graph corresponds to zero as a Laurent polynomial.  Thus any matching with the rightmost vertical edge of graph $H_j^\prime$ does not contribute to the enumerator $P(H_j)$.  Thus we must use the two rightmost horizontal edges, which have weight $y_{j}y_{j+2}$, and then compute $P(\overline{H_{j-1}})$ 
where $\overline{H_{j-1}}=\TT_{-j}\cup \dots \cup \TT_{j-1}$.  However, such a graph is the horizontal reflection of graph $H_{j-1}$ and so by analogous logic, we are now forced to use the two \emph{leftmost} horizontal edges in our matchings to get a nontrivial contribution.  Such edges have weight $y_{-j-1}y_{-j+1} = (-y_{j+1})(-y_{j-1}) = y_{j-1}y_{j+1}$.  Consequently, after two iterations, and two helpfully placed negative signs, we have the identity 
$$P(H_j) = (-y_{-j})(-y_{1-j}) P(H_{j-2})y_{j+1}y_{j+2}.$$ Induction thus yields the result for the case of $H_j$, i.e. a graph with an even number of tiles centered around $\TT_0\cup \TT_1$.

We now prove this result for the corresponding graphs with an odd number of tiles.  We only describe the proof for the case of those centered around $\TT_1$ since the proof for those centered around $\TT_{n+2}$ is analogous, only with messier notation.  
We let $\hat{H}_j= \TT_{-j}\cup \cdots \cup \TT_{j+2}$ and
for the moment ignore boundaries and periodicity; only use the assignments $y_1=1$, $y_{-1}=-1$, and take the limit of $$\frac{P(\hat{H}_j)}{y_{-j}\cdots y_{-2}y_{-1}y_0y_1y_2\cdots y_{j+2}}$$ as $y_0\rightarrow 0$.

We look at possible perfect matchings of graph $\hat{H}_j$, and note that if both horizontal edges or both vertical edges of tile $\TT_j$ are used, then we are reduced to computing $P({\hat{H}_j}^\prime \sqcup \TT_{j+2}) =  P({\hat{H}_j}^\prime)P(\TT_{j+2})$ with ${\hat{H}_j}^\prime = \TT_{-j}\cup \cdots \cup \TT_{j}$.  However, since ${\hat{H}_j}^\prime$ is centered around tile $\TT_0$, the Laurent polynomial $P({\hat{H}_j}^\prime)/y_0$ tends to zero as $y_0\rightarrow 0$.  We conclude that any matching of $\hat{H}_j$ resulting in a nontrivial contribution to 
$\lim_{y_0\rightarrow 0}\frac{P(H)}{y_{-j}\cdots y_{-2}y_{-1}y_0y_1y_2\cdots y_{j+2}}$ must utilize the two horizontal edges of tile $\TT_{j+1}$ with weight $y_{j}y_{j+2}$.  

However, this step reduces our calculation to that of the matching polynomial of 
$\overline{H^{j-1}} = \TT_{-j} \cup \dots \cup \TT_{j-1}$, which has an even number of tiles and is centered around $T_{-1}\cup T_0$.  By our earlier logic, we therefore can only have a nontrivial contribution to $P(\hat{H}_j)$ in the case where we use the two leftmost horizontal edges, which have weight $y_{-j-1}y_{-j+1}$.  This reduction results in subgraph $H_{j-2}$, which is centered around $T_0\cup T_1$, and so by induction 
$$P(\hat{H}_j) = (-y_{-j})y_{1-j}\bigg(y_{2-j}\cdots y_{j}\bigg)(-y_{j+1})y_{j+2}$$ which is the same as the denominator corresponding to $\hat{H}_j= \TT_{-j}\cup \cdots \cup \TT_{j+2}$.
\end{proof}

We thus have shown $(a)$, i.e. that after applying substitutions to the initial row of the $A_\infty$ lattice, we get the extended $A_n$ lattice with the proper horizontal periodicity.  We now wish to show the vertical periodicity of $(b)$.

\begin{Lem} \label{exciseGr}
Whenever a graph contains the tile $\TT_1$ (resp. $\TT_{n+2}$), we may excise the graph by removing a subgraph centered around $\TT_1$ (resp. $\TT_{n+2}$).  without changing the corresponding Laurent polynomial.\end{Lem}

Underneath initial variables $y_2=x_1, \dots, y_{n+1}=x_n$, any graph containing tile $\TT_0$ (resp. $\TT_{n+3}$), or other forbidden tiles $\TT_i$ with $i$ outside the range $\{2,\dots, n+1\}$ must contain either tile $\TT_1$ or tile $\TT_{n+2}$.  Consequently, this process of excision is sufficient to eliminate all tiles $\TT_i$ with $i \not\in \{2,\dots, n+1\}$.  (We will see shortly that the apparent problem of a graph containing both $\TT_1$ and $\TT_{n+2}$ where the subgraphs of excision overlap is not actually an issue.)

\begin{proof}
Without loss of generality, assume that graph $G$ contains tile $\TT_1$ and has the form $G=\TT_{2-j}\cup \dots \cup \TT_1 \cup \dots \cup \TT_j \cup \TT_{j+1} \cup \dots \cup \TT_{j+k}$ for $j, k\geq 1$.  The content of the claim is that graph $G$ and graph $G^\prime = \TT_{j+1}\cup \dots \cup \TT_{j+k}$ biject to the same cluster variable.

We categorize perfect matchings of $G$ based on whether or not two horizontal edges appear on the tile $\TT_{j+1}$.  If they do not, then that matching of $G$ decomposes into a matching of subgraph $\TT_{2-j}\cup \dots \cup \TT_j$, and a matching of $\TT_{j+2}\cup \dots \cup \TT_{j+k}$.  Thus we obtain
\begin{eqnarray*} \frac{P(\TT_{2-j}\cup \dots \cup \TT_j)\cdot P(\TT_{j+2}\cup \dots \cup \TT_{j+k})}{y_{2-j}\cdots y_{j+k}} &=&
\frac{P(\TT_{2-j}\cup \dots \cup \TT_j)}{y_{2-j}\cdots y_{j}}\cdot \frac{P(\TT_{j+2}\cup \dots \cup \TT_{j+k})}{y_{j+1}\cdots y_{j+k}} \\ &=&
1\cdot \frac{P(\TT_{j+2}\cup \dots \cup \TT_{j+k})}{y_{j+1}\cdots y_{j+k}}
\end{eqnarray*}
as a contribution to the Laurent polynomial corresponding to graph $G$, where the second equality follows from Lemma \ref{CenterOne}.

If on the other hand, the horizontal edges of $\TT_{j+1}$ \emph{are used}, the matching decomposes differently and we get a contribution of 
$$\frac{P(\TT_{2-j}\cup \dots \cup \TT_{j-1})}{y_{2-j}\cdots y_{j-1}}\cdot
\frac{y_jy_{j+2}}{y_jy_{j+1}y_{j+2}}\cdot \frac{P(\TT_{j+3}\cup \dots \cup \TT_{j+k})}{y_{j+3}\cdots y_{j+k}},$$
which equals
$$y_j\cdot \frac{y_{j+2}}{y_{j+1}y_{j+2}}\cdot \frac{P(\TT_{j+3}\cup \dots \cup \TT_{j+k})}{y_{j+3}\cdots y_{j+k}}$$ by Lemma \ref{CenterOne}.  Comparing the sum of these two Laurent polynomials to the cluster variable corresponding to graph $G^\prime$ finishes the proof.
\end{proof}

With Lemma \ref{exciseGr} proved, we turn our attention back to the proof of Proposition
\ref{CaseAn}. 

\begin{proof}  [Proof of Proposition \ref{CaseAn}]
Recall that the extended $A_n$-lattice is horizontally periodic. We thus restrict
our attention to a region that lies between the two columns of positive ones and below the strip
of initial variables.

\vspace{1em} $\begin{array}{cccccccccccc}
x_1^{(1)} &~~~& x_3^{(1)} &~~~& x_5^{(1)} &~~~& x_7^{(1)} &~~~& x_{9}^{(1)} &~~~& x_{11}^{(1)} \\
 ~~~& x_2^{(2)} &~~~& x_4^{(2)} &~~~& x_6^{(2)}   &~~~& x_8^{(2)}  &~~~& x_{10}^{(2)} & ~~~ \\
x_1^{(3)} &~~~& x_3^{(3)} &~~~& x_5^{(3)} &~~~& x_7^{(3)} &~~~& x_{9}^{(3)} &~~~& x_{11}^{(3)} \\
 ~~~& x_2^{(4)} &~~~& x_4^{(4)} &~~~& x_6^{(4)}   &~~~& x_8^{(4)}  &~~~& x_{10}^{(4)} & ~~~ \\
x_1^{(5)} &~~~& x_3^{(5)} &~~~& x_5^{(5)} &~~~& x_7^{(5)} &~~~& x_{9}^{(5)} &~~~& x_{11}^{(5)} \\
 ~~~& x_2^{(6)} &~~~& x_4^{(6)} &~~~& x_6^{(6)}   &~~~& x_8^{(6)}  &~~~& x_{10}^{(6)} & ~~~ \\
x_1^{(7)} &~~~& x_3^{(7)} &~~~& x_5^{(7)} &~~~& x_7^{(7)} &~~~& x_{9}^{(7)} &~~~& x_{11}^{(7)} \\
 ~~~& x_2^{(8)} &~~~& x_4^{(8)} &~~~& x_6^{(8)}   &~~~& x_8^{(8)}  &~~~& x_{10}^{(8)} & ~~~ \\
x_1^{(9)} &~~~& x_3^{(9)} &~~~& x_5^{(9)} &~~~& x_7^{(9)} &~~~& x_{9}^{(9)} &~~~& x_{11}^{(9)} \\
 ~~~& x_2^{(10)} &~~~& x_4^{(10)} &~~~& x_6^{(10)}   &~~~& x_8^{(10)}  &~~~& x_{10}^{(10)} & ~~~ \\
x_1^{(11)} &~~~& x_3^{(11)} &~~~& x_5^{(11)} &~~~& x_7^{(11)} &~~~& x_{9}^{(11)} &~~~& x_{11}^{(11)} \\
 ~~~& x_2^{(12)} &~~~& x_4^{(12)} &~~~& x_6^{(12)}   &~~~& x_8^{(12)}  &~~~& x_{12}^{(6)} & ~~~
\end{array}$ \\ \begin{center} The first twelve rows of this region for $A_{11}$.
\end{center} \vspace{1em}

\noindent Cluster variable $y_i^{(j)}$ contains tiles $\TT_{i-j+1}$ through $\TT_{i+j-1}$ for a total
of $2j-1$ tiles.  Since we are considering only cluster variables between the two columns of
positive ones, we have graphs centered at $\TT_i$ for $2 \leq i \leq n+1$.  Thus our graph either \vspace{1.0em}

1) Contains no forbidden tiles. \\

2) Contains forbidden tiles including $\TT_1$ and tiles to the left. \\

3) Contains forbidden tiles including $\TT_{n+2}$ and tiles to the right. \\

4) Contains both sets of forbidden segments.

\vspace{1.0em} \noindent If there are forbidden tiles on only one side, then Lemma \ref{exciseGr}
allows us to edit out the forbidden segment.  However, if there are forbidden strips on both sides we would
encounter a problem if the segments we were deleting overlapped.  So consider a graph that
contains tiles $\TT_1\cup\dots \cup \TT_{n+2}$, along with $a$ tiles to the left and $b$ tiles to the right.  By Lemma \ref{exciseGr}, when
we get rid of the $a+b$ tiles on the two ends, we would also be deleting $a+b+2$ tiles from the
middle (including $\TT_1$ and $\TT_{n+2}$).  For the first $n+1$ rows of our lattice, the number of
total tiles is $\leq 2n+1$ and thus $a+b+n+2 \leq 2n+1$ which implies that $a+b \leq n-1 < n$.
Thus there cannot be an overlap in these first $n+1$ rows.  Furthermore, recalling our indexing $y_{i+1}=x_i$, we get as an added bonus
that the graph corresponding to $y_{i+1}^{(n+1)}=x_i^{(n+1)}$, which consists of $2n+1$ tiles, bijects to the same cluster variable as the graph of the single tile $\TT_{n-i-1}=T_{n-i}$ (resp. $\TT_{n-i}=T_{n+1-i}$) for $i$ even and $n$ even (resp. $n$ odd).

Consequently, we get a combinatorial interpretation for the first $n+1$ rows.  We will refer to the region
of the extended $A_n$-lattice that lies between the two columns of positive ones and in the first
$n+1$ rows underneath the initial variables as the $A_n$-lattice.  The diagonals of the
$A_n$-lattice inherit the properties of the $A_\infty$-lattice and this implies that all the
cluster variables do in fact appear in our lattice.  In particular, their denominators are in bijection with
positive roots of $A_n$'s root system and we have the desired corresponding graph for each of
them.  We now can complete the substitutions by letting $y_i =x_{i-1}$, for $i \in \{2,\dots, n+1\}$, which induces $y_{i}^{(j)}=x_{i-1}^{(j)}$ for $i\in \{2,\dots ,n+1\}$.  In particualr, the extended $A_n$ lattice reduces to doubly-periodic copies of the $A_n$ lattice containing graphs involving only tiles $\TT_2=T_1,\dots, \TT_{n+1}=T_n$. 
Thus Propositon \ref{CaseAn} is proven.  As a corollary of this argument, the diamond
condition implies that the $(n+2)$nd and $(n+3)$rd rows consists of the initial cluster variables
written in reverse order.
\end{proof}

\begin{Rem} Consider a new lattice $\{z_{i}^{(j)}\}$ consisting
of connected subsets of $\mathcal{T}_{A_n}$ such that $T_i \in z_i^{(j)}
\iff T_i$ appears in the graph associated to $x_i^{(j)}$ and add columns consisting of empty
sets on the left-hand and right-hand sides of this lattice.  This
lattice satisfies a tropical-like diamond condition where one of
the four following hold.
\begin{eqnarray*}
a = b \cup c \mathrm{~and~} d = b \cap c \\
a = b \cap c \mathrm{~and~} d = b \cup c \\
b = a \cup d \mathrm{~and~} c = a \cap d \\
b = a \cap d \mathrm{~and~} c = a \cup d
\end{eqnarray*}
\end{Rem}

\begin{Rem}
Such lattices are known as frieze patterns, and were studied by Conway and Coxeter \cite{ConCox} in the 1970's.  Such patterns have also been studied in connection with cluster algebras in work of Caldero \cite{Caldero2} and work of Propp \cite{MarkPropp}.  These lattices are also special cases of the bipartite belt described in \cite{ClusIV}; each row of the lattice corresponds to a seed of the belt.
\end{Rem}

\begin{Rem}
Hugh Thomas \cite{HTPers} brought it to the author's
attention that one can also derive the above lattices via the algorithm for
constructing the Auslander-Reiten quiver \cite{AssocAlg} starting from
projective representations; in particular the pattern of denominator vectors agrees with the dimension vectors of the indecomposables in the AR quiver.
\end{Rem}

\section{$C_n$}

The Lie algebra $C_n$ has the following Dynkin diagram
$$\bullet\Rightarrow
 \bullet\line(1,0){3}\bullet\line(1,0){3}\bullet\line(1,0){3}\bullet\line(1,0){3}
 \bullet\line(1,0){3} \dots \dots \line(1,0){3}\bullet$$
\noindent and thus the bipartite exchange matrix is:
$$\begin{bmatrix}
 0  &  2 &  0  &  0  & \dots & 0  & 0 \\
-1  &  0 & -1  &  0  & \dots & 0  & 0 \\
 0  & 1 &  0  & 1  & \dots & 0  & 0 \\
 0  &  0 & -1  &  0  & \dots & 0  & 0 \\
\dots & \dots & \dots & \dots & \dots & \dots & \dots \\
0  &  0 & 0  & 0  & \dots & (-1)^{n+1} & 0 \\
\end{bmatrix}.$$
\noindent To build the corresponding graphs we let $\mathcal{T}_{C_n}$ be
identical to $\mathcal{T}_{A_n}$ except that tile $T_1$ now has
weights of $x_2$ and $x_2$ opposite each other instead of a lone
weighted edge. This change to $T_1$ corresponds to the change to
the exchange polynomial associated to label $1$ in the seed of
this cluster algebra.

\vspace{1em}\begin{center} \includegraphics[width = 3in , height
= 0.6in]{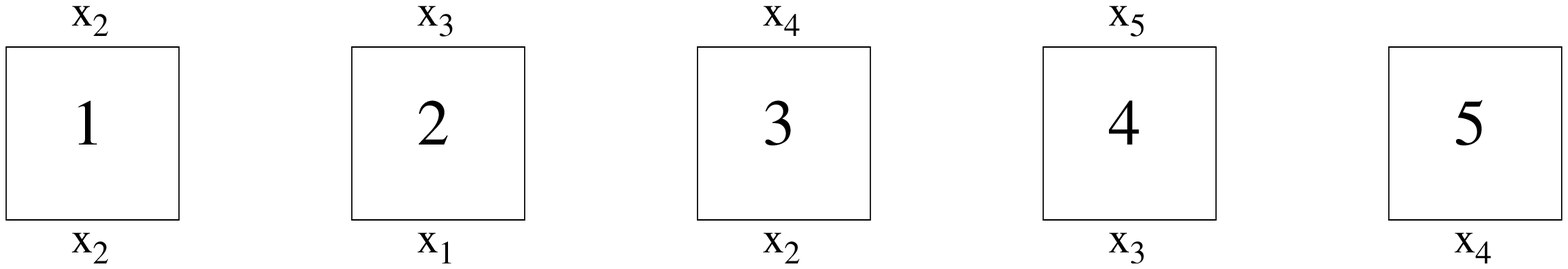}\\ Tiles for $C_5$.
\end{center} \vspace{1em}

We use gluing rule \ref{Gluing} again which leads us to a
collection similar to $\mathcal{G}_{A_n}$ except now tile $T_1$
can connect to tile $T_2$ on either side. Thus the collection of
possible graphs, $\mathcal{G}_{C_n}$ corresponds to the sets of
the form

$$\{T_i,T_{i+1},T_{i+2},\dots, T_{j-1},T_j\}$$ for $1 \leq i < j
\leq n$ or multisets of the form

$$\{T_i,T_{i-1},T_{i-2},\dots, T_3, T_2, T_1, T_2, T_3, T_{j-1},T_j\}$$ for $2 \leq i \leq j
\leq n$.  This collection $\mathcal{G}_{C_n}$ has the same
cardinality as the collection of non-initial cluster variables for
a cluster algebra of type $C_n$ and thus the collection of positive
roots for a root system of type $C_n$, as in the last case
\cite{ClustII,Kac}.

\begin{Prop} \label{CaseCn} The set of graphs $\mathcal{G}_{C_n}$ is in bijection
with the set of non-initial cluster variables for a coefficient-free cluster algebra
of type $C_n$ such that the statement of Theorem \ref{vargraph}
holds. \end{Prop}

This can be proved quickly by using the folding procedure as in
\cite{ClustII}.  We identify $A_{2n-1}$ with $C_n$ by letting $x_k
= x_{n+1-k}$ for $k \in \{1,\dots, n-1\}$.  We let $x_n = x_1$ and
let $x_k = x_{k-n+1}$ for $k \in \{n+1,\dots, 2n-1\}$. Our lattice
will contain repeats but we can restrict our list to the right half, including the central axis, to obtain the correct number of graphs.  Thus Proposition
\ref{CaseAn} implies Proposition \ref{CaseCn}.

\vspace{1em}\begin{center} \includegraphics[width = 3in ,
height = 3in]{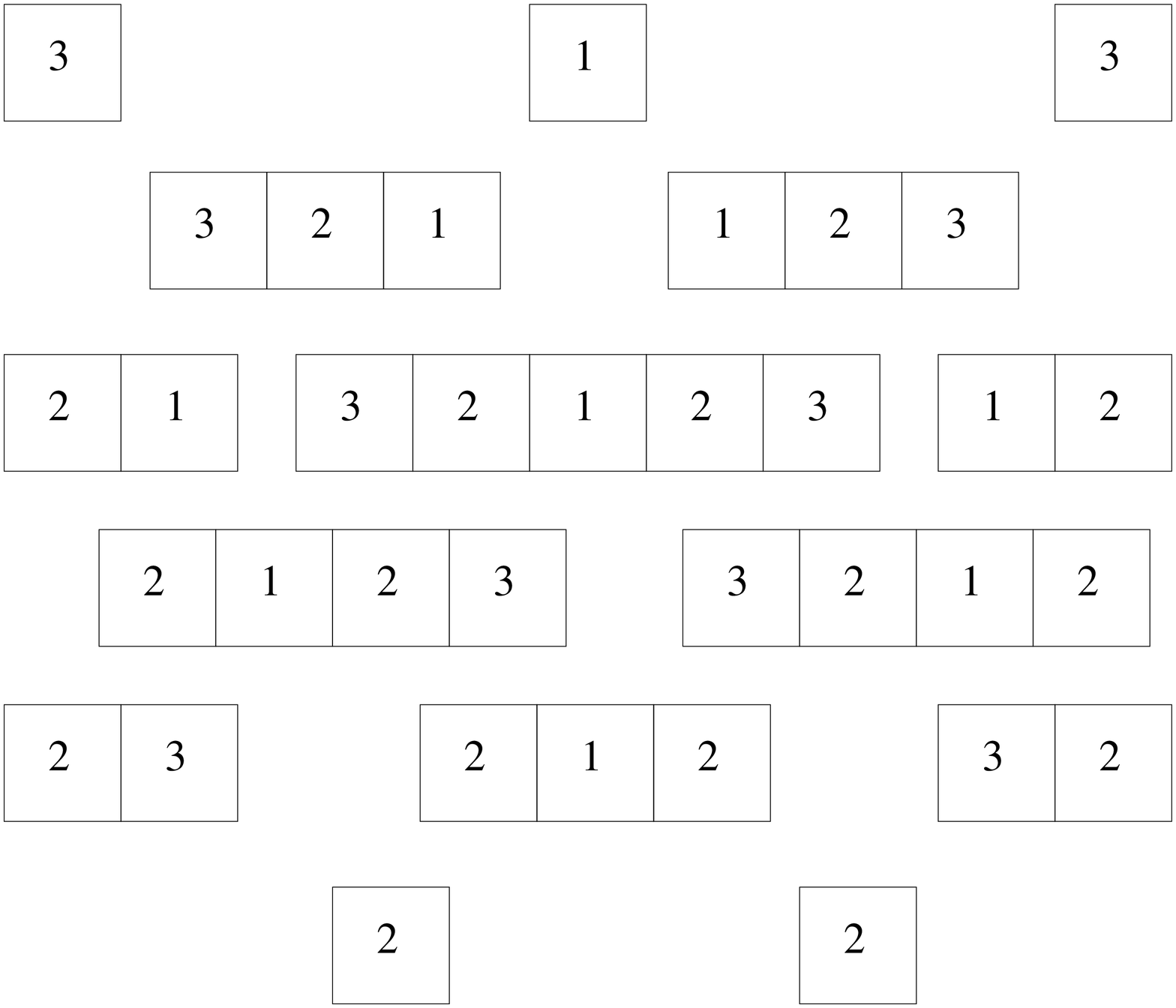} \\ The collection $\mathcal{G}_{C_3}$ with duplicates.
\end{center} \vspace{1em}

\section{$B_n$ and $D_n$} In the previous two cases, all of the exchange polynomials had degree two or less.
For the cases of $B_n$ and $D_n$, exactly one of the exchange
polynomials has degree three.  We will deal with such exchanges by including hexagons as potential tiles.  We start with the case of $B_n$, which is a folded version of the simply-laced $D_n$ case.  By folding, our proofs will require less notation and as we will see, the $D_n$ case has a symmetry such that we can easily derive this case from the results for $B_n$.

In the case of $B_n$, the Dynkin diagram is
$$\bullet\Leftarrow
 \bullet\line(1,0){3}\bullet\line(1,0){3}\bullet\line(1,0){3}\bullet\line(1,0){3}
 \bullet\line(1,0){3} \dots \dots \line(1,0){3}\bullet$$
\noindent and thus the bipartite exchange matrix is:
$$\begin{bmatrix}
 0  &  1 &  0  &  0  & \dots & 0  & 0 \\
-2  &  0 & -1  &  0  & \dots & 0  & 0 \\
 0  & 1 &  0  & 1  & \dots & 0  & 0 \\
 0  &  0 & -1  &  0  & \dots & 0  & 0 \\
\dots & \dots & \dots & \dots & \dots & \dots & \dots \\
0  &  0 & 0  & 0  & \dots & (-1)^{n+1} & 0 \\
\end{bmatrix}.$$
\noindent We will now use the notation $T_1$ through $T_n$ to refer to a collection of tiles, $\mathcal{T}_{B_n}$, related to $B_n$.  We construct $\mathcal{T}_{B_n}$ from
$\mathcal{T}_{A_n}$ by first replacing $T_2$ with a hexagon having
weights $1$, $x_1$, $1$, $x_1$, $1$, and $x_3$ in clockwise order starting from the top. We let $T_1$ be a trapezoid with a single weighted edge of $x_2$ on its northern side.  Note that $T_1$ is homeomorphic to its previous definition.  Then for all $i > 2$ we define $T_i$, for type $B$, as a counter-clockwise rotation of the $A_n$-tile $T_i$, including the boundary tile $T_n$ which have a single weighted edge of $x_{n-1}$ on its eastern side.

\vspace{1em}\begin{center} \includegraphics[width = 2.45in ,
height = 1.6in]{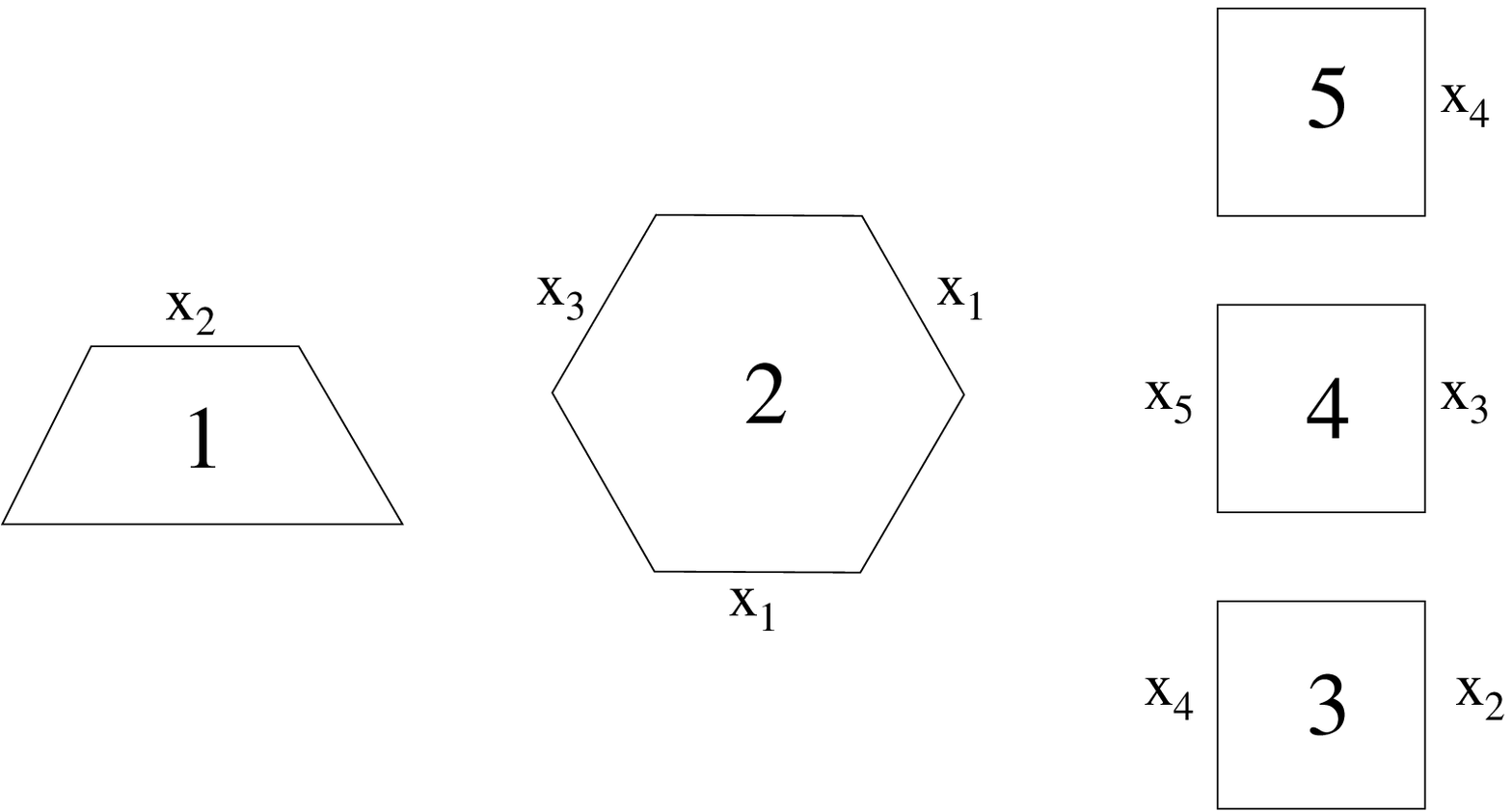}\\ Tiles for $B_5$.
\end{center} \vspace{1em} 

The gluing rule will be more complicated now that hexagons are
involved.  As a first approximation, the set of graphs $\mathcal{G}_{B_n}$ will include any
graphs that can be constructed from $\mathcal{T}_{B_n}$ while
conforming to Rule \ref{Gluing}.  Again we are not allowing rotation or reflections of the tiles so they must be connected in the orientaitions as described above.  Any such graph will resemble either a \emph{tower} of tiles $T_a$ through $T_b$ for $3\leq a \leq b \leq n$, a \emph{base} involving hexagon $T_2$ with or without trapezoid $T_1$ on its western side, or may be a complex of a tower beginning with $T_3$ on top of a base.

In addition, we enlarge the set of $\mathcal{G}_{B_n}$ by allowing any graphs that obey the following second rule:

\begin{Rule} \label{Hexagons} The trapezoidal tile $T_1$ may appear twice if and only if the lift of the graph to $\mathcal{G}_{B_\infty}$ (i.e. $n$ arbitrarily large) has one of the following three forms:
\end{Rule}

\vspace{1em}\begin{center} 
\includegraphics[width = 1.1in ,
height = 0.8in]{BB121.eps}
\hspace{2em}
\includegraphics[width = 1.1in ,
height = 1.8in]{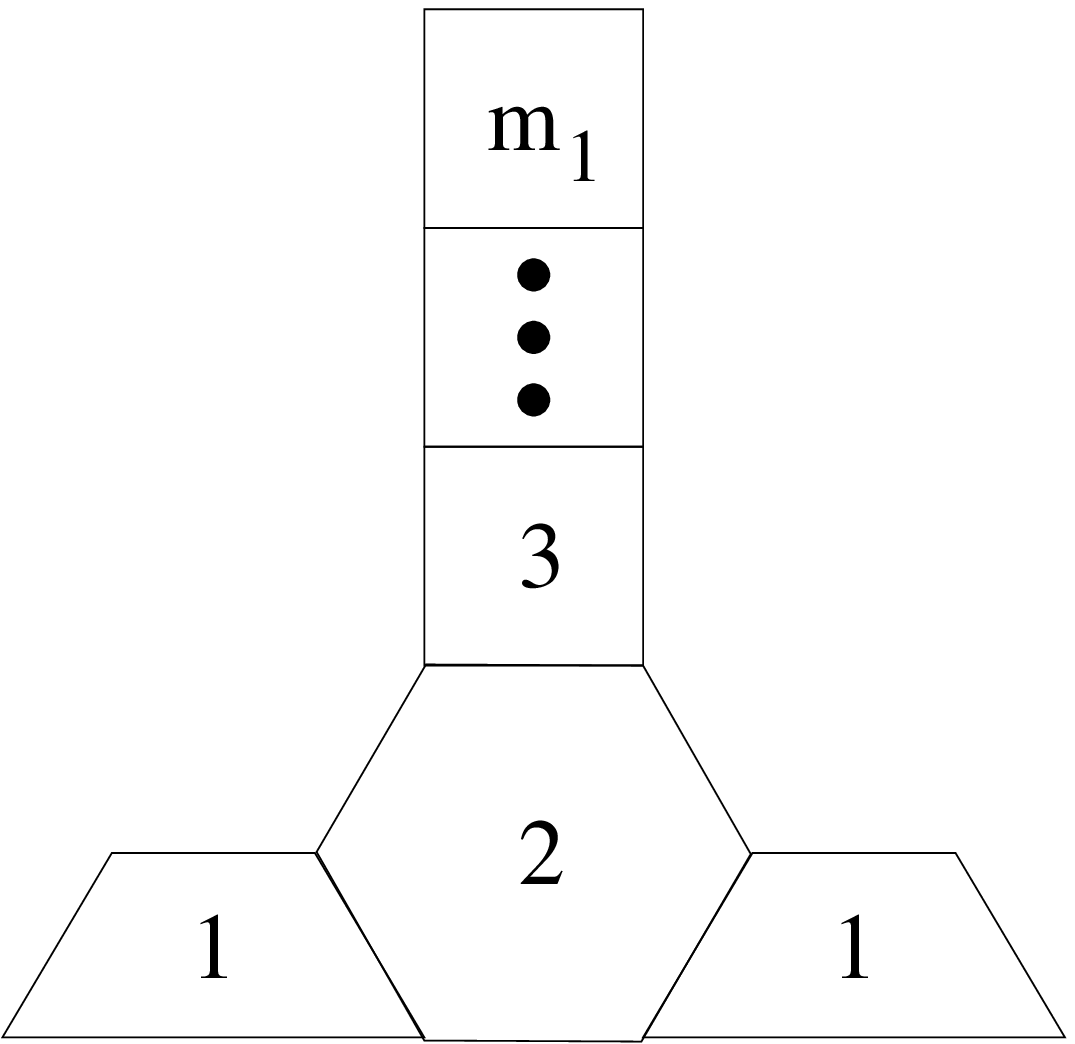}
\hspace{2em}
\includegraphics[width = 1.1in ,
height = 1.1in]{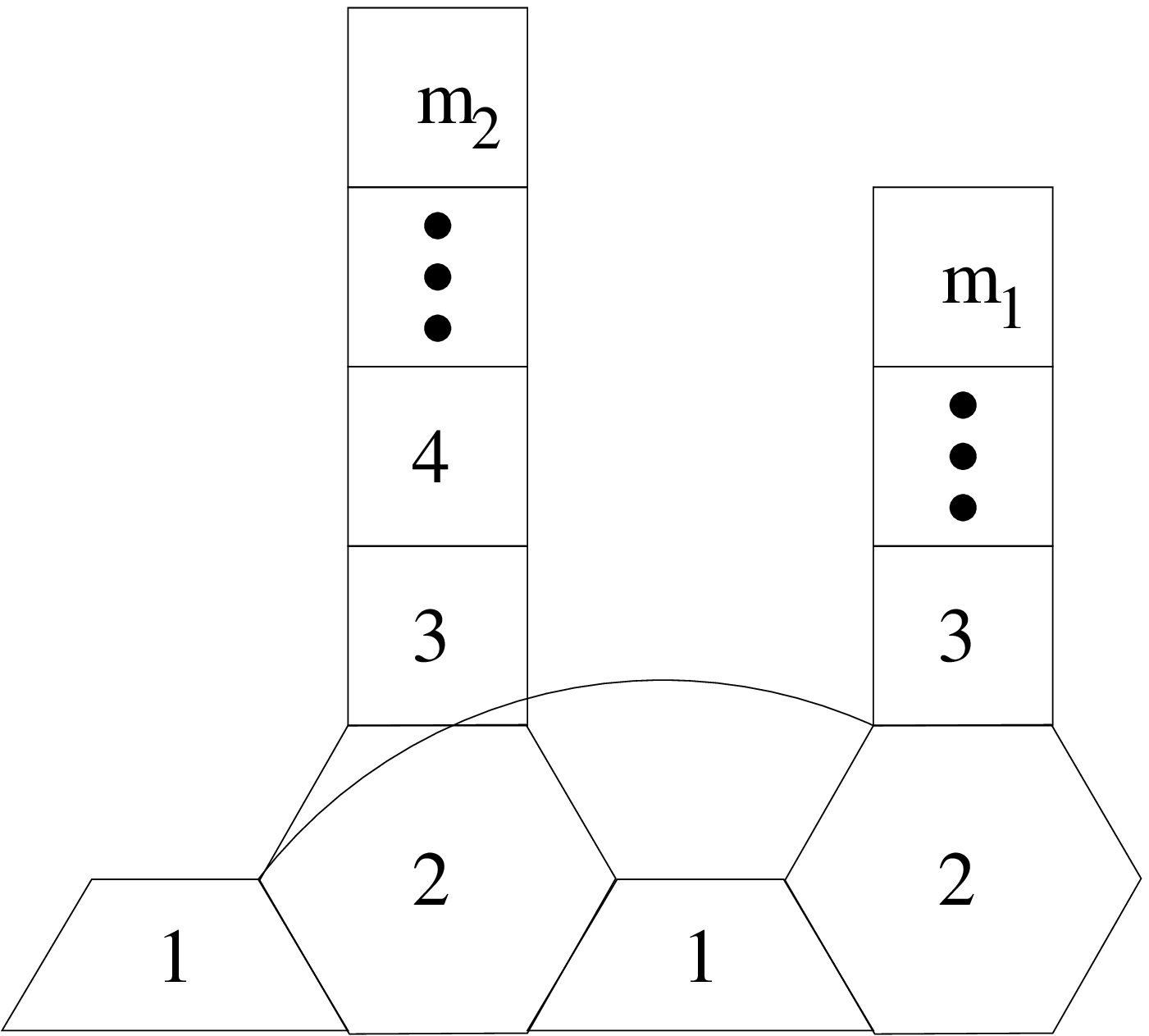}
\end{center}\vspace{1em}

\noindent where $3 \leq m_1 < m_2$ and $m_1,~m_2$ are both odd.

\begin{Rem} Notice that Rule $1$ is now broken when we connect 
trapezoid $T_1$ to a hexagon $T_2$ on its left.  Furthermore, in the last of these cases, we have adjoined an additional arc which had not been allowed or required in previous examples.  However, there is precedent for using such additional arcs, see Section $3$ of \cite{MusPropp}.  

It will develop that we can project these graphs, consisting of two towers, down to the $B_n$-lattice by excision around tile $T_{n+1}$ just as in the $A_n$ case.  Thus the fact that these are an odd number of tiles in each tower, with the larger tower on the left, greatly limits the set of such graphs.\end{Rem}

One can check that the collection of graphs
$\mathcal{G}_{B_n}$ obeying Rule $1$ \emph{or} Rule $2$ has the cardinality equal to the number of positive roots for $B_n$.  Further we will
prove, after excision, that Theorem \ref{vargraph} is satisfied by these defintions.

\begin{Prop} \label{CaseBn} The set of graphs $\mathcal{G}_{B_n}$ is in bijection
with the set of non-initial cluster variables for a coefficient-free cluster algebra
of type $B_n$ such that the statement of Theorem \ref{vargraph}
holds.
\end{Prop}

\begin{proof} Analogous to the $A_n$ case, we will first prove the result for the $B_\infty$ case, that is we assume that $n$ is arbitrarily large so that for $i\geq 3$ tile $T_i$ always has exactly two weighted edges ($x_{i-1}$ on its east and $x_{i+1}$ on its west). This greatly simplies the proofs by allowing easier notation and bypassing case-by-case analysis.  We will later discuss how to obtain the result for a specific $B_n$ from such graphs.  We create a semi-infinite lattice whose entries satisfy a deformed diamond condition.

\vspace{1em}$\begin{array}{cccccccccccc}
x_1 &~~~& x_3 &~~~& x_5 &~~~& \dots &~~~&  &~~~&  \\
 ~~~& x_2 &~~~& x_4 &~~~& x_6   &~~~& \dots  &~~~&  & ~~~ \\
x_1^{(1)} &~~~& x_3^{(1)} &~~~& x_5^{(1)} &~~~& \dots &~~~&  &~~~&  \\
 ~~~& x_2^{(2)} &~~~& x_4^{(2)} &~~~& x_6^{(2)}   &~~~& \dots  &~~~&  & ~~~ \\
x_1^{(3)} &~~~& x_3^{(3)} &~~~& x_5^{(3)} &~~~& \dots &~~~&  &~~~&  \\
 ~~~& x_2^{(4)} &~~~& x_4^{(4)} &~~~& x_6^{(4)}   &~~~& \dots  &~~~&  & ~~~ \\
 \dots & \dots & \dots & \dots
\end{array}$ \\ \begin{center} The lattice for $B_\infty$.
\end{center}\vspace{1em}

\noindent Without the boundary on the right, any collection of four variables

\begin{center}$\begin{array}{ccc}
&a& \\
b&&c \\
&d&
\end{array}$\end{center}

\noindent such that $b \not = x_1^{(j)}$ and $c \not = x_2^{(j)}$ will satisfy $ad - bc = 1$.  A  diamond such that $c = x_2^{(j)}$ will satisfy the truncated condition $ad - c = 1$, and a
diamond which contains $b = x_1^{(j)}$ will satisfy the relation $ad - b^2c = 1$.  As before, we say that a Laurent polynomial bijects to graph $G$, which we denote as $x_i^{(j)}\leftrightarrow G$, if $x_i^{(j)} = P(G)/{x_1^{\alpha_1}\cdots x_n^{\alpha_n}}$ where $P(G)$ is the matching polynomial of $G$ and $\alpha_i$ encodes the number of occurences of tile $T_i$ in graph $G$.

Given this setup along with the initial assignments of
$x_i^{(0)} = x_i$ for $i\geq 1$, we directly verify that
\begin{center}
$x_1^{(1)} ~\longleftrightarrow~ $ \includegraphics[width = 0.4in , height = 0.4in]{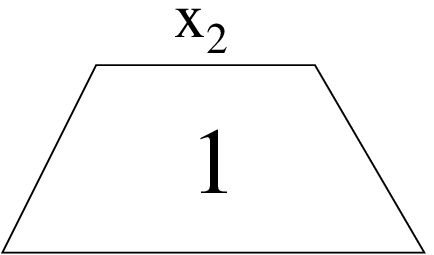}
\hspace{0.1in}, \hspace{0.6in} $x_2^{(2)} ~\longleftrightarrow~ $ \includegraphics[width = 0.7in ,
height = 0.7in]{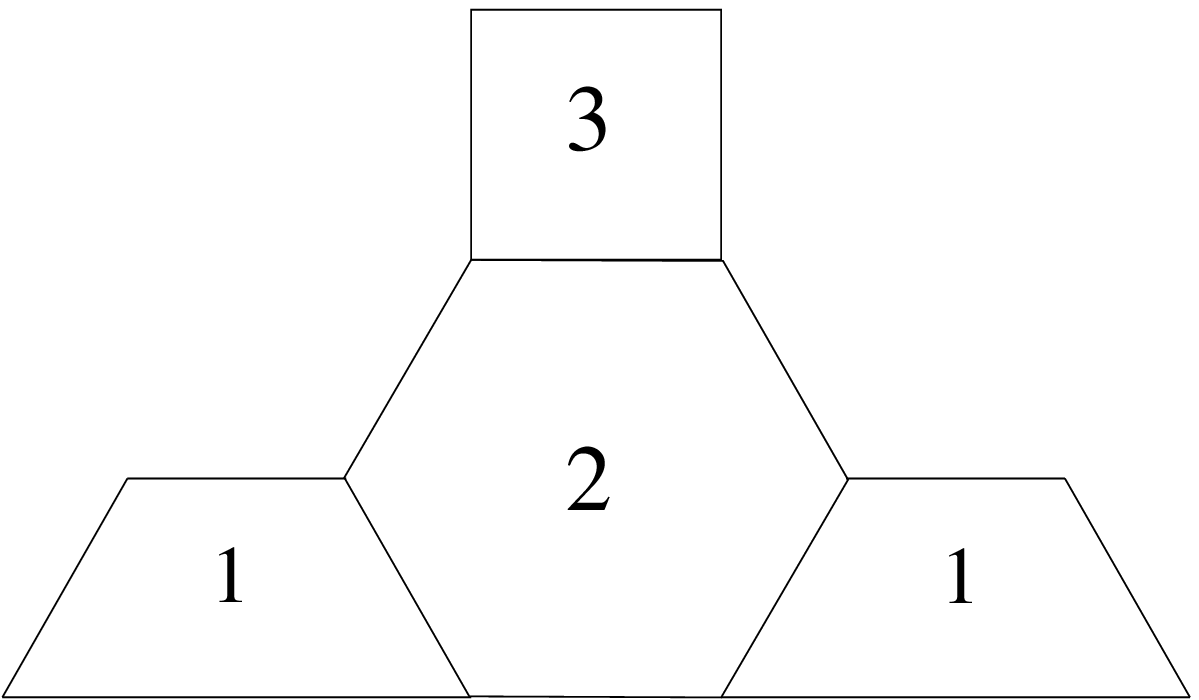} \hspace{0.05in},
\\ \vspace{1.5em} $x_1^{(3)} ~\longleftrightarrow~ $ \includegraphics[width = 0.7in , height =
0.7in]{BB123} 
\hspace{0.1in}, \hspace{0.3in}
 $x_2^{(4)} ~\longleftrightarrow~ $ \includegraphics[width = 0.7in , height = 0.7in]{BB121235}
\hspace{0.05in},
\end{center}\vspace{1.5em}
where the weights of the edges are as dictated by the definitions of tiles $T_1$ through $T_n$.
Additionally for $i-j \geq 2$, the only initial variables used to determine $x_i^{(j)}$ are
$\{x_3,x_5,x_7,\dots\}$ and thus we recover the regular diamond pattern used in the $A_\infty$
case.  Consequently, we immediately obtain

\begin{Lem} \label{AAregion}

$x_i^{(j)} ~\longleftrightarrow~ $ \includegraphics[width = 1.5in , height =
0.3in]{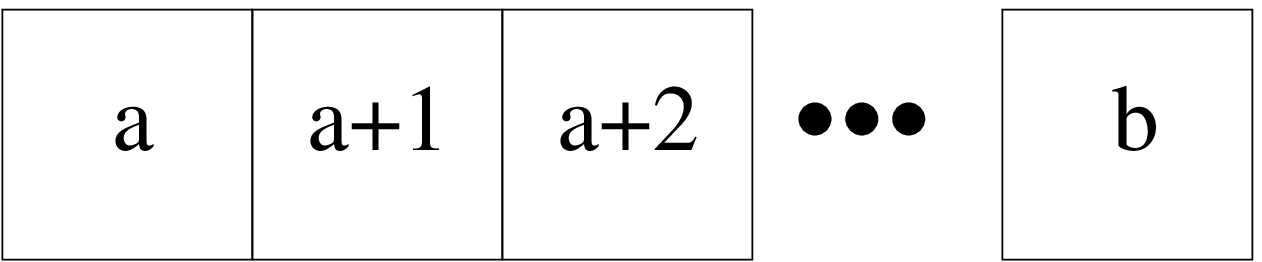} \vspace{1.5em} \hspace{1em} for $a = i-j+1, ~b = i+j-1$ when $i - j \geq
2$. We also find that $3 \leq a \leq b$.
\end{Lem}

We proceed with the rest of the proof in three steps.  The first two steps are proved inductively
by using Lemma \ref{AAregion} as well as $\{x_1^{(1)}, x_2^{(2)}, x_1^{(3)}, x_2^{(4)} \}$ as a
base case. We will prove the inductive step via the usual diamond condition $ad - bc = 1$ which
will hold for the diagonals $i - j = 0$ and $i - j = - 2$ while $i \geq 2$.

\begin{Lem} \label{xii}

$x_i^{(i)} ~\longleftrightarrow~ $ \includegraphics[width = 1in , height = 1in]{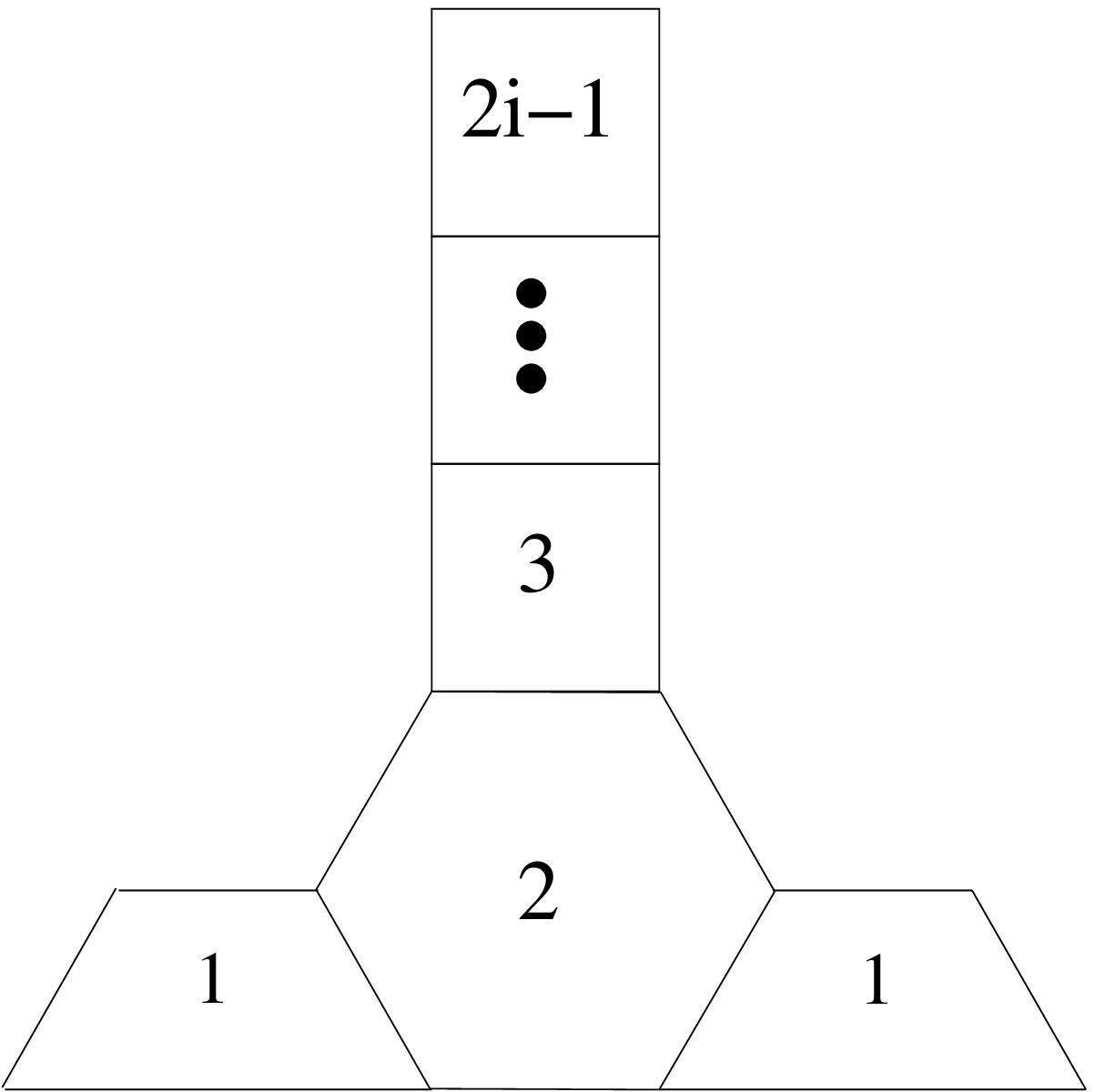}
\hspace{0.3em}  for $i \geq 2$.

\end{Lem}

\begin{Lem} \label{xii+2}

$x_i^{(i+2)} ~\longleftrightarrow~ $ \includegraphics[width = 1.5in , height = 1in]{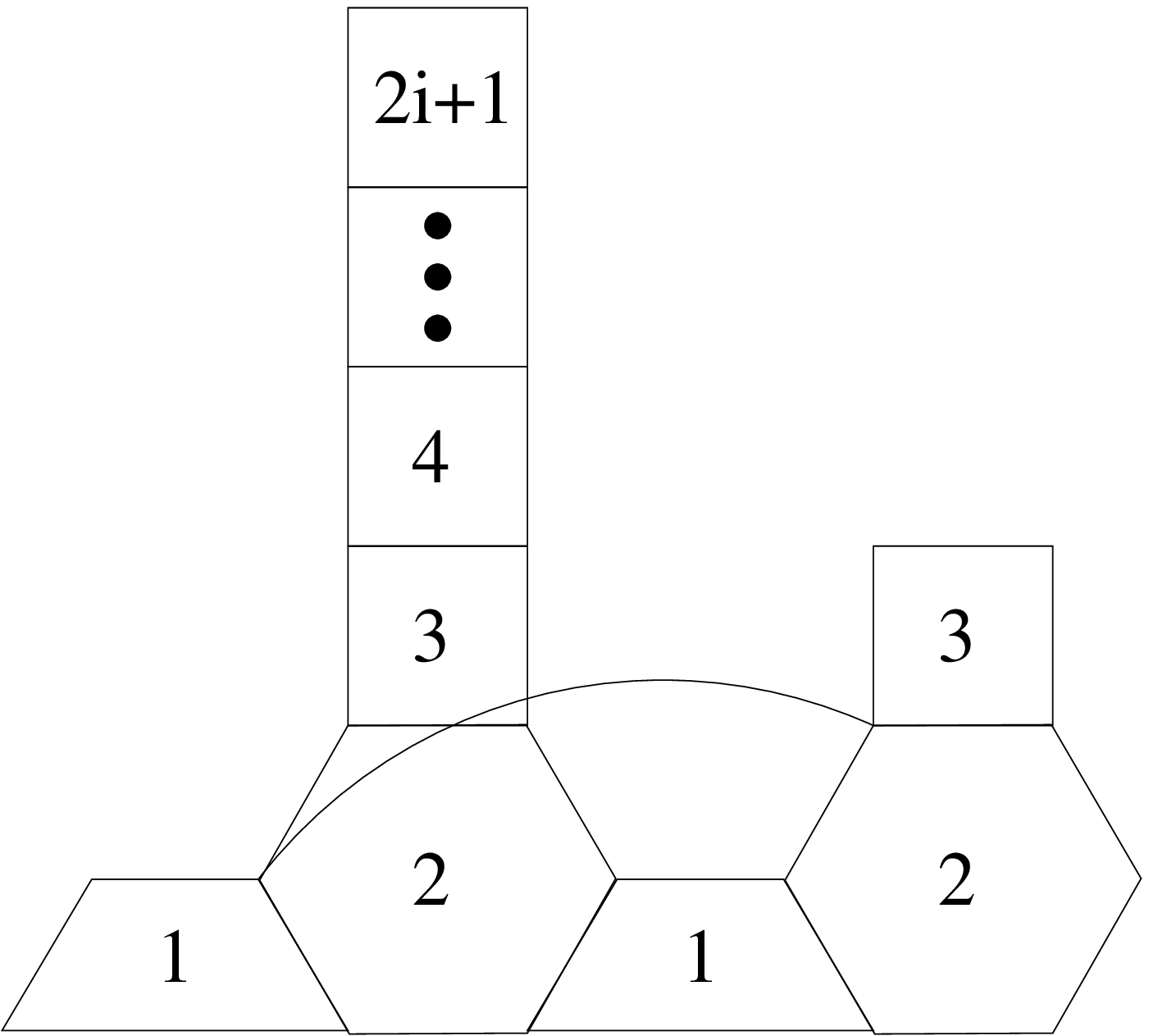}
\hspace{0.3em}  for $i \geq 2$.

\end{Lem}

The proof of these Lemmas will prove Proposition \ref{CaseBn} for all $x_i^{(j)}$ such that $i - j
\geq -2$. We must now use variants of the diamond conditions ($ad - c = 1$ and $ad- b^2c = 1$) to
extend down columns $x_1^{(j)}$ and $x_2^{(j)}$ respectively.  But to continue to have new entries
to use as $c$ in the relation we must continually extend down diagonals as we extend down the
columns. Consequently, we proceed to prove the following three results for $j = 3$, then for $j =
4$, and so on by induction.

\begin{Lem} \label{restofdiagonals}

$x_1^{(2j+1)} ~\longleftrightarrow~ $ \includegraphics[width = 1in , height = 1in]{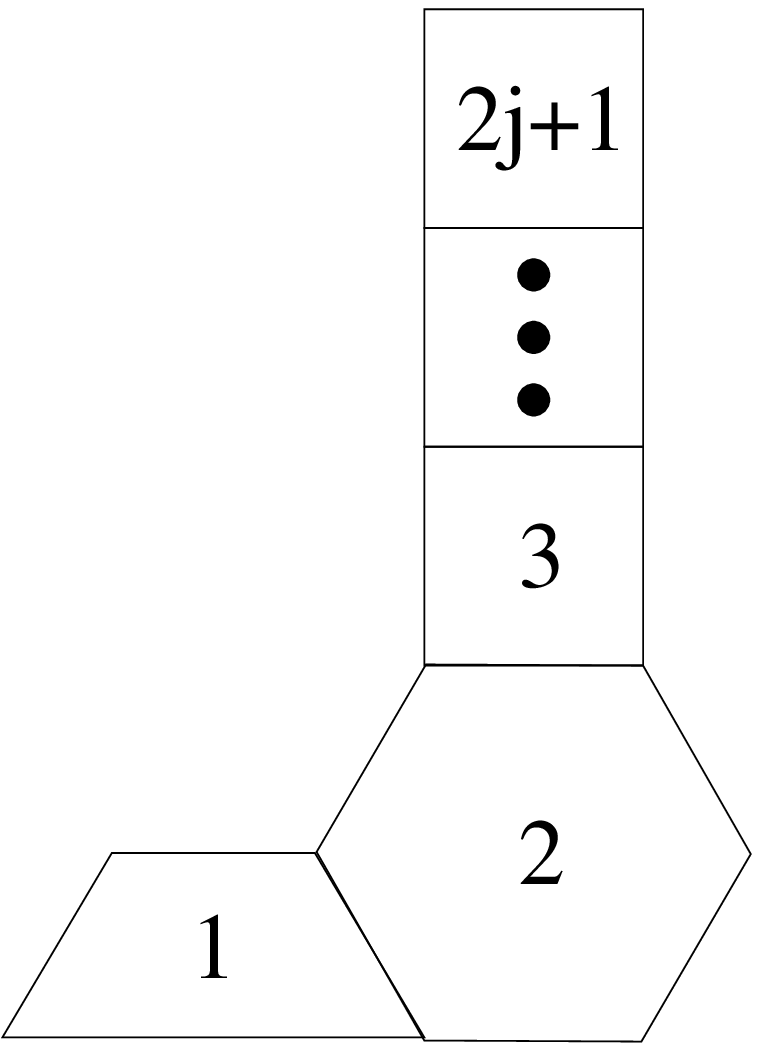}
\hspace{0.3in} for $j \geq 1$,\\ \vspace{1em}

$x_2^{(2j)} ~\longleftrightarrow~ $\includegraphics[width = 1.5in , height = 1in]{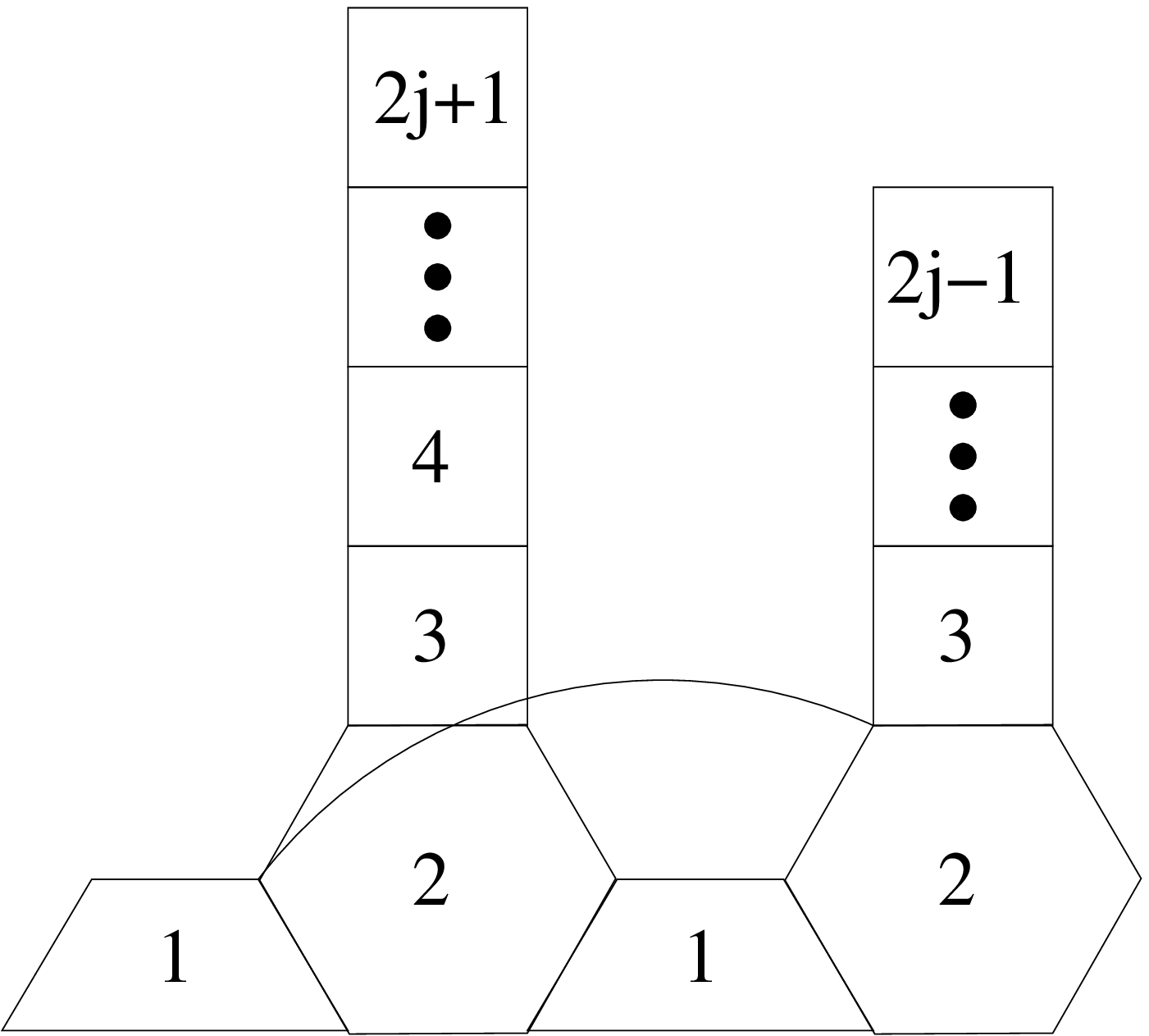}
 \hspace{0.3in} for $j \geq 2$, \\ \vspace{1em}

$x_i^{(i+2j)} ~\longleftrightarrow~ $ \includegraphics[width = 1.5in , height = 1in]{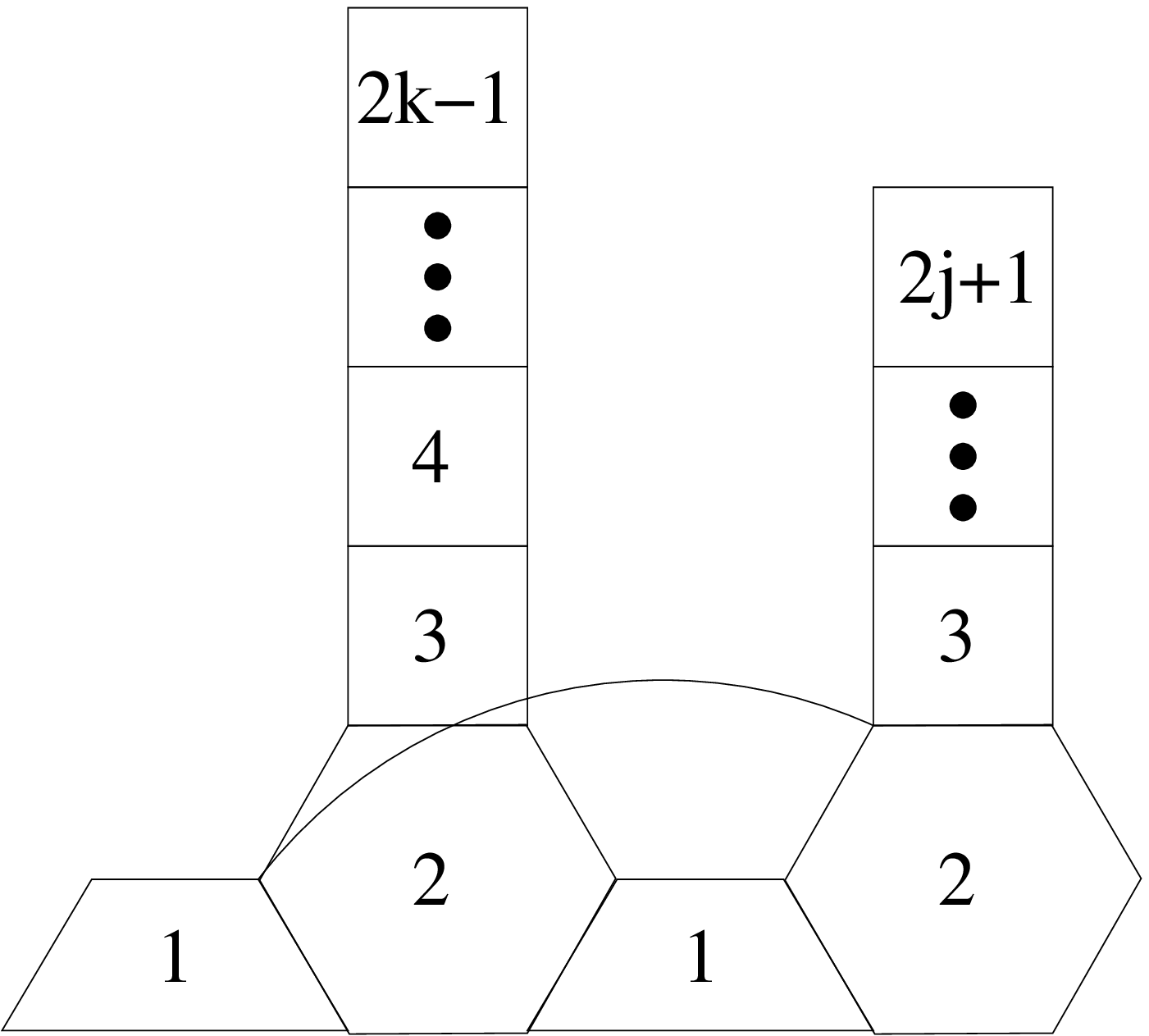}
\hspace{0.3in} where $k = i+j$ for $i \geq 2$ and $j \geq 1$.
\end{Lem}

\vspace{2.5em}\begin{center} \includegraphics[width = 2.5in , height = 2.5in]{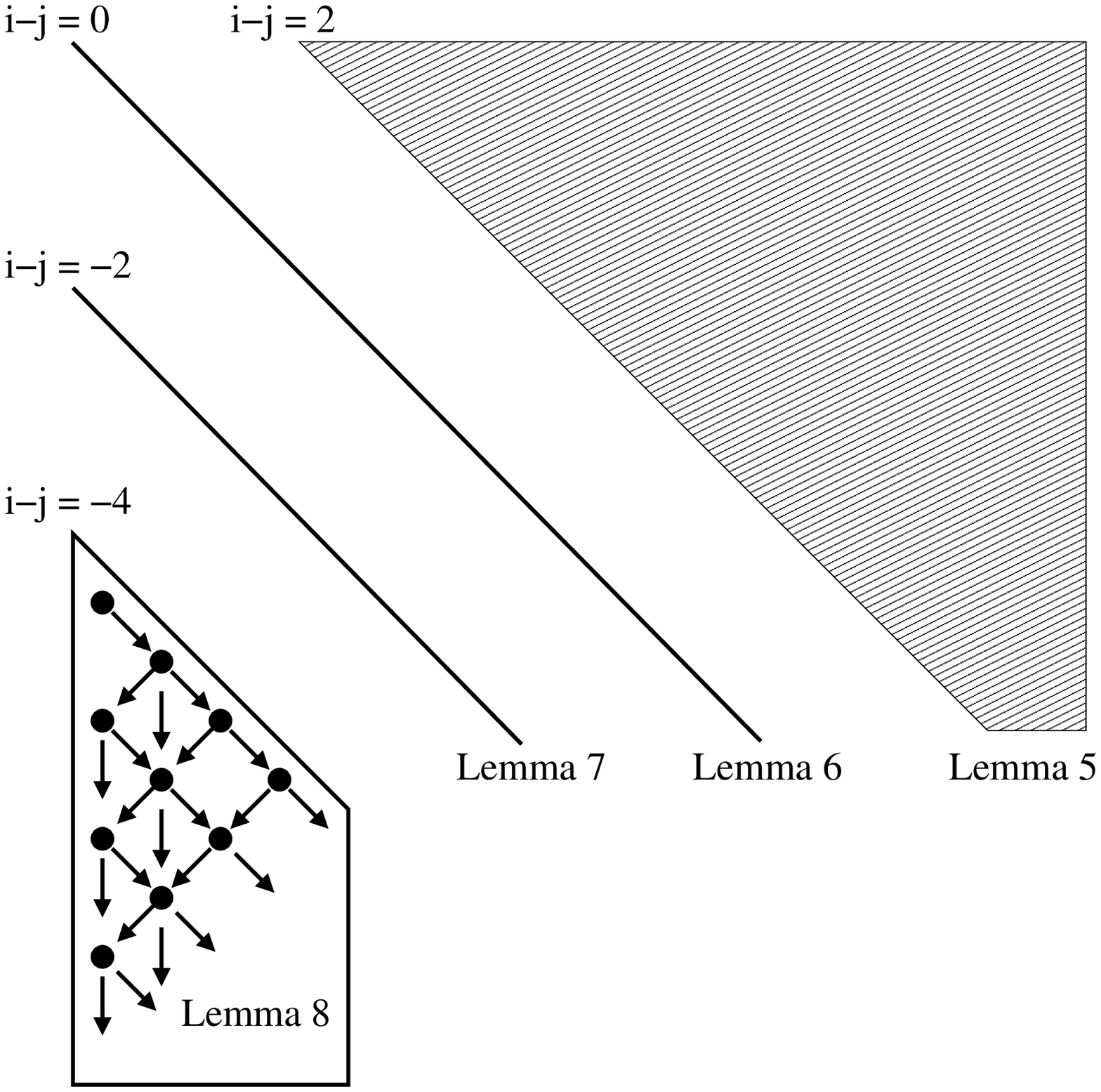}\\
A model of how these Lemmas fit together and relate to the $B_\infty$-lattice.
\end{center} \vspace{0.5em}

\noindent With the proof of Proposition \ref{CaseBn} now broken down into manageable chunks, we
proceed to prove the Lemmas.

\vspace{1em} \noindent\bf Proof of Lemma \ref{xii}.\rm \hspace{0.05in} By Lemma \ref{AAregion}, we have
that the northeast portion of the lattice is filled in, and we use the entries on diagonal $i-j=2$ and the base cases of $x_1^{(1)}$ and $x_2^{(2)}$ to extend to the rest of diagonal $i-j=0$ by the diamond condtion.  Assuming that cluster variables are

\vspace{1em}

$a ~\longleftrightarrow~ $ \includegraphics[width = 1.5in , height = 0.3in]{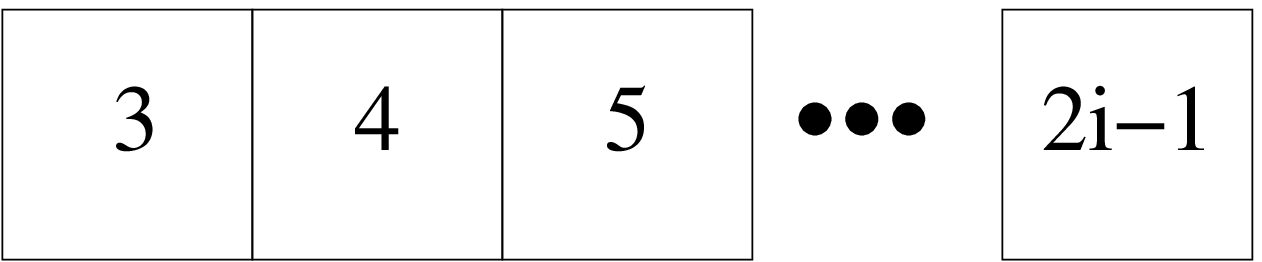} \\
\vspace{1em}

$b ~\longleftrightarrow~ $\includegraphics[width = 1in , height = 1in]{Bxii.eps} \\ \vspace{1em}

$c ~\longleftrightarrow~ $ \includegraphics[width = 1.5in , height = 0.3in]{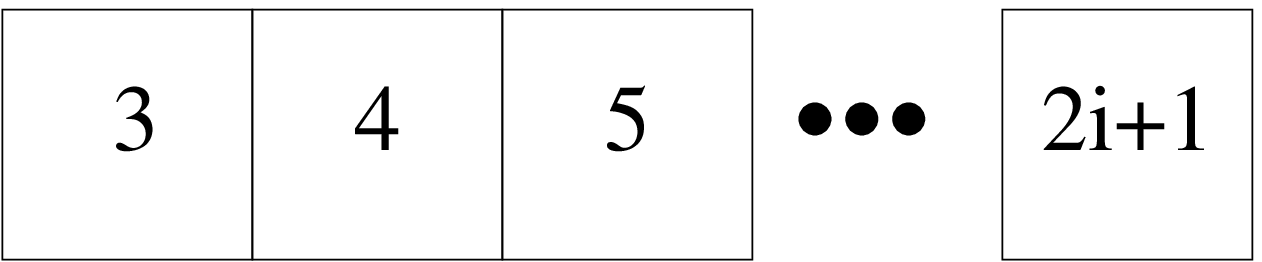}
\vspace{1em}

\noindent we wish to show that $d ~\longleftrightarrow~ $ \includegraphics[width = 1in , height =
1in]{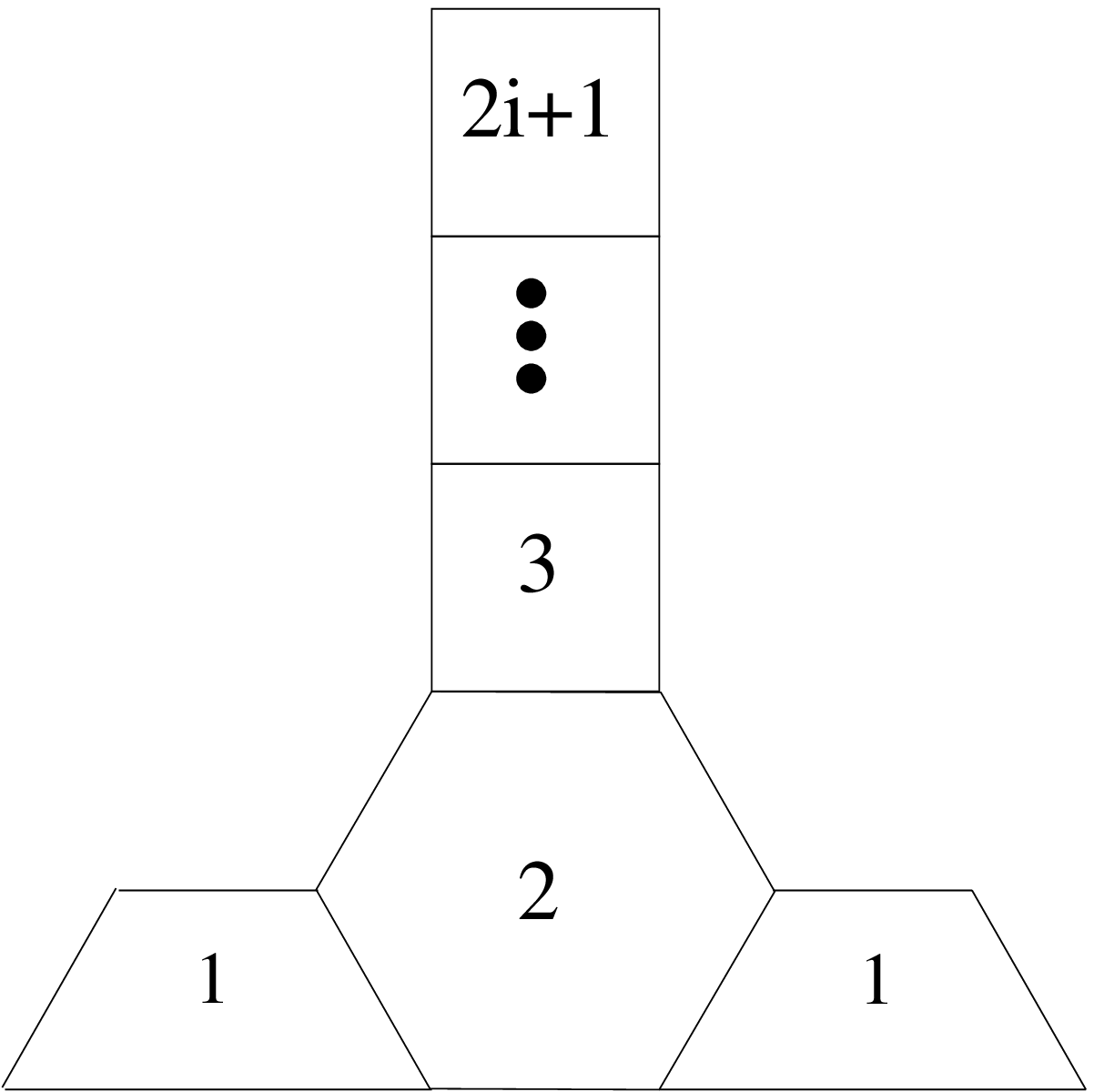}, \hspace{1.5em} given the diamond relation $ad = bc + 1 $.  First of all, we see
that the occurrences of tiles match up on each side of the equal sign, which implies that the denominators agree appropriately.
It suffices to show the weighted number of matchings also match up accordingly.  Any pair
of matchings of the graphs \includegraphics[width = 1.0in , height = 0.2in]{Lem6a.eps} and
\includegraphics[width = 0.8in , height = 0.8in]{Lem6d.eps} can be decomposed into a
pair of matchings on graphs \includegraphics[width = 1.0in , height = 0.2in]{Lem6c.eps}
and
\includegraphics[width = 0.8in , height = 0.8in]{Bxii.eps} except for one matching.  The logic is identical to that of Lemma \ref{diag}.  Here we swap the tops of the towers and note if the top horizontal edge of the hexagon is used (instead of the NW and NE diagonal edges), then completely swapping the two towers is permissible.

\vspace{0.5em} \noindent This extraneous indecomposable pairing is the pair 
\includegraphics[width = 1in , height = 0.9in]{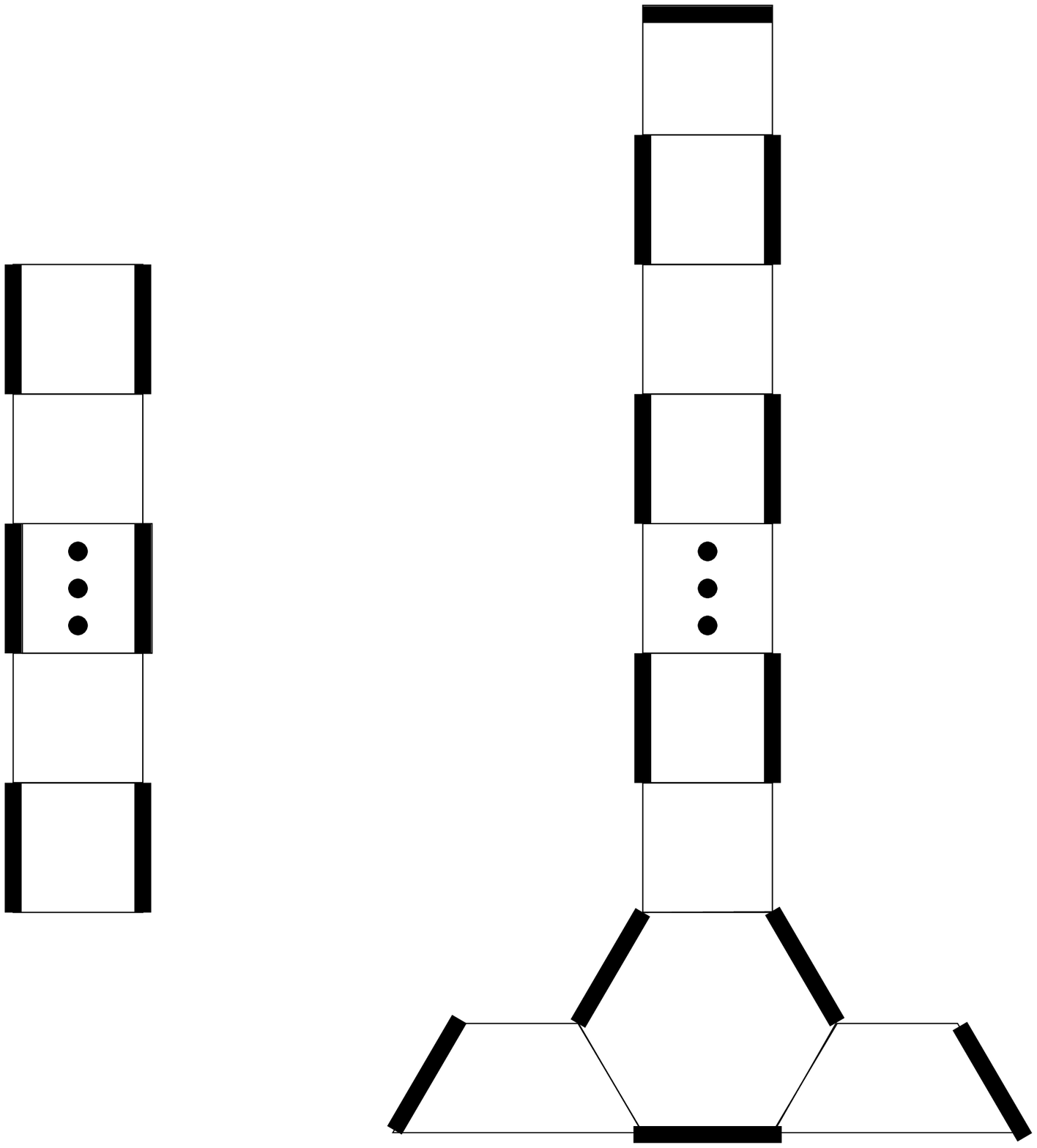} and it has exactly the correct
weight $(x_1^2x_2x_3^2x_4^2\cdots x_{2i-1}^2x_{2i}x_{2i+1})$ and thus the Lemma is proved.

\vspace{1em} \noindent\bf Proof of Lemma \ref{xii+2}.\rm \hspace{0.05in}  The proof of this
Lemma is analogous except that we shift the diamond pattern so that

\vspace{1em}

$a ~\longleftrightarrow~ $ \includegraphics[width = 1in , height = 1in]{Bxii.eps} \hspace{0.1in},
\hspace{0.6in} $d ~\longleftrightarrow~ $ \includegraphics[width = 1.5in , height =
1in]{Bxii2.eps} \\ \vspace{1em}

$b ~\longleftrightarrow~ $\includegraphics[width = 1.5in , height = 1in]{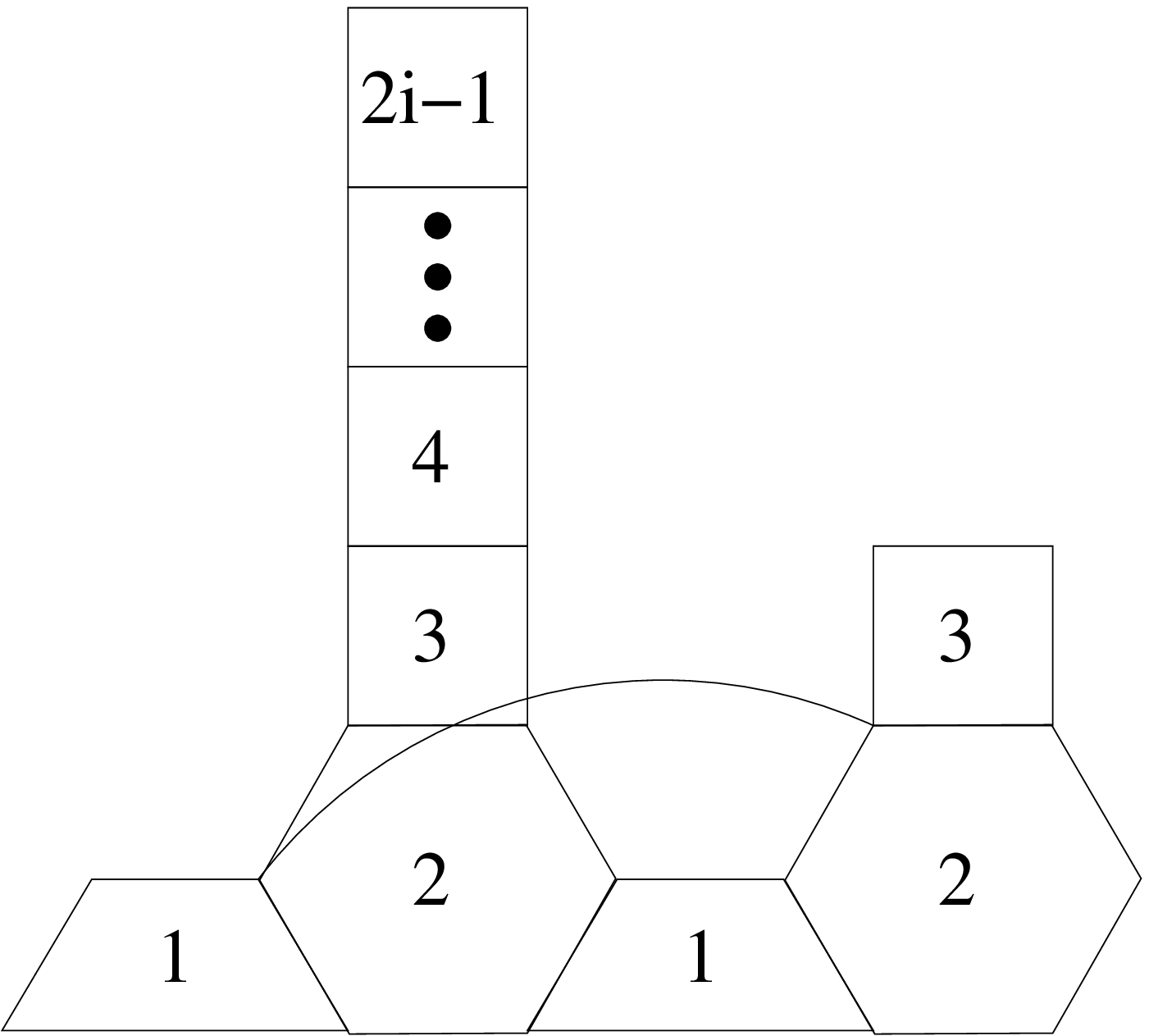}
\hspace{0.1in}, \hspace{0.6in} $c ~\longleftrightarrow~ $ \includegraphics[width = 1in , height =
1in]{Lem6d}\hspace{0.1in}.

\vspace{1em} \noindent Again, we have a bijection by swapping the tops of the (left) towers, and there is exactly one extraneous pair of matchings: 

\vspace{1em}\includegraphics[width = 1in, height = 0.9in]{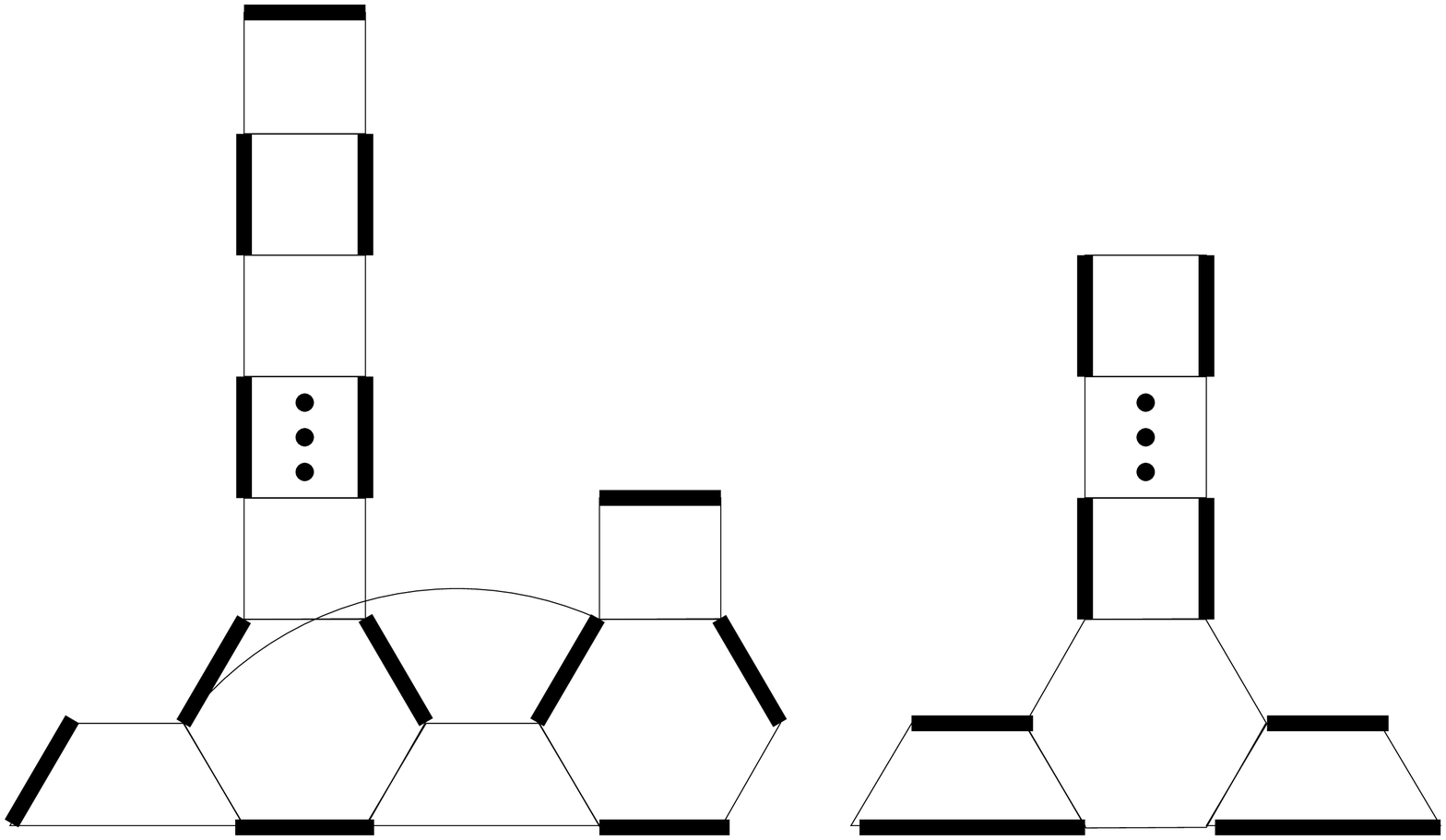}. \hspace{2em}  This pair has precisely the correct weight of $x_1^4x_2^3x_3^3x_4^2x_5^2\cdots x_{2i-1}^2x_{2i}x_{2i+1}$.

\vspace{1em} \noindent\bf Proof of Lemma \ref{restofdiagonals}.\rm \hspace{0.05in}  The first
part is proven by the following observations.  If we

\vspace{1em} \noindent let graph $G_1$ be \includegraphics[width = 1in , height =
0.7in]{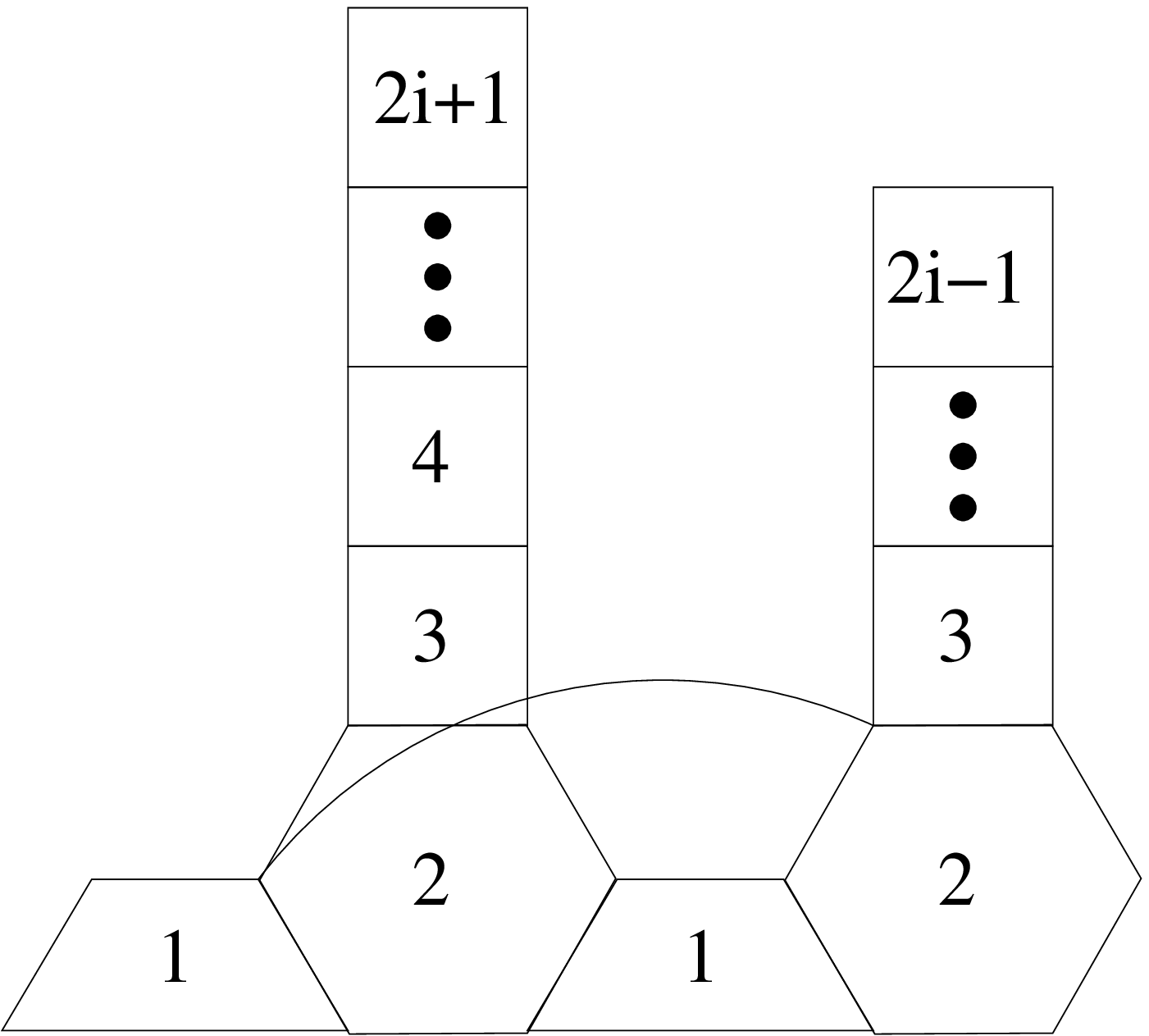}, $G_2$ be
\includegraphics[width = 1in , height = 0.7in]{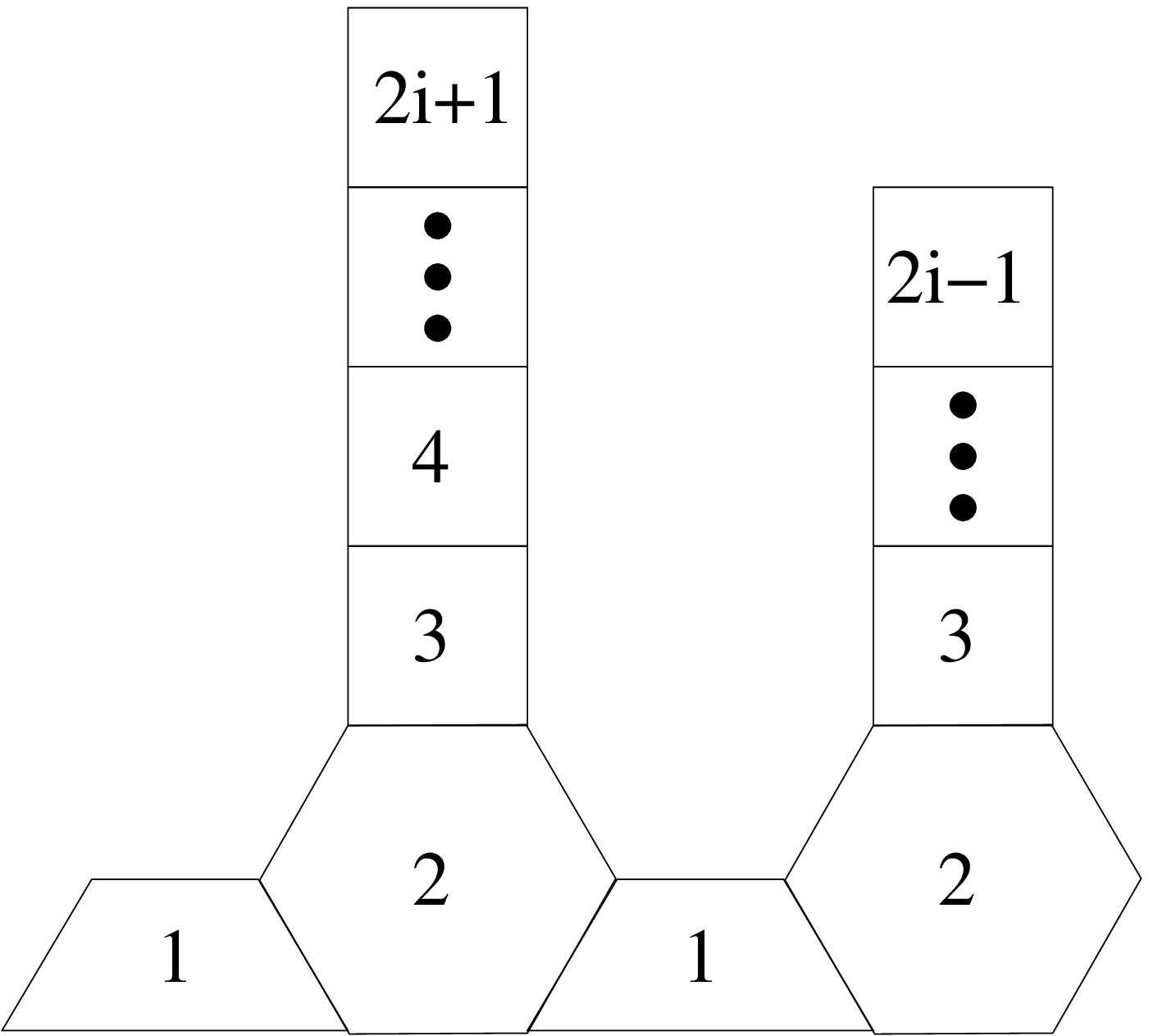},
$T_1$ be \includegraphics[width = 0.7in , height = 0.7in]{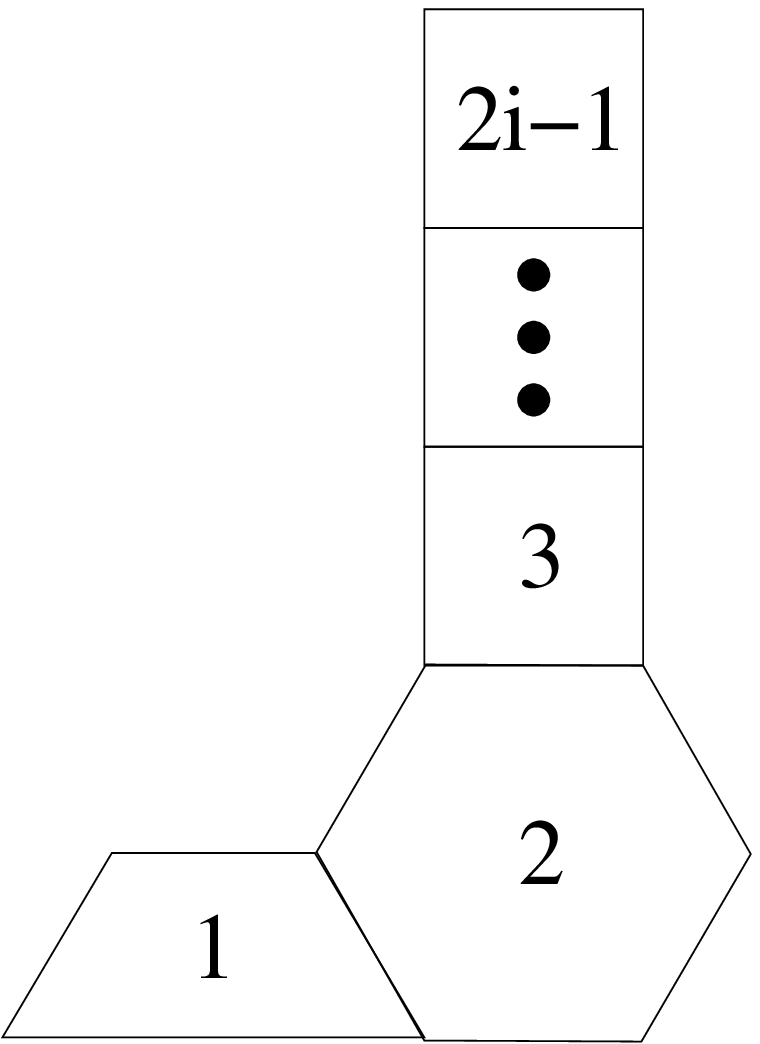}, \hspace{1em}

\noindent $T_2$ be \includegraphics[width = 0.7in , height = 0.7in]{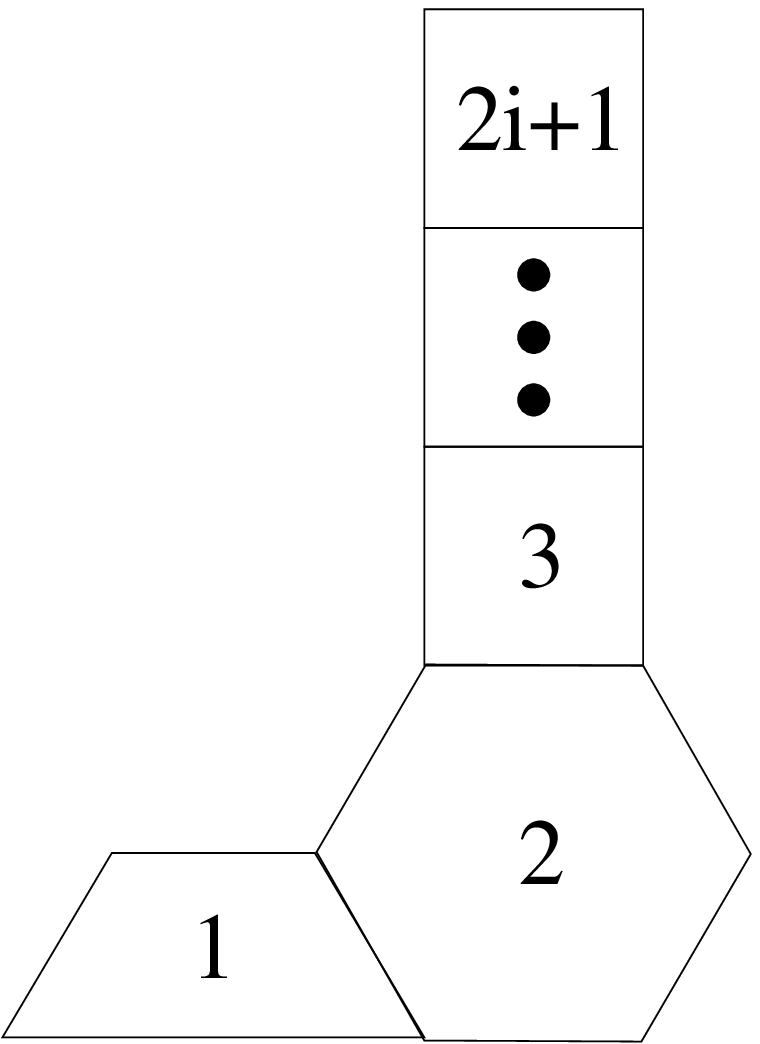}, $H_1$ be
\includegraphics[width = 1.3in , height = 0.65in]{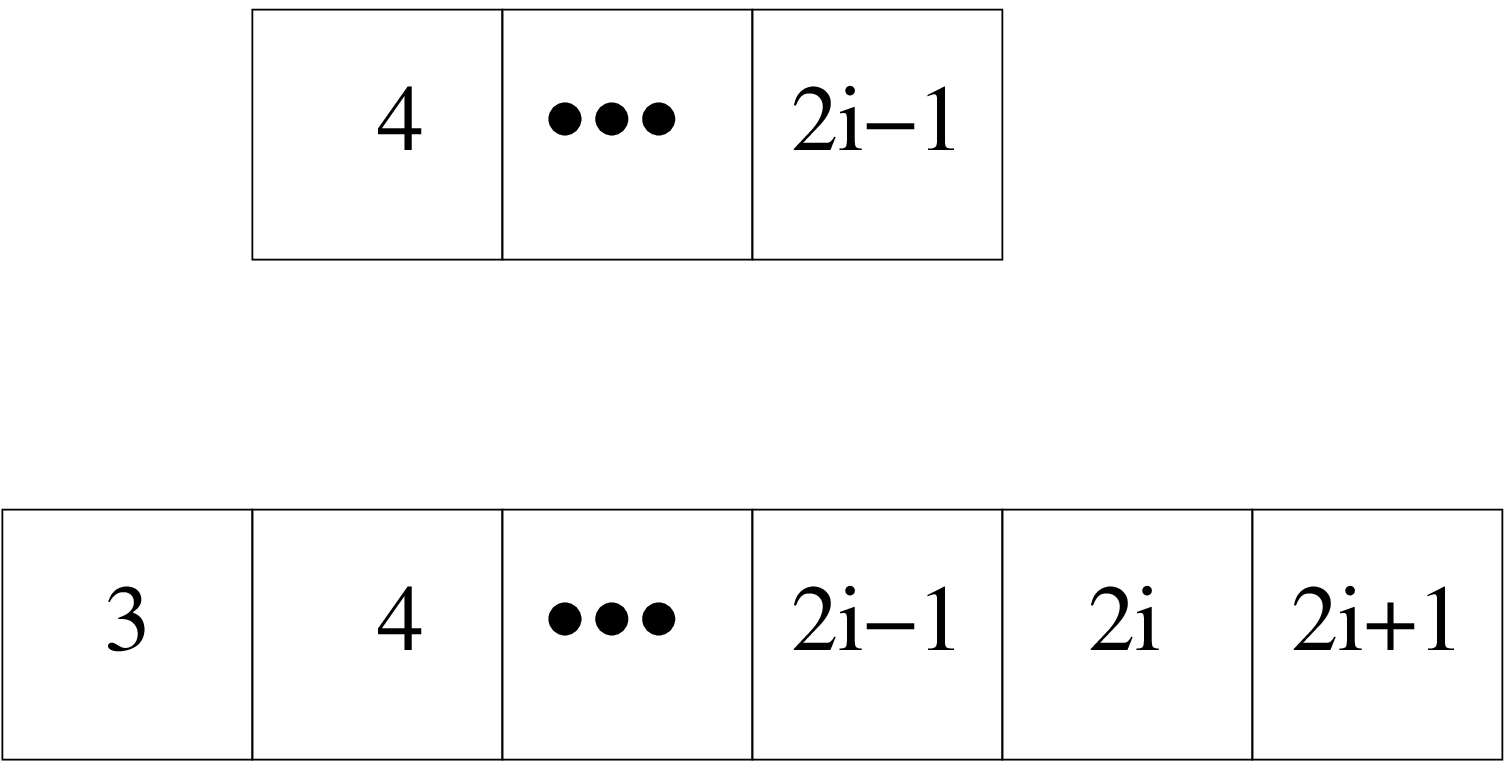},
and let $H_2$ be \includegraphics[width = 1.3in , height = 0.65in]{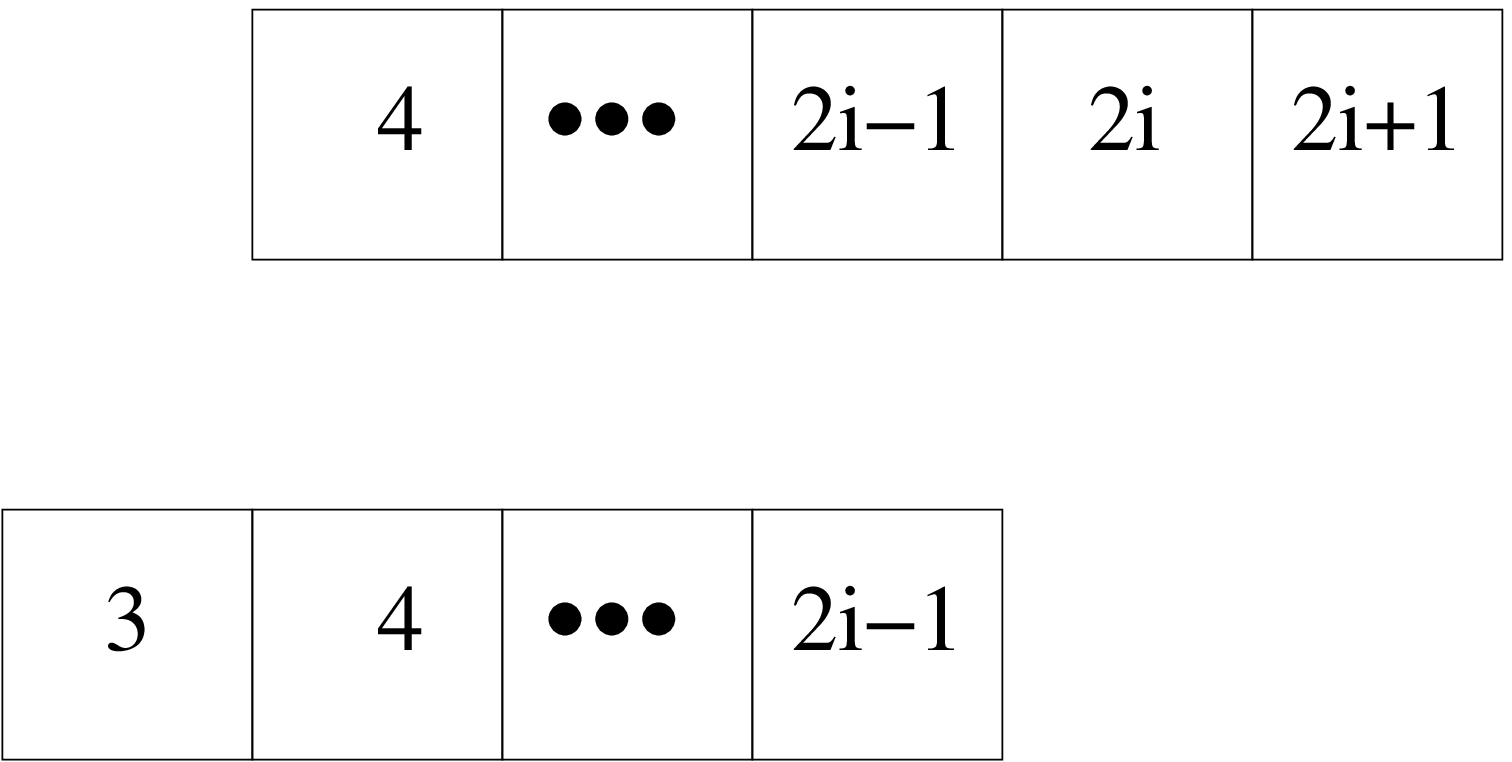} \hspace{0.1in}
then \begin{eqnarray*} P(G_1) &=& P(G_2) + x_1^2x_2x_3P(H_1) \\
P(T_1)P(T_2) &=& P(G_2) + x_1^2x_2x_3P(H_2). \end{eqnarray*}

Putting these two  equalities together we obtain
$$P(G_1) = P(T_1)P(T_2) + x_1^2x_2x_3\bigg(P(H_2) -P(H_1)\bigg).$$

\noindent Most matchings of $H_2$ correspond to a matching of $H_1$ by the usual procedure of swapping the right-hand sides.  The extraneous
matching of $H_2$ has the form

\begin{center}\includegraphics[width = 1.3in , height = 0.65in]{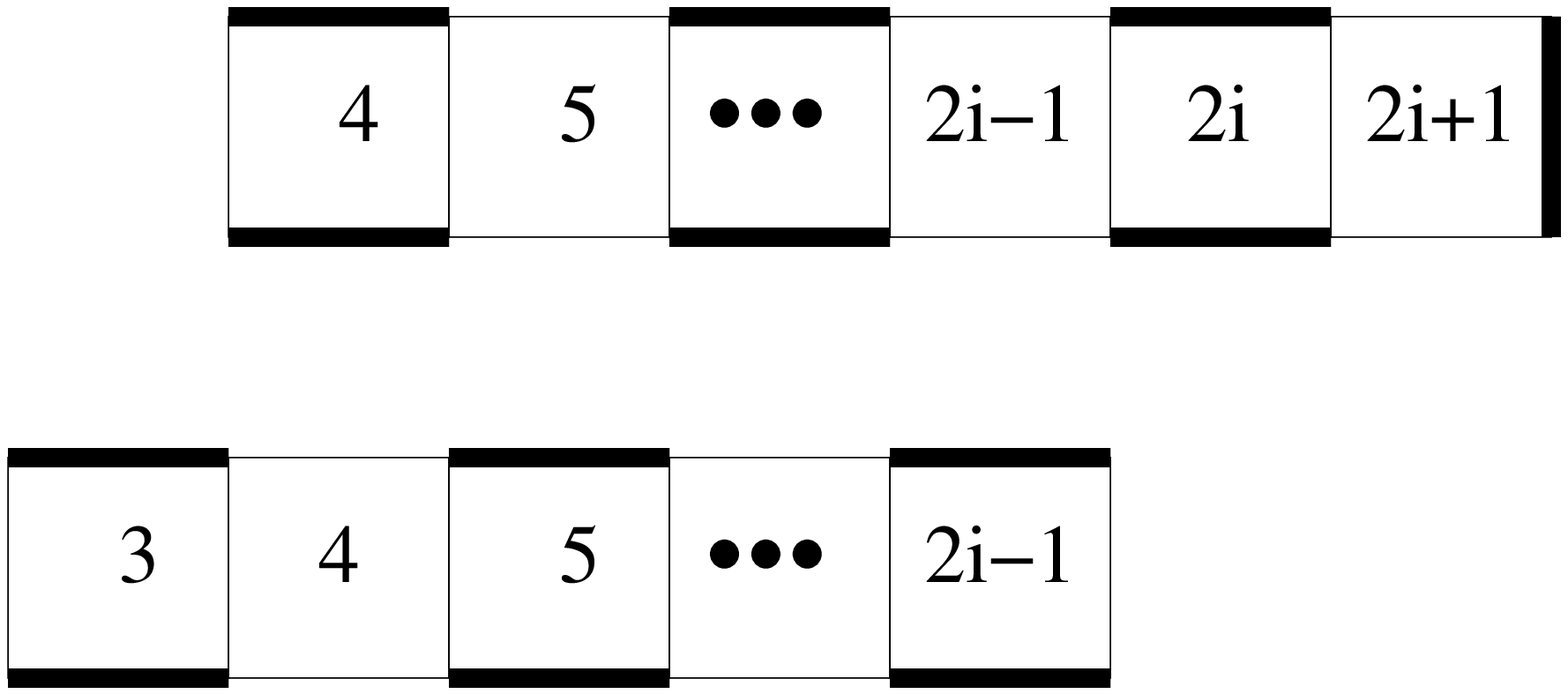}, \end{center} 
contributing a factor of $x_3x_5x_5x_7\cdots x_{2i-1}x_{2i+1}\cdot x_2x_4x_4x_6\cdots
x_{2i-2}x_{2i}$, and yielding the identity $$P(G) = P(T_1)P(T_2) = x_1^2x_2^2x_3^2\cdots x_{2i-1}^2x_{2i}x_{2i+1}.$$
Since $x_1^{(j-1)}x_1^{(j+1)} = x_2^{(j)}+1$ is satisfied by letting $T_1\longleftrightarrow x_1^{(j-1)}$ and $T_2\longleftrightarrow x_1^{(j+1)}$, part one of Lemma \ref{restofdiagonals} is proved.  

Part three is proved analogously to Lemma \ref{xii+2}.  In this case, we have a diamond where all four entries are graphs consisting of two towers on the maximal base of two trapezoids and two hexagons.  We inductively know the validity of these graphs for Laurent polynomials $a$, $b$, and $c$, so it is sufficient to verify the diamond condition if $d \longleftrightarrow$  \includegraphics[width = 1.5in , height = 1in]{Bxii2j.eps}.  For ease of notation, we temporarily let $G_a$, $G_b$, $G_c$, and $G_d$ be the graphs corresponding to these particular Laurent polynomials.  As before, we wish to present a bijection 
$$\phi: M(G_a) \times M(G_d) \setminus\{(m_a^\prime, m_d^\prime)\} \rightarrow M(G_b)\times M(G_c)$$ where $(m_a^\prime, m_d^\prime)$ is a specific pair of matchings.  Map $\phi$ starts by swapping the two left towers of $G_a$ and $G_d$ if able; this is analogous to the earlier cases.  However, in the case where these towers cannot be swapped (because the alternating pattern continues down into the base) map $\phi$ then attempts to swap the right towers of $G_a$ and $G_d$.  There is exactly one pair of matchings where both attempts at swapping fails.  This is the pair of matchings where the alternating patterns appear on both towers down through the bases; by inspection, such a pair has precisely the weight of the extraneous monomial.

This leaves part two as the crux of the proof and the step
which utilizes the diamond relation $ad - b^2c = 1$ which makes $B_\infty$ different from the
previous cases.

We wish to show that the assignments \vspace{1em}

$a ~\longleftrightarrow~ $ \includegraphics[width = 1in , height = 1in]{Lem8G1.eps}
\hspace{0.1in}, \hspace{0.6in} $d ~\longleftrightarrow~ $
\includegraphics[width = 1.5in , height = 1in]{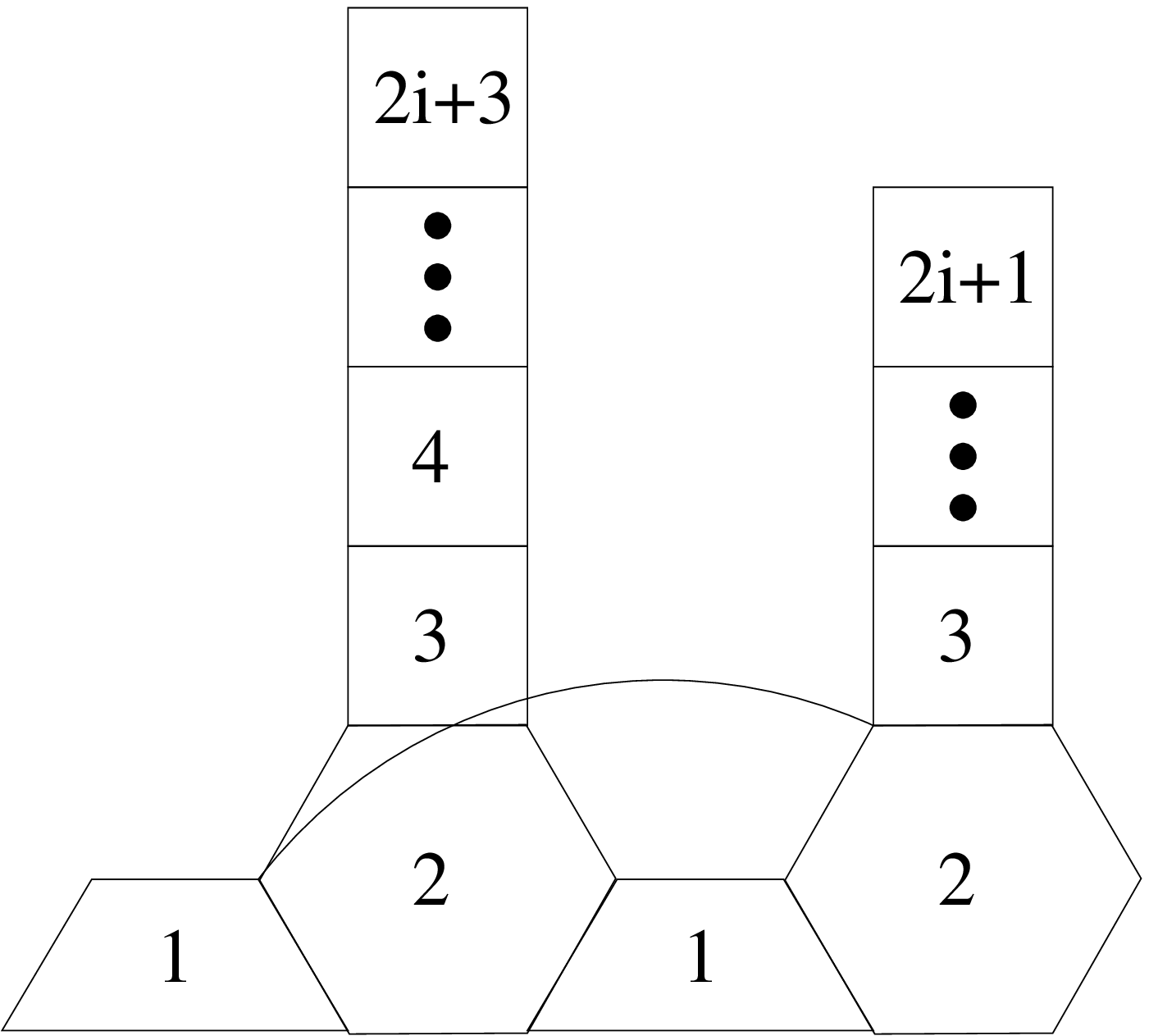} \\ \vspace{1em}

$b ~\longleftrightarrow~ $\includegraphics[width = 1in , height = 1in]{Lem8T2.eps} \hspace{0.1in},
\hspace{0.6in} $c ~\longleftrightarrow~ $ \includegraphics[width = 1in , height =
1in]{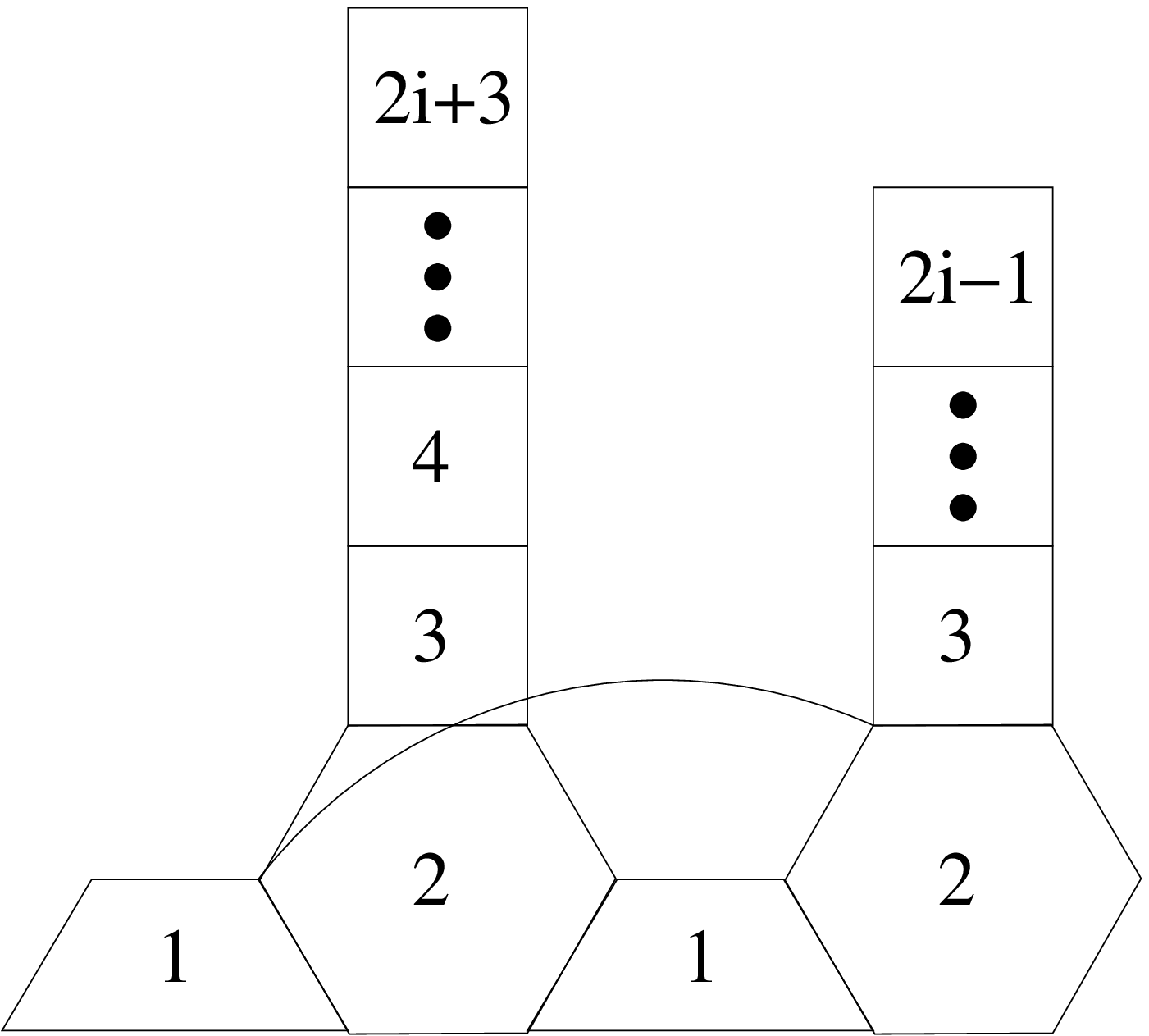}\hspace{0.1in} 

\vspace{2em}\noindent satisfies $ad-b^2c=1$.  We shall use the fact that if $b ~\longleftrightarrow~ $\includegraphics[width = 1in
, height = 1in]{Lem8T2.eps} \hspace{0.1in}, then $b^2 ~\longleftrightarrow~
$\includegraphics[width = 1.5in , height = 1in]{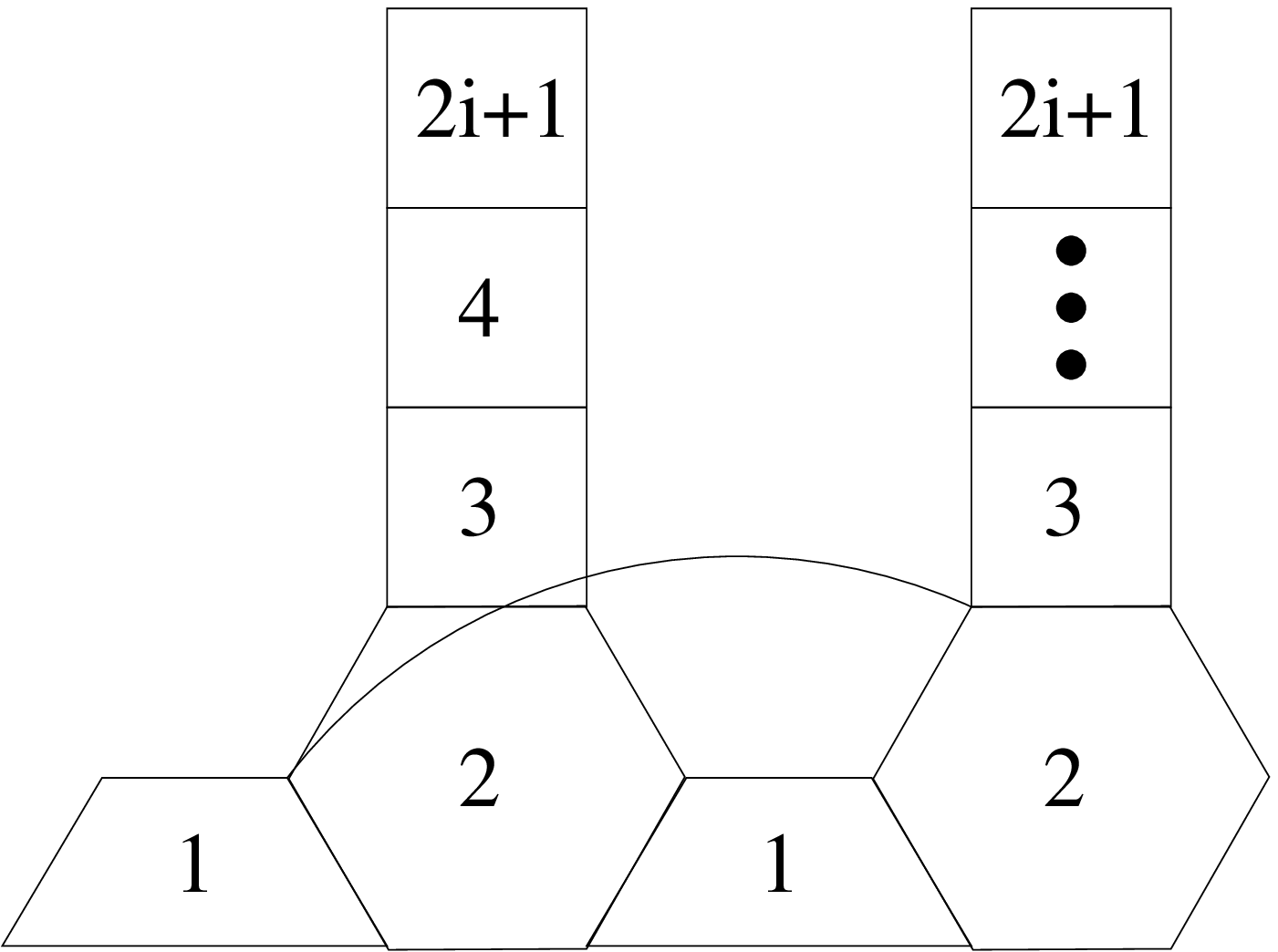}.

\noindent Clearly the denominator corresponds correctly.  The number of weighted matchings, and
thus the numerator, is also correct since there is a weight-preserving bijection between matchings of
\vspace{1em}$\includegraphics[width = 1.05in , height = 0.7in]{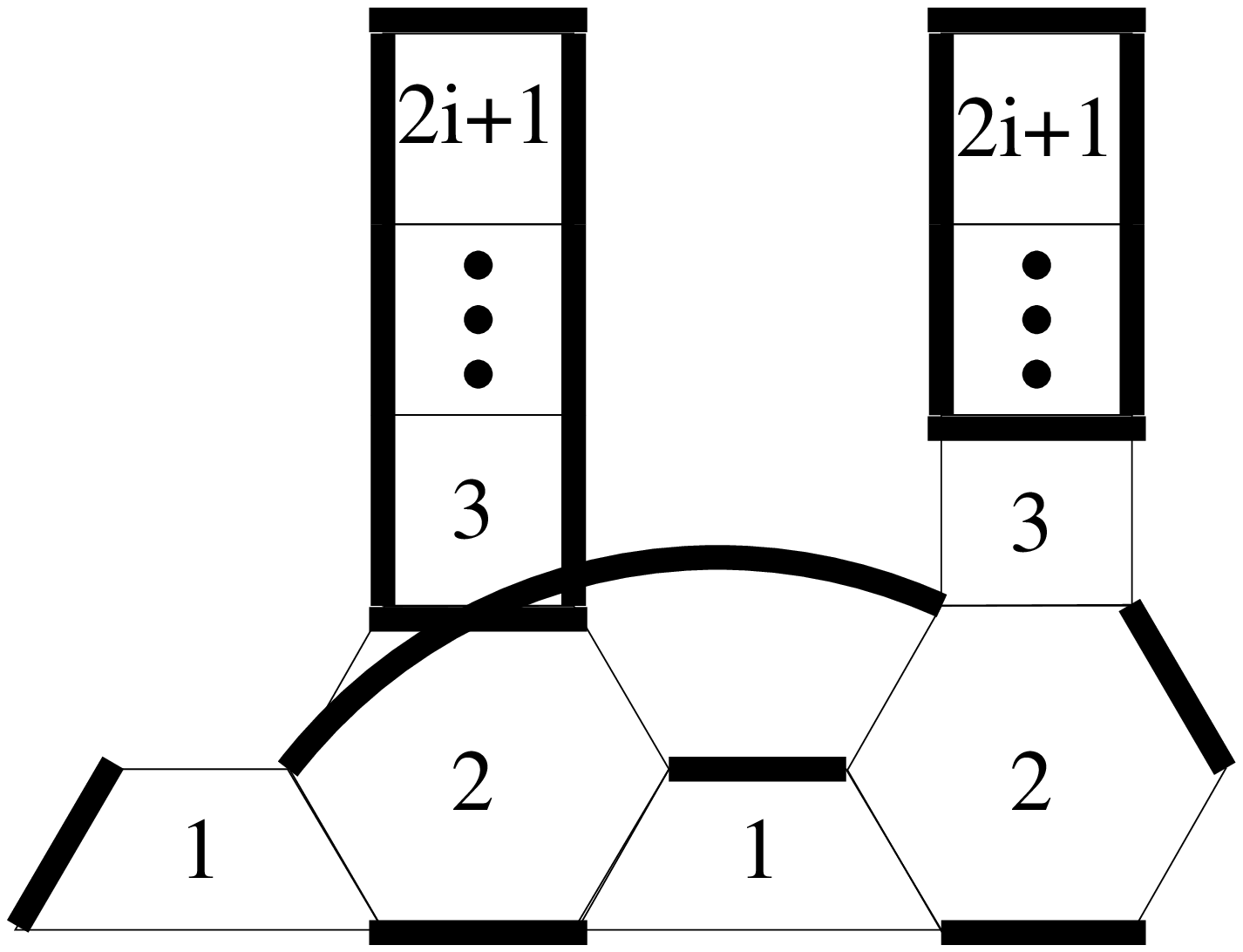} \hspace{0.1in}$ and matchings of 
$\includegraphics[width = 1.05in , height = 0.7in]{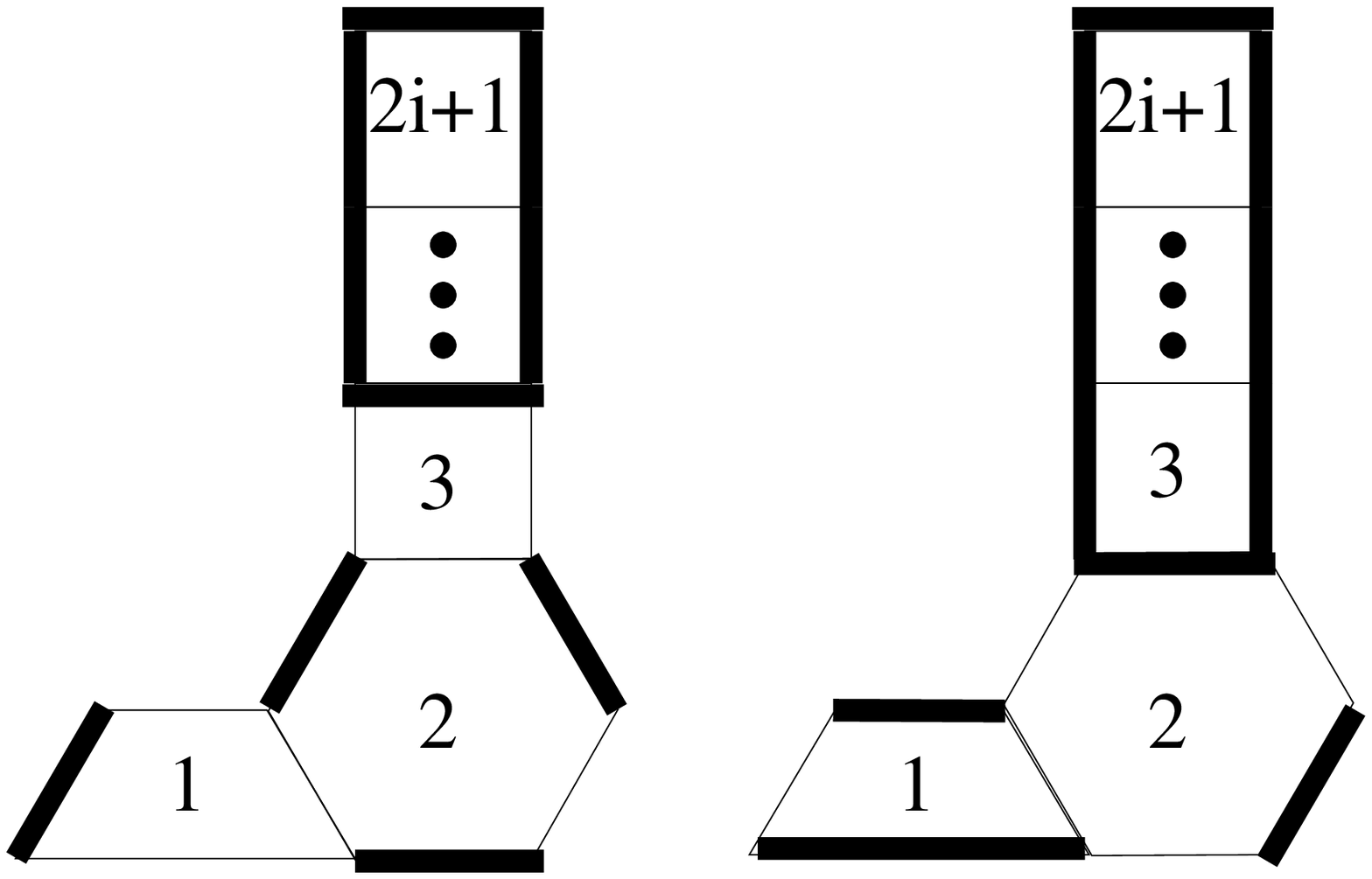}$. With this substitution, a superposition argument analogous to that which just proved part three demonstrates the validity of the standard diamond relation $ad - (b^2)c = 1$.

With the last step completed, these three Lemmas prove that the cluster variables of the
$B_\infty$-lattice correspond exactly to the desired graphs.  There is a pattern inherent in the
NW to SE and NE to SW diagonals once again.  This time, this pattern manifests itself (in the
region where $i \leq j$) by dictating the choice of right tower (NW to SE) and the choice of
left tower (NE to SW).  Recall that Lemma \ref{AAregion} already described the pattern in the
region where $i > j$ using grid graphs consisting of tiles $T_a\cup\dots \cup T_b$ where $a \geq 3$.  We now turn to the problem of restricting to specific $B_n$.  We use the
same methodology as in the $A_n$ case.

Given that we did not include the right-hand boundary, graphs can contain tile $T_i$ for arbitrarily large $i$.  We thus want to essentially let $x_{n+1}=1$ to force tile $T_n$ to have the proper weights, as a tile in $\mathcal{T}_{B_n}$ as opposed to $\mathcal{T}_{B_\infty}$.  However, unlike the $A_n$ case, we cannot simply apply this substitution and use horizontal periodicity to make sure the diamond condition holds throughout.  The problem is that the left-hand boundary satisfies a diamond conditon of a different form.  Nonetheless, the diagonals to the northeast, those where $i - j \geq
2$, contain graphs which are connected grid graphs, towers, of the form $T_a \cup \dots \cup T_b$ where $3 \leq a \leq b$, and any neighboring four entries satisfy the same diamond condition as the $A_n$ case.  Thus the logic of Lemma \ref{CenterOne} carries over and we can excise connected subgraphs centered at tile $T_{n+1}$ ($\tilde{T}_{n+2}$ under the old notation).  In particular, we obtain the following excision.

\begin{Lem} If $b$ satisfies $n+1 \leq b \leq 2n-(a-2)$ then the graphs
$T_a\cup \dots \cup T_b$ and $T_a \cup \dots \cup T_{2n+2-b}$ biject to the same Laurent polynomial.
\end{Lem}

This Lemma follows directly from Lemma \ref{CenterOne} after replacing $T_{n+1}$ by the equivalent tile $\tilde{T}_{n+2}$.  The restriction of $b \leq 2n-(a-2)$ must be added here since we have a boundary on the left side, i.e. hexagon $T_2$ cannot be excised during this procedure.  Notice that at the extreme, $T_3 \cup \dots \cup T_{2n-2}$ bijects to the same Laurent polynomial as $T_3$, and further $T_3 \cup \dots \cup T_{2n-1}$ is centered around $T_{n+1}$ and thus bijects to the Laurent polynomial $1$, the same as the empty graph.  

Using this Lemma, we are able to determine the northeast corner of the finite $B_n$ lattice.  Using these graphs during the inductive step of Lemma \ref{xii} in lieau of the arbitrarily large towers of $T_3\cup \dots \cup T_{2j-1}$ allows us to fill in the next diagonal of the $B_n$ lattice where the towers sitting on the base of $T_1 \cup T_2 \cup T_1$ will consist exclusively of the tiles between $T_3$ and $T_n$.  Simiarly, the recursive steps of Lemmas \ref{xii+2} and \ref{restofdiagonals} also follow with these truncated graphs, with tiles between $T_1$ and $T_n$ used instead.

Applications of Lemma \ref{xii+2} and then successive applications of the third part of Lemma
\ref{restofdiagonals} will allow this interpretation to extend down SW to all but a SE corner of
the lattice.  Note that the diagonals again determine the left-hand and right-hand towers of the
$B_n$-lattice since this property is inherited from the $B_\infty$-lattice.  We compute the SE corner by starting at the bottom with initial row $T_2$, $T_4$, $T_6$, $\dots$, and propogating \emph{upwards} via the diamond condition.  We get via Lemma \ref{AAregion}, now with tile $T_2$, instead of $T_3$, as the smallest allowable tile, all but a single diagonal of the SE corner.  The final diagonal has the form $(Tower~1)\cup T_2 \cup T_1 \cup T_2 \cup T_1$, proven by applying Lemma \ref{xii+2} upwards.   

Notice that in the end, we obtain a lattice where the NW to SE diagonals dictate the right towers and NE to SW diagonals dictate the left towers.  There is one caveat: the \emph{empty} tower, $Tow_\emptyset$ is now allowed.  Thus one has to determine from context whether a graph consisting of a single tower is of the form $Tow_\emptyset \cup Tow_R$ or $Tow_L \cup Tow_\emptyset$.  Alternatively, we can picture the SE corner as sitting directly above the NE corner of the lattice to form a half-diamond.  Thus Proposition \ref{CaseBn} is proven.
\end{proof}

On the next page, we give the lattice corresponding to $B_6$.  Notice that there are six graphs in the northeast corner and four graphs in the southeast corner which are also graphs corresponding to positive roots and cluster variables for the $A_6$ case.  In fact, if we decrease each label by $2$ and horizontally reflect the southeast corner, we can fit these two pieces together to obtain the $A_4$ lattice of graphs exactly.

\vspace{1em}

Also comparing with the $B_4$ lattice we notice boundary behavior.  For example the second entry of the third column is now $T_3 \cup T_4$ instead of $T_3 \cup T_4\cup T_5$ and the second to last element of the fourth column is $T_2 \cup T_3$ instead of $T_2 \cup T_3 \cup (T_4\cup T_5 \cup T_6)$.

\vspace{10em}

\newpage $\begin{array}{cccccc}
\includegraphics[width = 0.2in , height = 0.2in]{BB1.eps}  &~~~& \includegraphics[width = 0.2in , height = 0.2in]{BB3.eps} &~~~& \includegraphics[width = 0.2in , height = 0.2in]{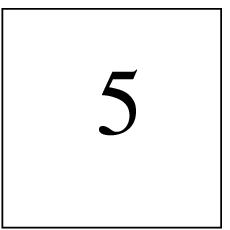} &~~~ \\
 ~~~& \includegraphics[width = 0.6in , height = 0.45in]{BB1213.eps} &~~~& \includegraphics[width = 0.2in , height = 0.5in]{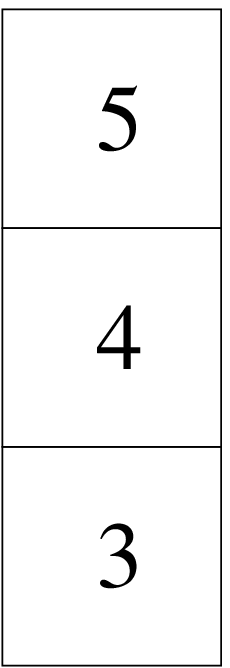} &~~~& \includegraphics[width = 0.2in , height = 0.4in]{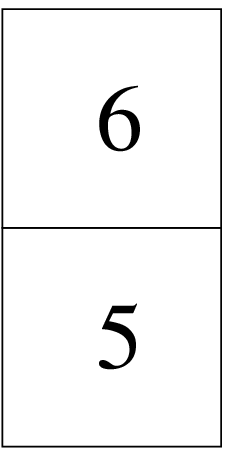}   \\
\includegraphics[width = 0.4in , height = 0.45in]{BB123.eps}  &~~~& \includegraphics[width = 0.6in , height = 0.6in]{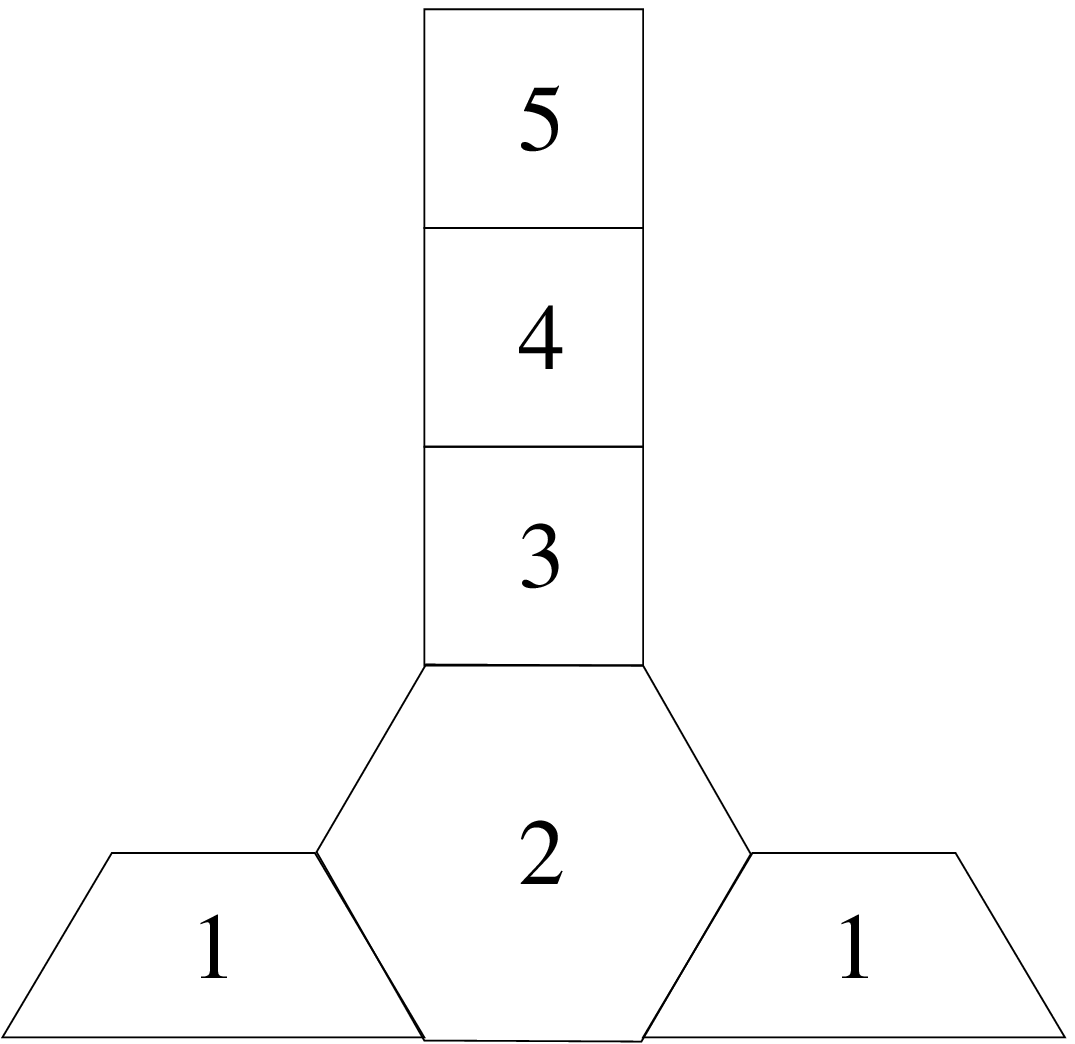} &~~~& \includegraphics[width = 0.2in , height = 0.7in]{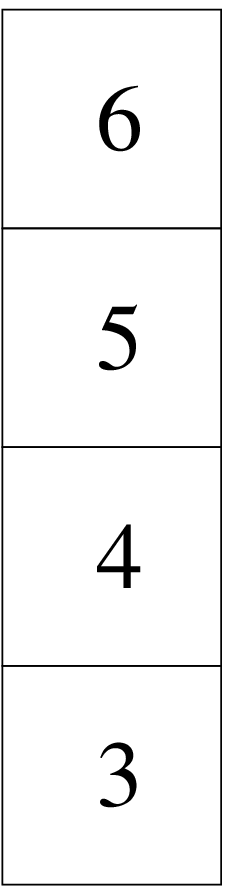} &~~~ \\
 ~~~& \includegraphics[width = 0.8in , height = 0.7in]{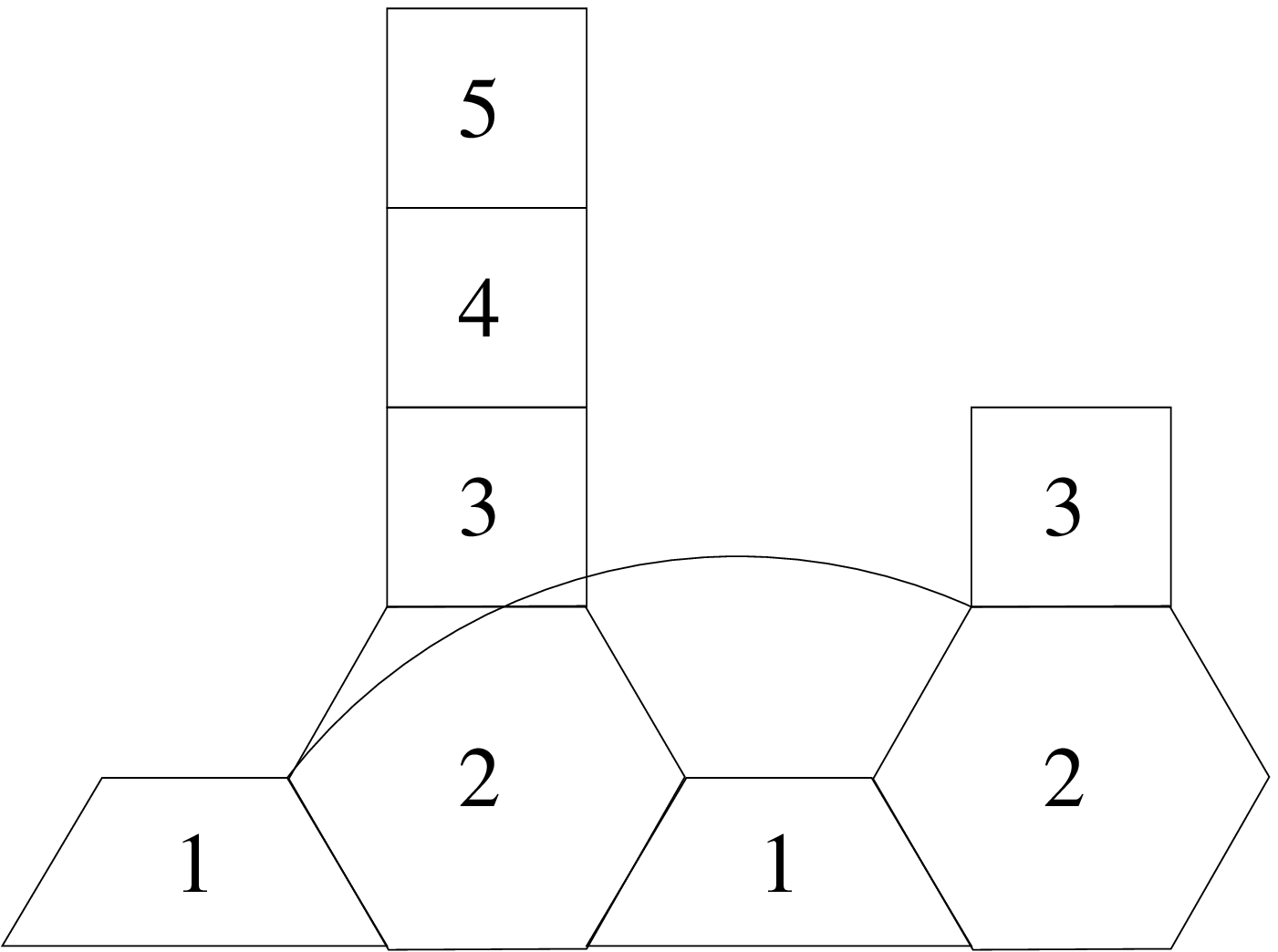} &~~~& \includegraphics[width = 0.6in , height = 0.8in]{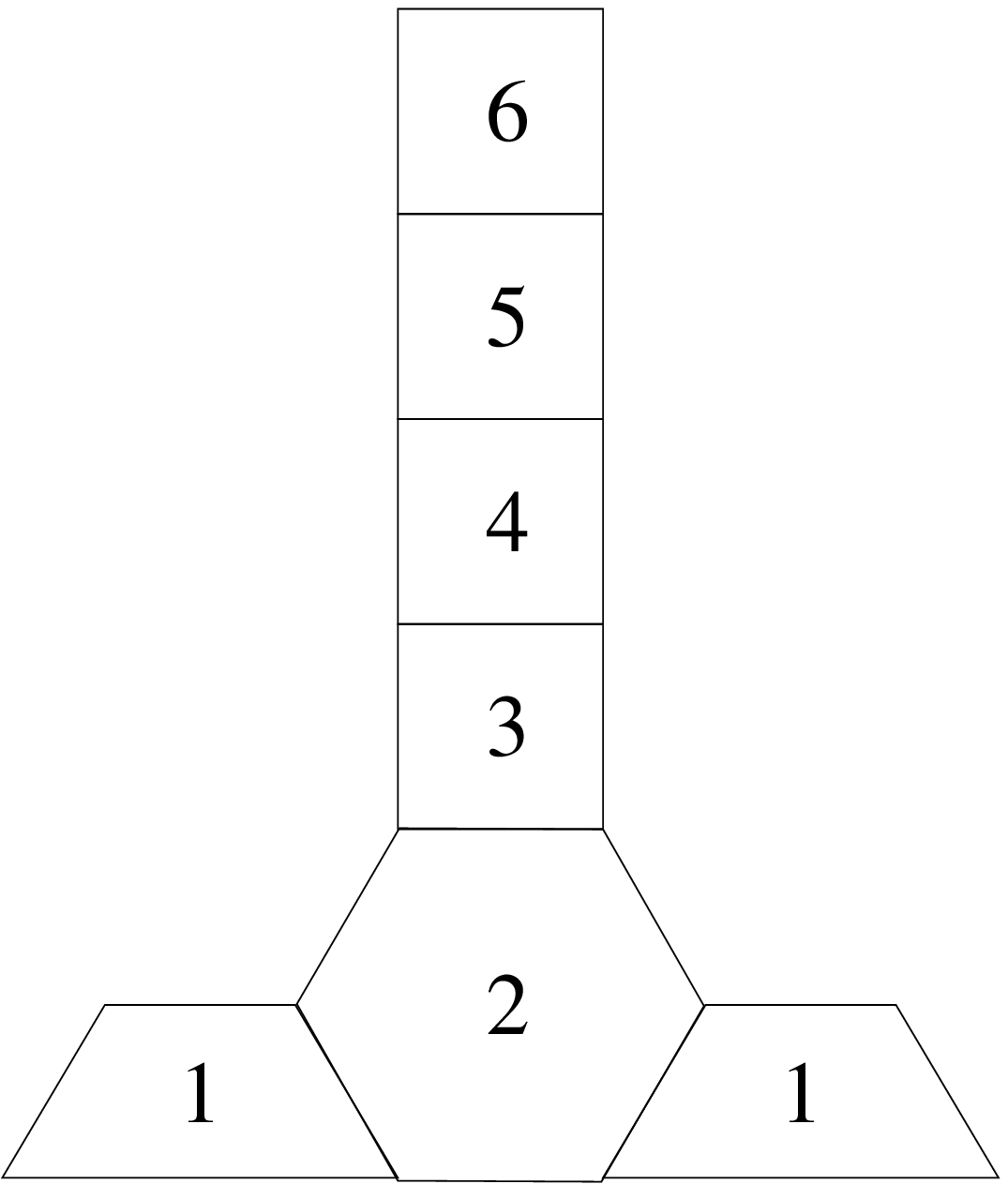} &~~~& \includegraphics[width = 0.2in , height = 0.4in]{BB34.eps}   \\
\includegraphics[width = 0.4in , height = 0.7in]{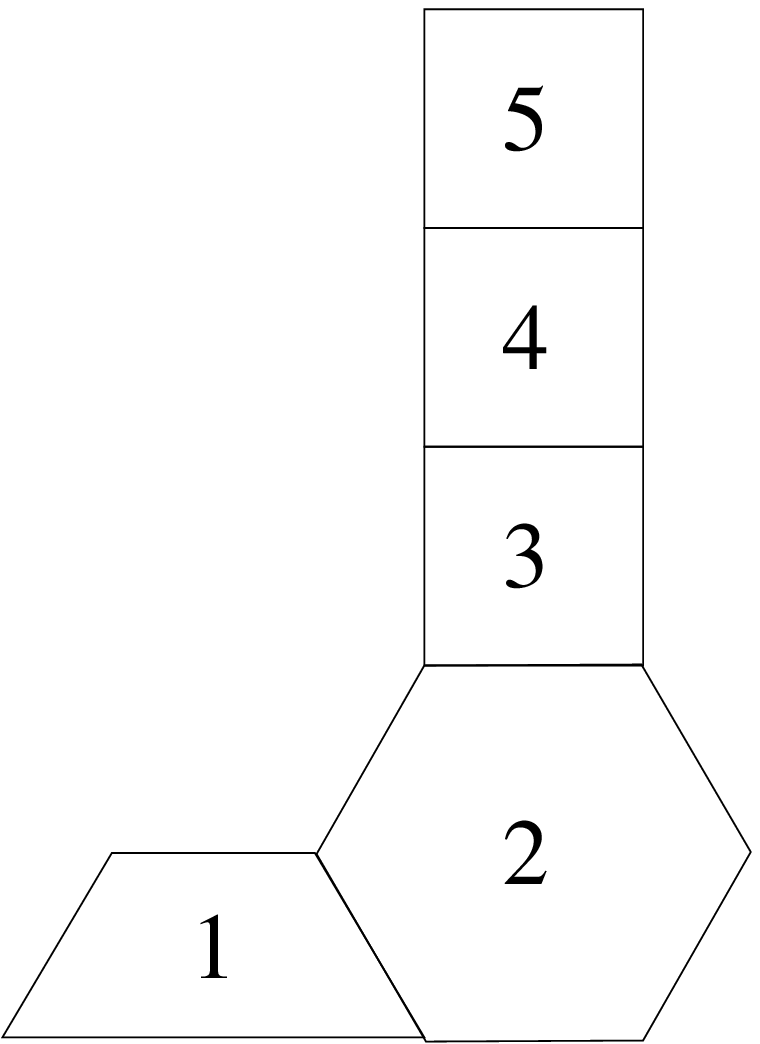}  &~~~& \includegraphics[width = 0.8in , height = 0.8in]{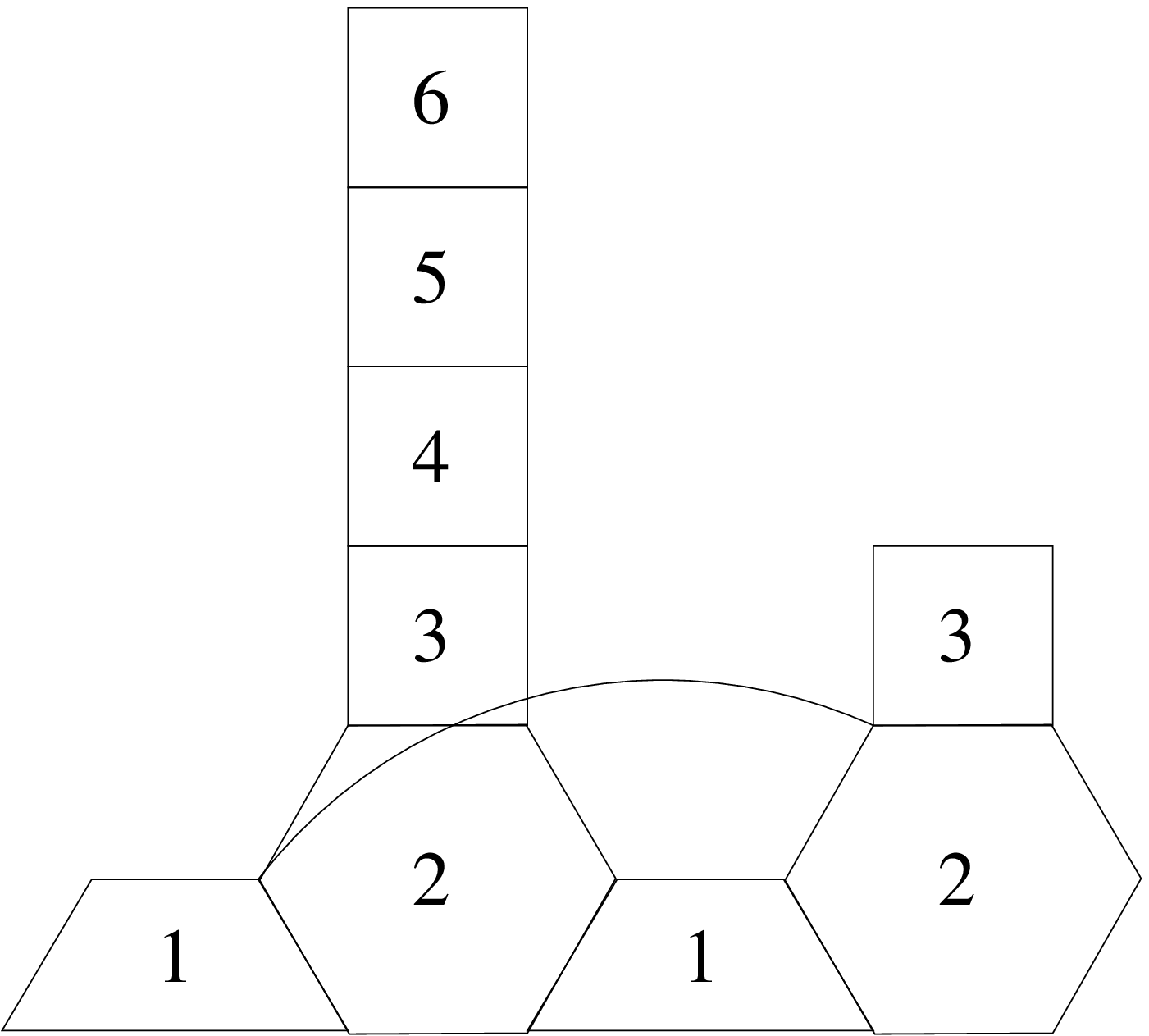} &~~~& \includegraphics[width = 0.6in , height = 0.6in]{BB1214.eps} &~~~ \\
 ~~~& \includegraphics[width = 0.8in , height = 0.8in]{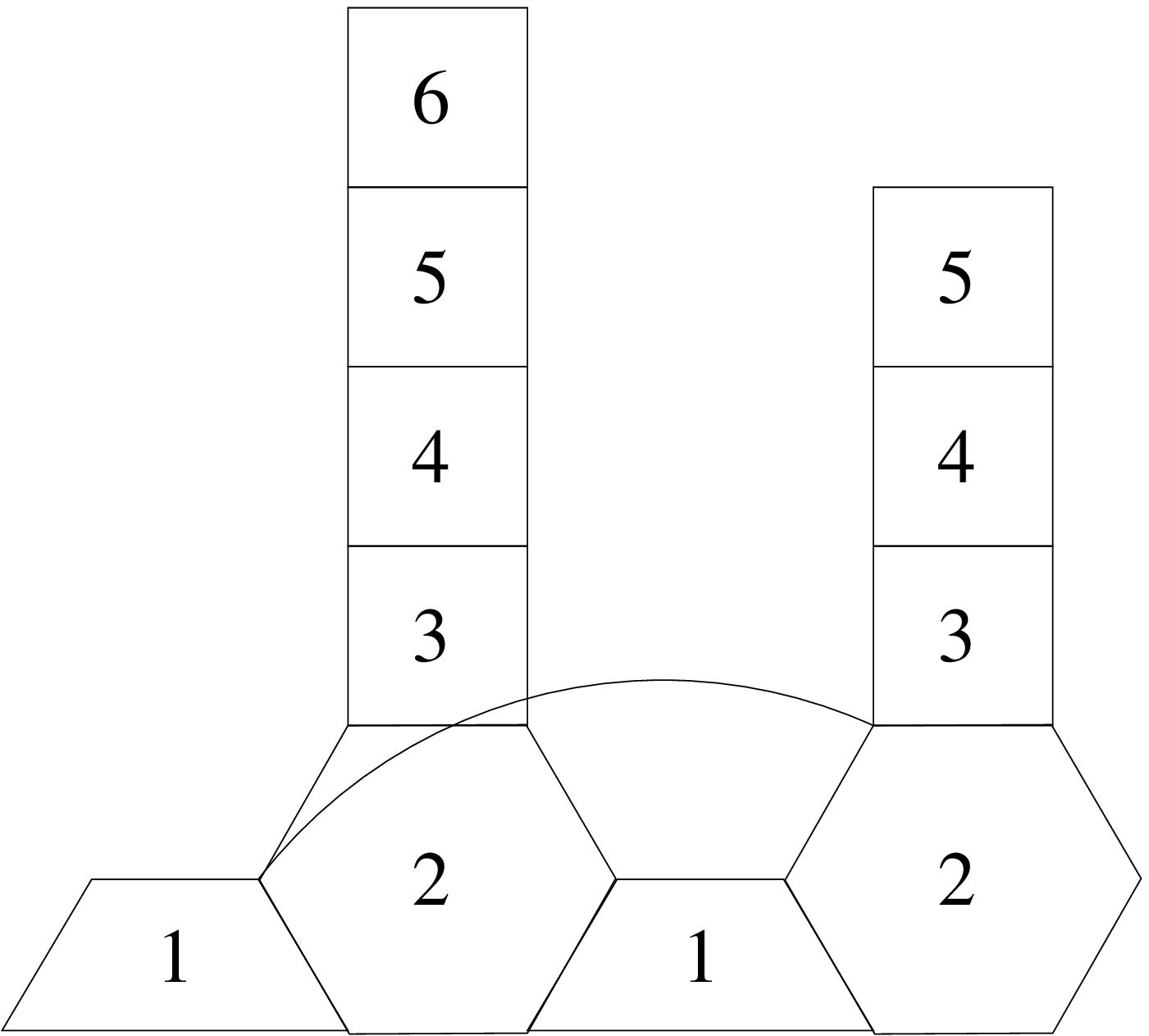} &~~~& \includegraphics[width = 0.8in , height = 0.8in]{BB121234.eps} &~~~& \includegraphics[width = 0.6in , height = 0.4in]{BB121.eps}   \\
 \includegraphics[width = 0.4in , height = 0.8in]{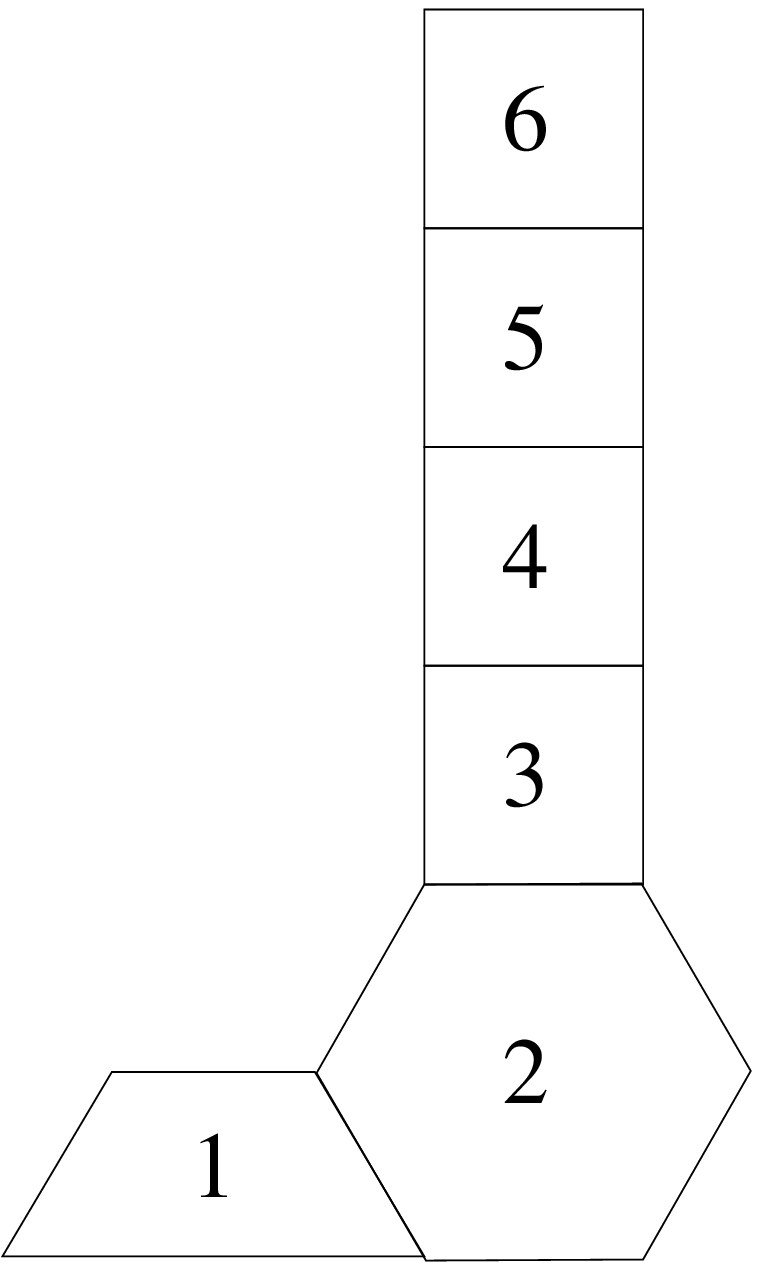}  &~~~& \includegraphics[width = 0.8in , height = 0.7in]{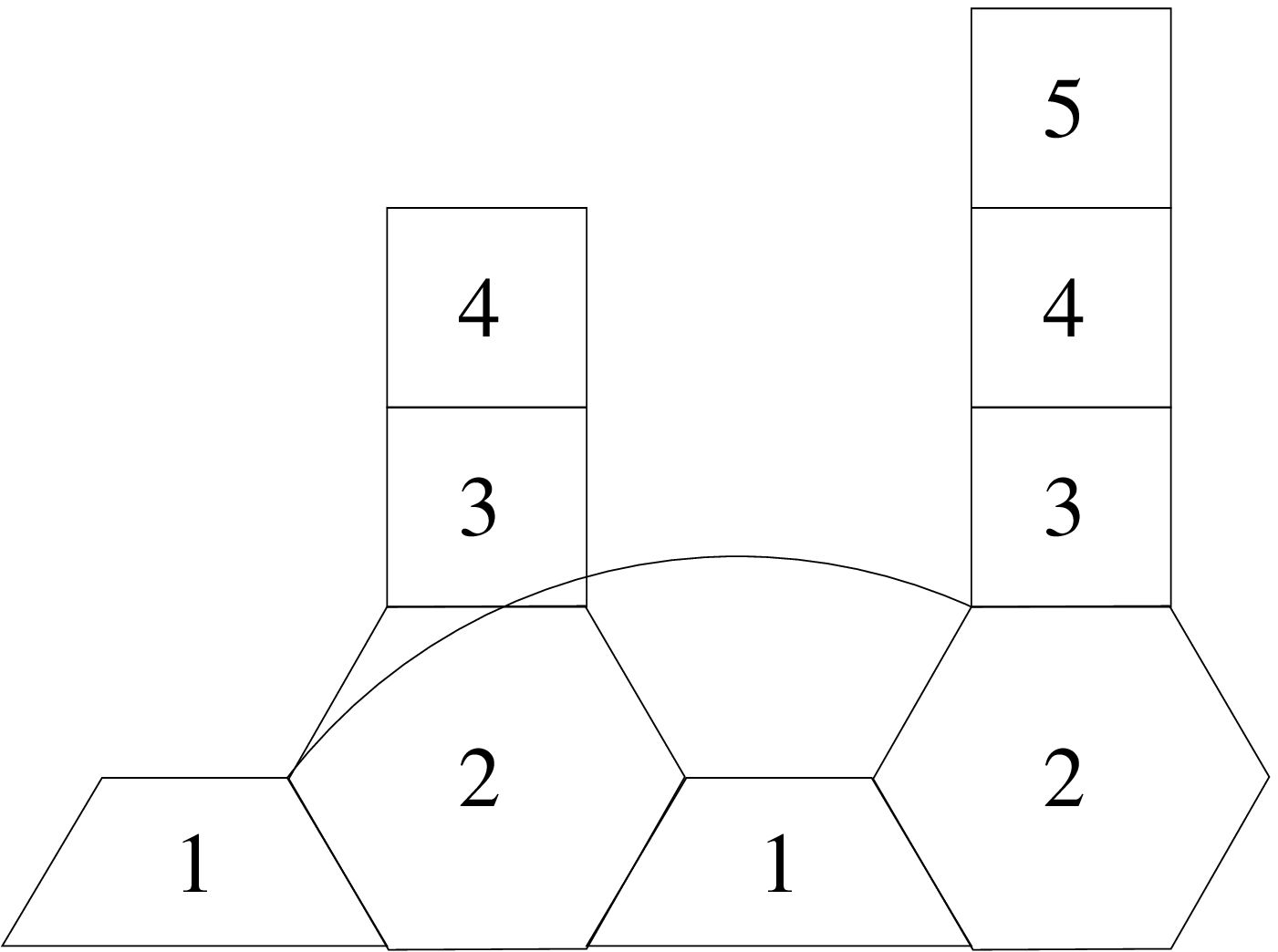} &~~~& \includegraphics[width = 0.8in , height = 0.4in]{BB121230.eps} &~~~ \\
 ~~~& \includegraphics[width = 0.8in , height = 0.8in]{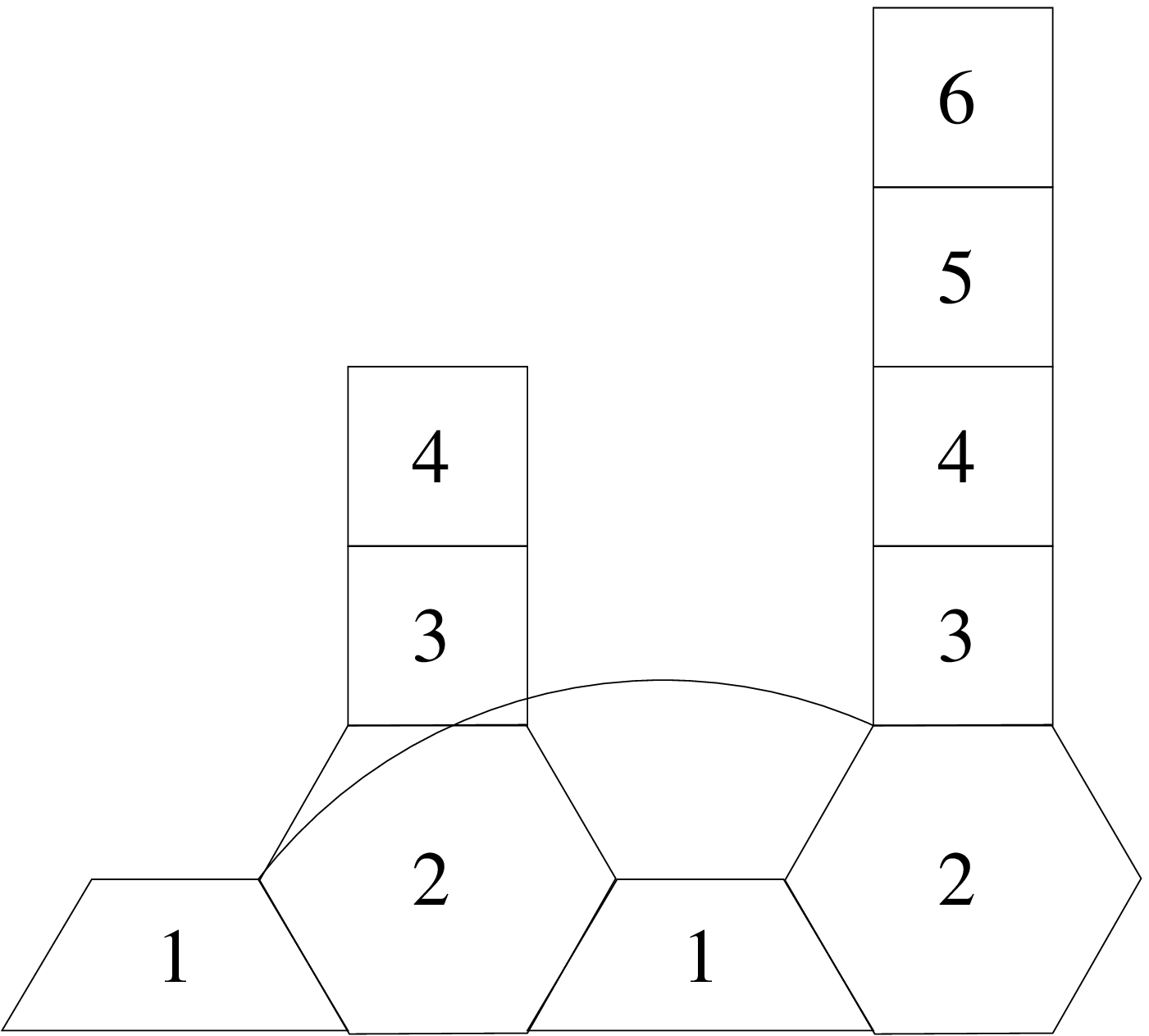} &~~~& \includegraphics[width = 0.8in , height = 0.7in]{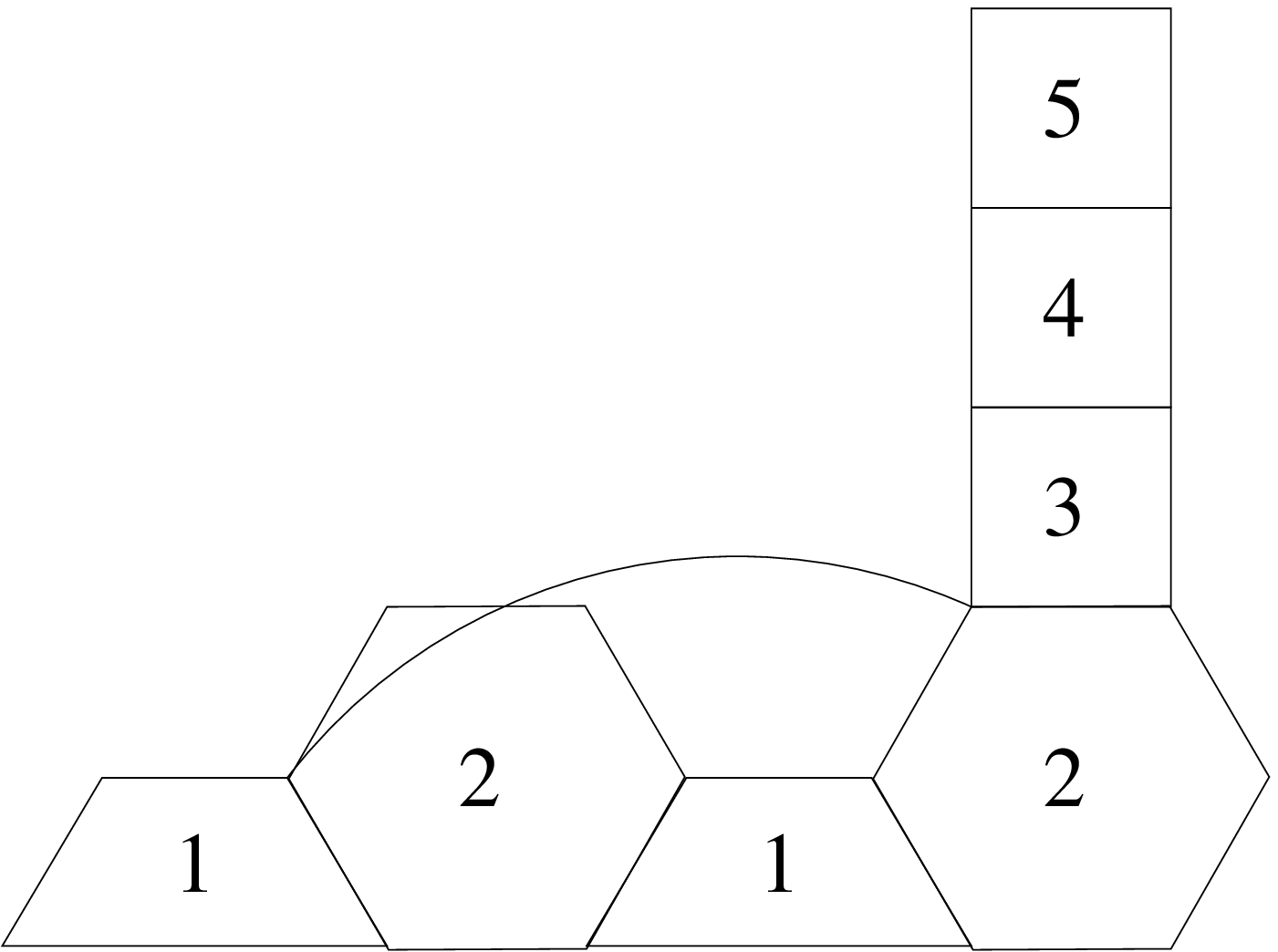} &~~~& \includegraphics[width = 0.2in , height = 0.4in]{BB23.eps}   \\
 \includegraphics[width = 0.4in , height = 0.6in]{BB124.eps}  &~~~& \includegraphics[width = 0.8in , height = 0.8in]{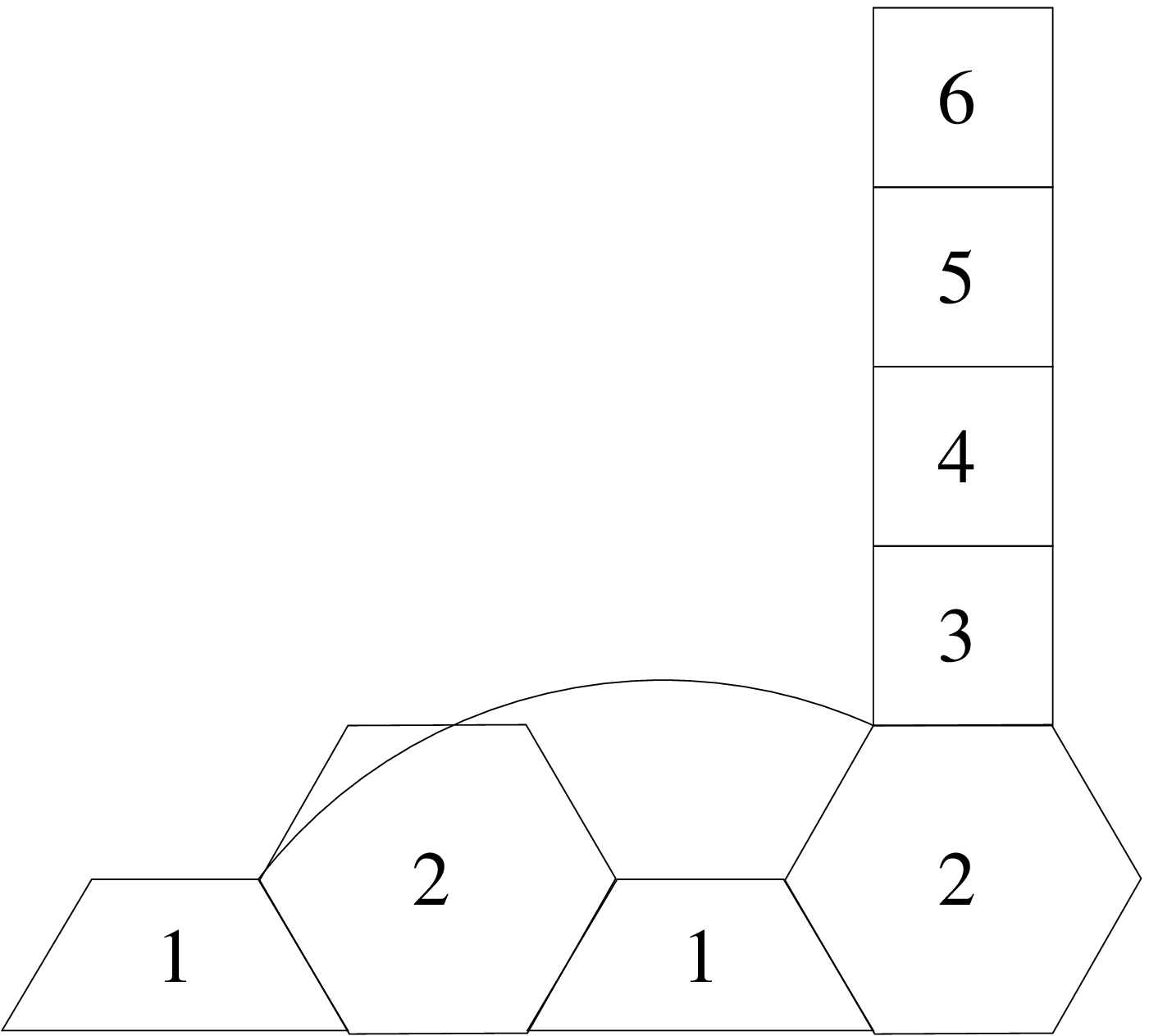} &~~~& \includegraphics[width = 0.2in , height = 0.7in]{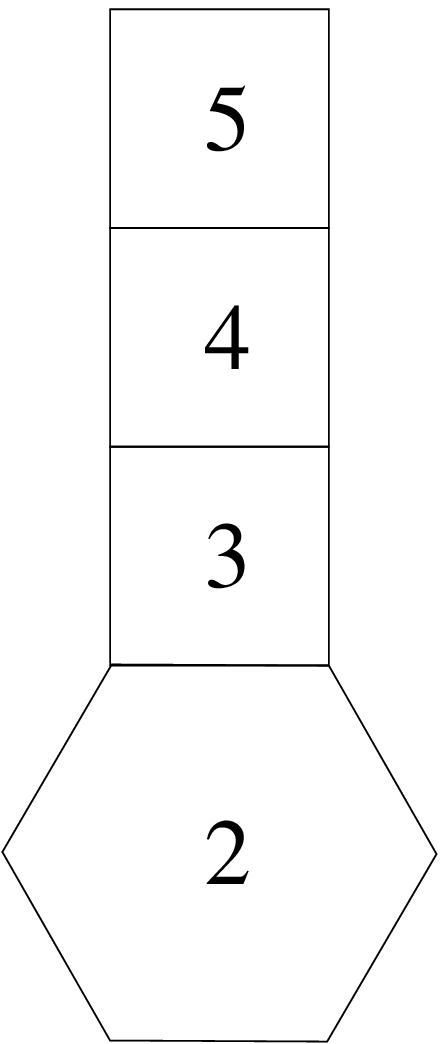} &~~~ \\
 ~~~& \includegraphics[width = 0.8in , height = 0.6in]{BB121240.eps} &~~~& \includegraphics[width = 0.2in , height = 0.8in]{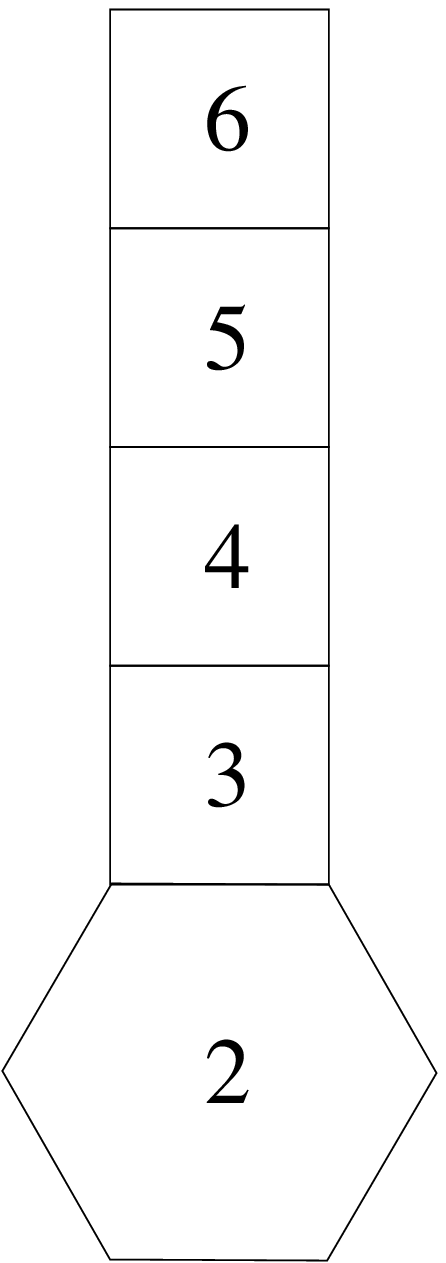} &~~~& \includegraphics[width = 0.2in , height = 0.4in]{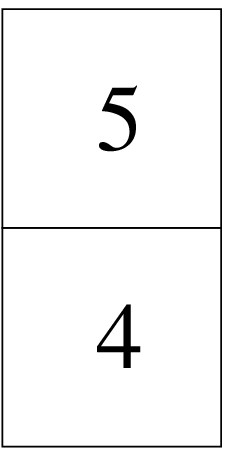}   \\
 \includegraphics[width = 0.4in , height = 0.4in]{BB12.eps}  &~~~& \includegraphics[width = 0.2in , height = 0.6in]{BB24.eps} &~~~& \includegraphics[width = 0.2in , height = 0.5in]{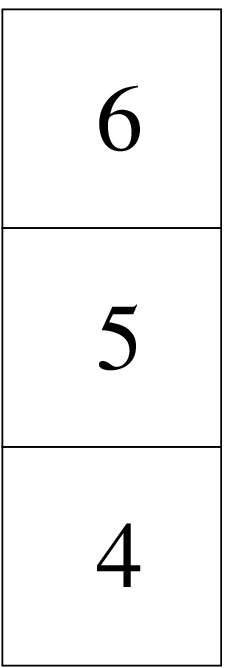} &~~~ \\
 ~~~& \includegraphics[width = 0.25in , height = 0.25in]{BB2.eps} &~~~& \includegraphics[width = 0.2in , height = 0.2in]{BB4.eps} &~~~& \includegraphics[width = 0.2in , height = 0.2in]{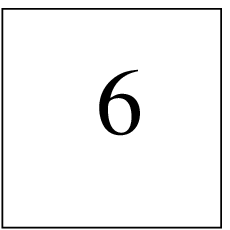}
\end{array}$  
\newpage

We were able to analyze $C_n$ based on $A_{2n-1}$ using a folding
procedure.  Analogously we can analyze $D_n$ using $B_{n-1}$ and
an \emph{unfolding} procedure.  We label the Dynkin diagram for $D_n$ starting with $1$ and $\overline{1}$ on the left, and label the rest in a line from $2$ to $n-1$. 

\begin{center}\includegraphics[width = 2.4in , height = 0.8in]{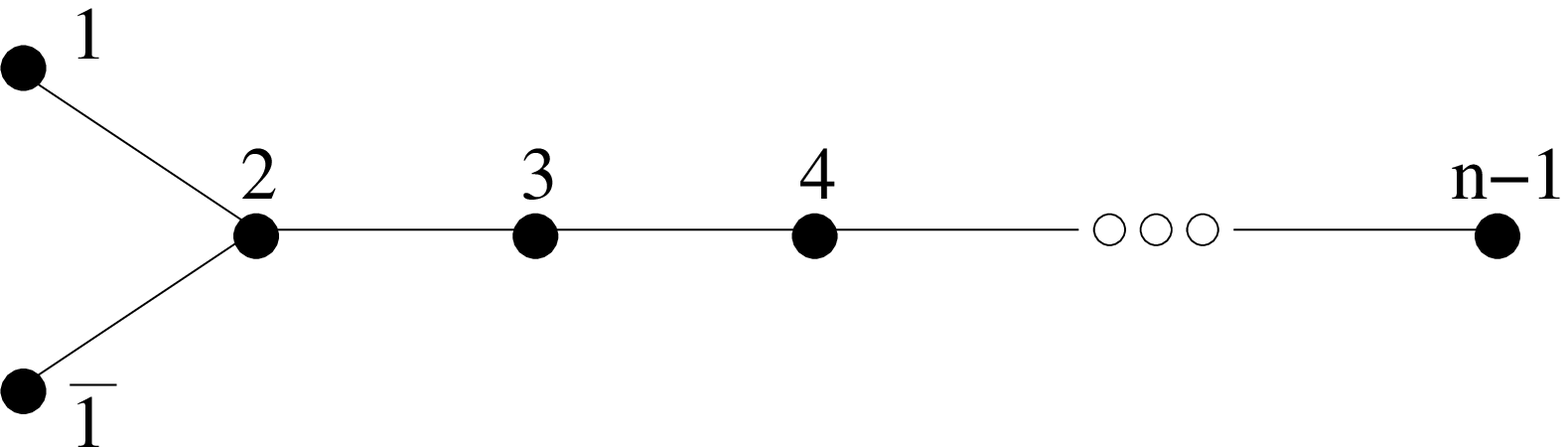}\end{center}

Indexing the rows and columns in the order $\{1,\overline{1},2,3,\dots, n-1\}$, the corresponding exchange matrix is therefore
$$\begin{bmatrix}
0 & 0  &  1 &  0  &  0  & \dots & 0  & 0 \\
0 & 0  &  1 &  0  &  0  & \dots & 0  & 0 \\
-1 & -1  &  0 & -1  &  0  & \dots & 0  & 0 \\
 0 & 0  & 1 &  0  & 1  & \dots & 0  & 0 \\
 0 & 0  &  0 & -1  &  0  & \dots & 0  & 0 \\
\dots & \dots & \dots & \dots & \dots & \dots & \dots & \dots \\
0 & 0  &  0 & 0  & 0  & \dots & (-1)^{n} & 0 \\
\end{bmatrix}.$$
We split the odd and even initial variables into the first two rows, in a zig-zagging pattern, just as before.  We then mutate in the order $1$, $\overline{1}$, $3$, $5$, $\dots$, $n$ (resp. $n-1$) if $n$ is odd (resp. even) to get the the third row, followed by mutation via $2$ then $4$, $6$, $\dots$ $n-1$ (resp. $n$) if $n$ is odd (resp. even) to get the fourth row.

The advantage of such an ordering is that the mutated exchange matrix, which we use to encode the binomial exchanges, is always the same, up to sign.  We notice that the analogue of the diamond condition for this case is $ad-bc=1$ if $b = x_i^{(j)}$ with $i \geq 2$ and 
\begin{eqnarray}
\label{Dexch1} x_{2}^{(j-1)}x_{2}^{(j+1)} - x_{1}^{(j)}x_{\overline{1}}^{(j)}x_3^{(j)} &=& 1 \\ 
\label{Dexch2} x_{1}^{(j-1)}x_{1}^{(j+1)} - x_{2}^{(j)} &=& 1 \\
\label{Dexch3} x_{\overline{1}}^{(j-1)}x_{\overline{1}}^{(j+1)} - x_{2}^{(j)} &=& 1
\end{eqnarray}
on the western boundary.

We let $\mathcal{T}_{D_n}$ be $\mathcal{T}_{B_{n-1}} \cup \{T_{\overline{1}} \}$ where $T_{\overline{1}}$ is the same tile as $T_1$ except with a different label.  We also change tile $T_2$ so that it is still a hexagon, but has weights $1, x_1, 1, x_{\overline{1}}, 1,$ and $x_3$ going around clockwise from the top.  Following the arguments of Lemmas $5$, $6$, $7$, and $8$ result in the same graph theoretic interpretation and lattice structure.  We use Rule $3$ which is analogous to Rule $2$

\begin{Rule}
Notice that when we apply Rule $1$ to set of tiles $\mathcal{T}_{D_n}$, we get a set of graphs consisting of a base of $T_2$ or $T_1 \cup T_2$ adjoining a tower of $T_a \cup \dots \cup T_b$, as before.  
We enlarge the set of graphs by allowing a base of $T_{\overline{1}} \cup T_2$ (with or without an accompanying tower), and also allow \emph{both} tile $T_1$ and tile $T_{\overline{1}}$ to appear if and only if the lift of the graph to $\mathcal{G}_{D_\infty}$ (i.e. $n$ arbitrarily large) is of the following three forms:

\begin{center} 
\includegraphics[width = 1.1in ,
height = 0.8in]{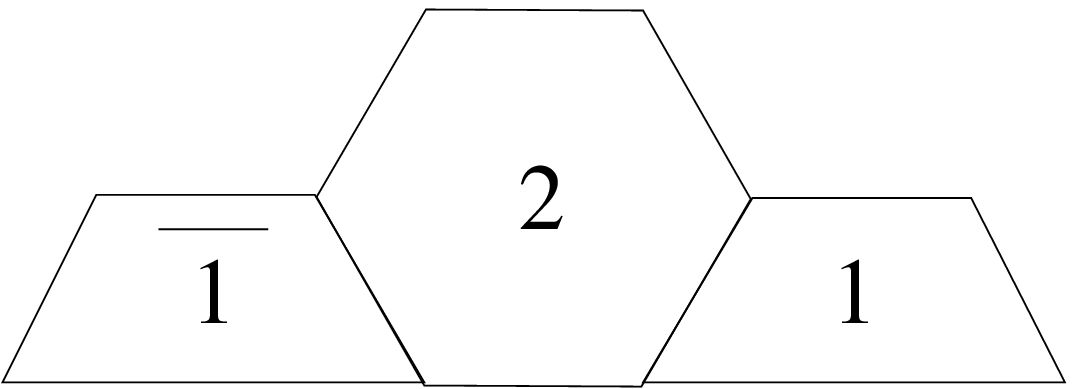}
\hspace{2em}
\includegraphics[width = 1.1in ,
height = 1.8in]{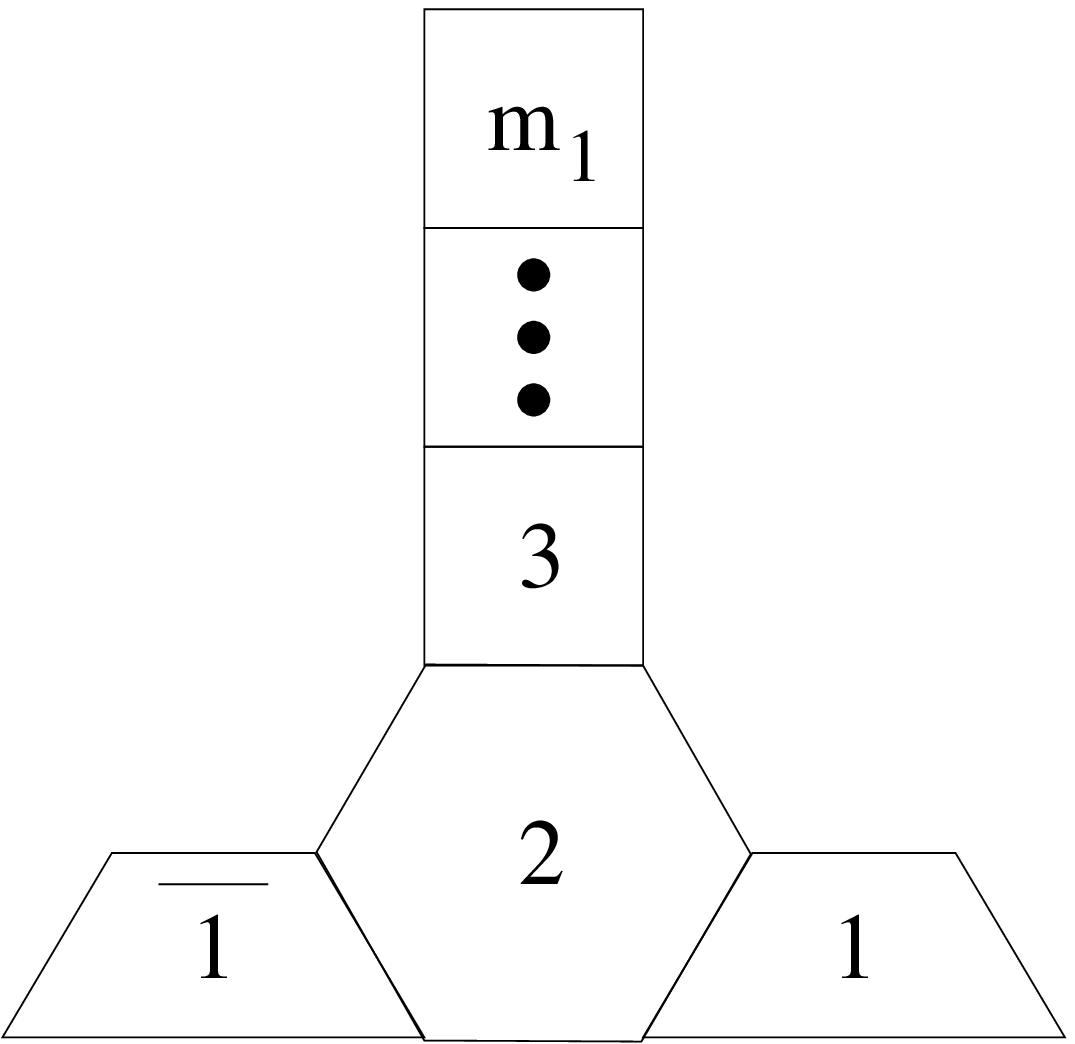}
\hspace{2em}
\includegraphics[width = 1.1in ,
height = 1.1in]{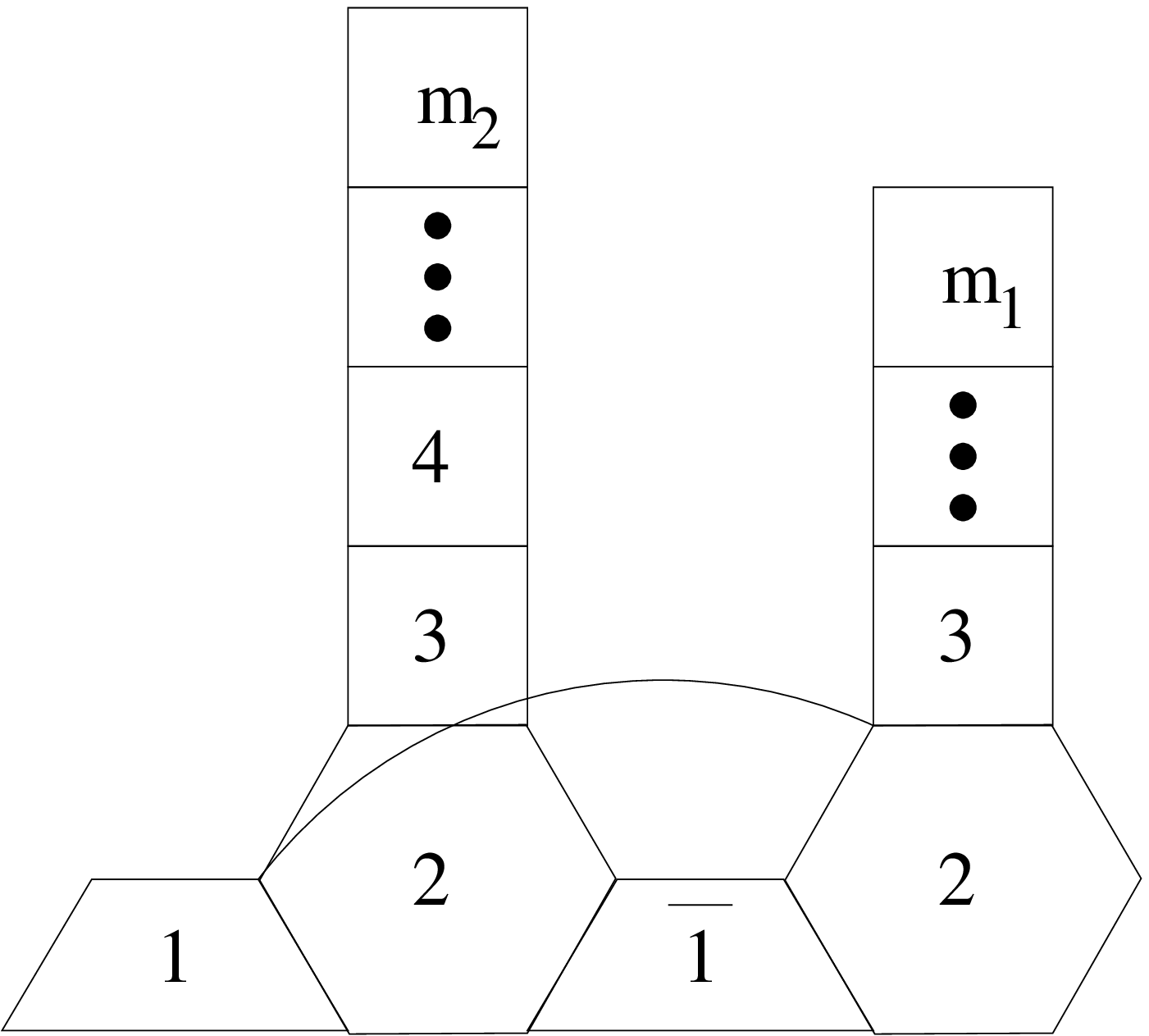}
\end{center}
where $3 \leq m_1\leq m_2$ and $m_1,~m_2$ both odd.
\end{Rule}

Let $\mathcal{T}_{D_n}$ be defined as above and $\mathcal{G}_{D_n}$ be the set of graphs constructed according to Rules $1$ and $3$.  In particular, this construction will be quite analogous to that of $\mathcal{G}_{B_{n-1}}$.  

\begin{Prop} \label{CaseDn} 
The set $\mathcal{G}_{D_n}$  is in bijection
with the set of non-initial cluster variables for a coefficient-free cluster algebra
of type $D_n$ such that the statement of Theorem \ref{vargraph}
holds.
\end{Prop}

\begin{center} \includegraphics[width = 3.5in ,
height = 1.5in]{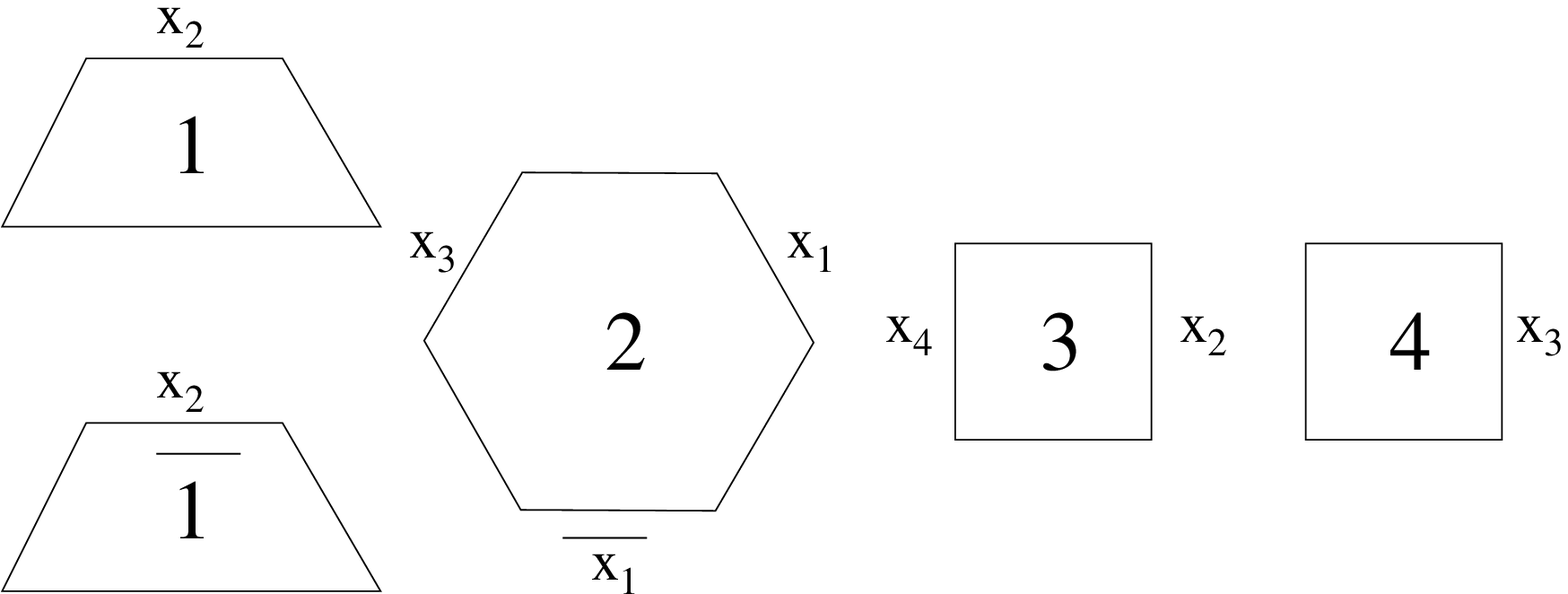} \\ Tiles for $D_5$.
\end{center}

\begin{Rem}
As indicated, the proof follows from the exact same logic as Lemmas $5$ through $8$.  The only caveat is that as a consequence of the proof, that $x_1^{(j)}$ will sometimes be a tower on base $T_1 \cup T_2$, and sometimes contain base $T_{\overline{1}} \cup T_2$.  In particular, $x_1^{(j)}$ contains $T_1$ if and only if $j$ is odd, and so we get an alternating behavior.
\end{Rem}

On the next page, we give the lattice for $\mathcal{G}_{D_5}$.  We have the usual diamond condition for four entries in three consecutive rows and three consecutive columns, not including columns one.  We encode column one by placing $x_1^{(j)}$ on top of $x_{\overline{1}}^{(j)}$, and we have the exchange relations (\ref{Dexch1}), (\ref{Dexch2}), and (\ref{Dexch3}).

\newpage $\begin{array}{cccc}
\includegraphics[width = 0.2in , height = 0.2in]{BB1.eps}  &~~~&~~~&~~~~  \\
~~~~&  ~~~& \includegraphics[width = 0.2in , height = 0.2in]{BB3.eps} &~~~ \\
\includegraphics[width = 0.2in , height = 0.2in]{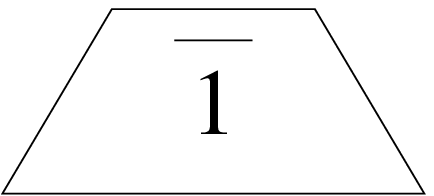}  &~~~&~~~&~~~~  \\
~~~& \includegraphics[width = 0.6in , height = 0.6in]{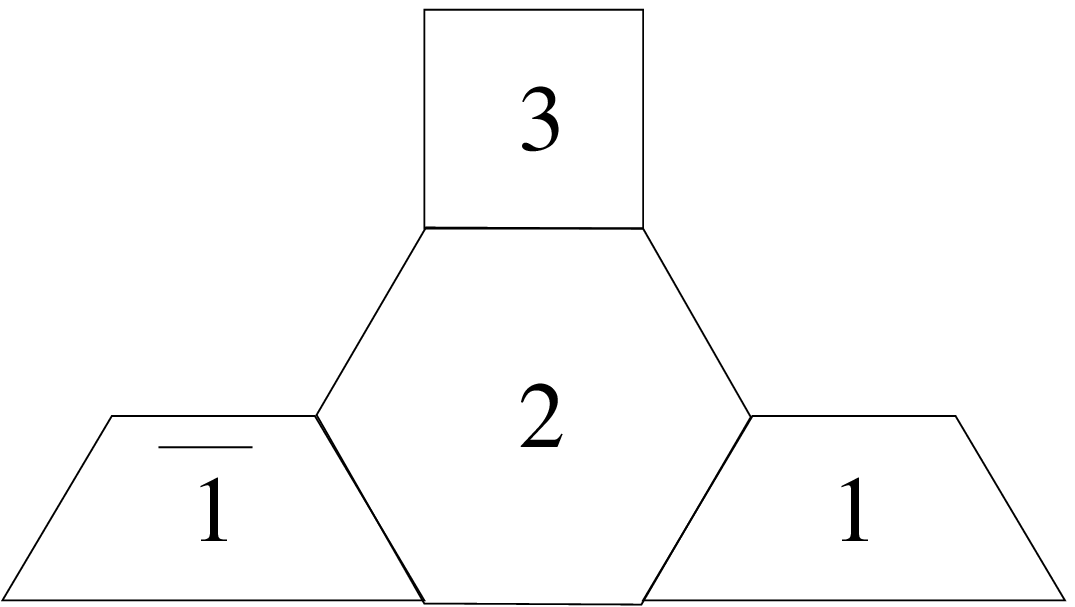} &~~~& \includegraphics[width = 0.2in , height = 0.4in]{BB34.eps}   \\
\includegraphics[width = 0.4in , height = 0.6in]{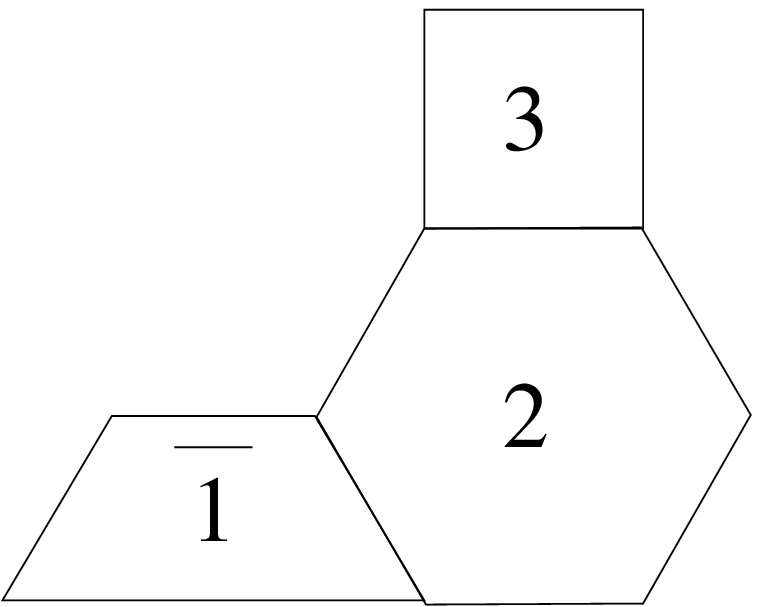}  &~~~&~~~&~~~~  \\
~~~~& ~~~& 
\includegraphics[width = 0.6in , height = 0.7in]{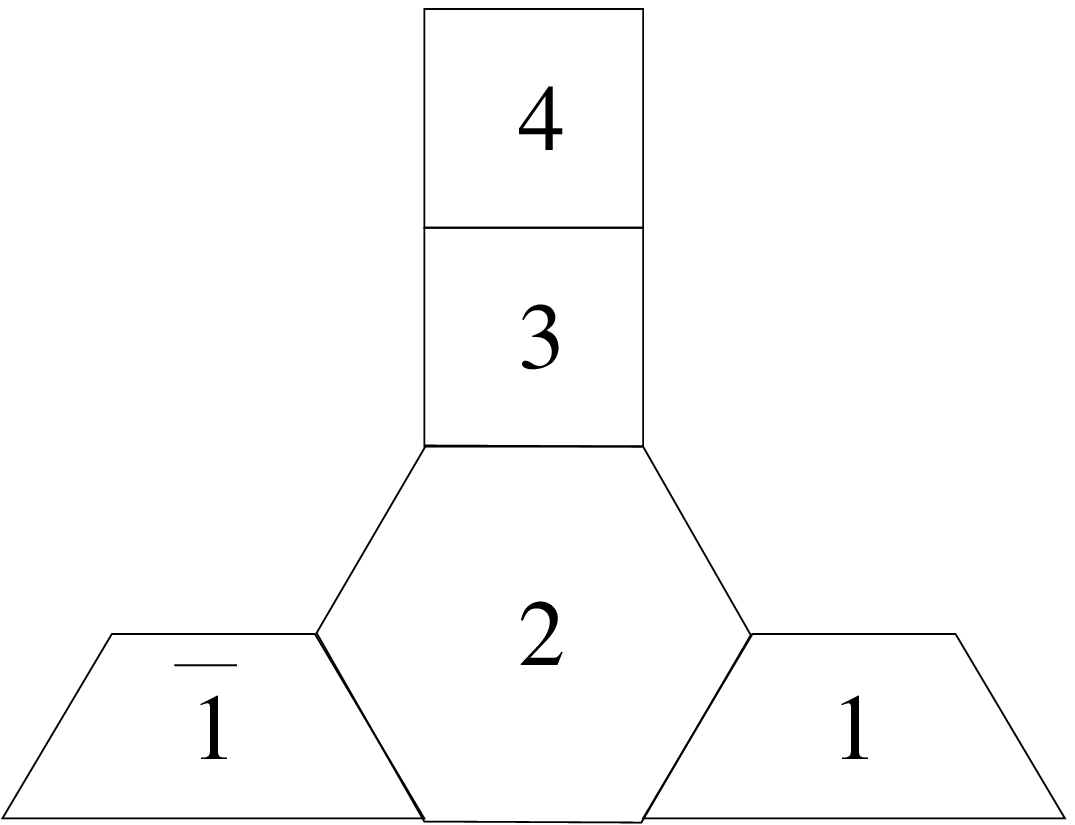} &~~~ \\
\includegraphics[width = 0.4in , height = 0.6in]{BB123.eps}  &~~~&~~~&~~~~  \\
~~~& \includegraphics[width = 0.8in , height = 0.7in]{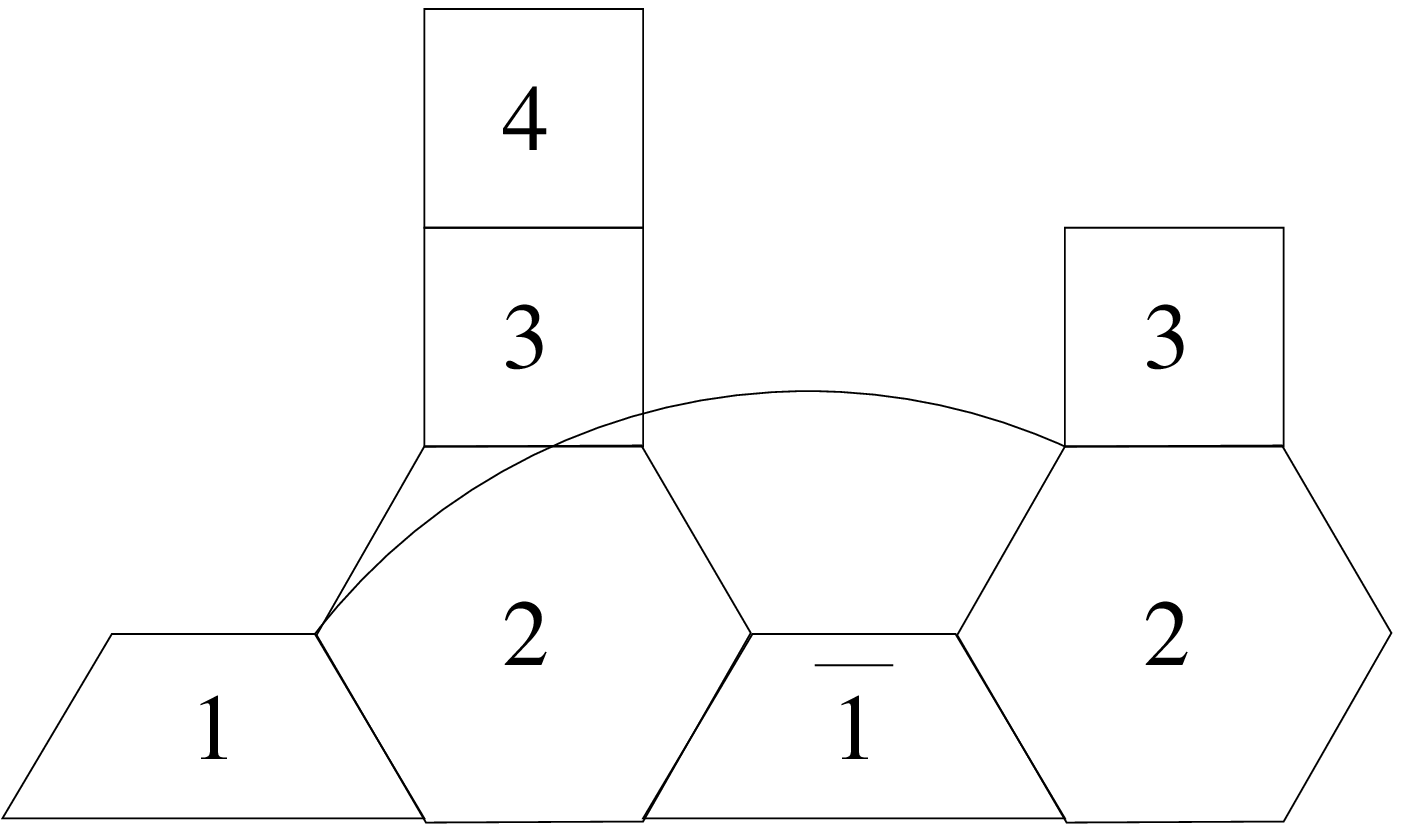} &~~~& \includegraphics[width = 0.6in , height = 0.4in]{DD121.eps}  \\
\includegraphics[width = 0.4in , height = 0.7in]{BB124.eps}  &~~~&~~~&~~~~ \\
~~~~& ~~~& \includegraphics[width = 0.8in , height = 0.6in]{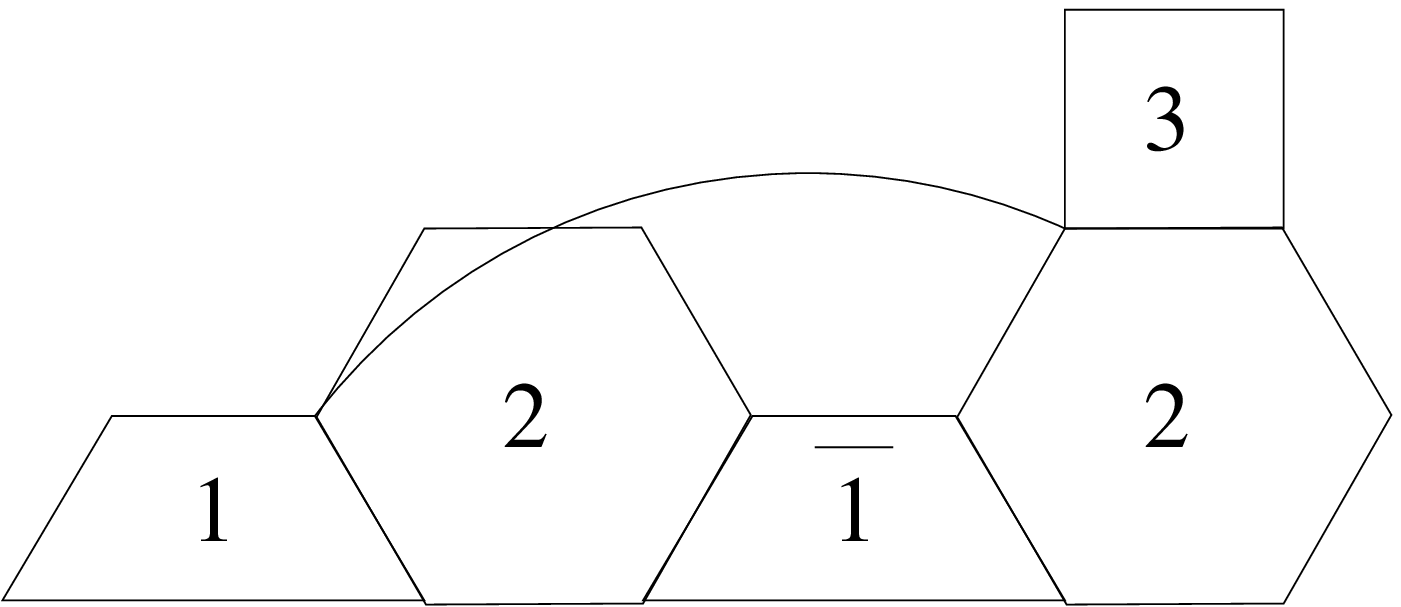} &~~~ \\
\includegraphics[width = 0.4in , height = 0.7in]{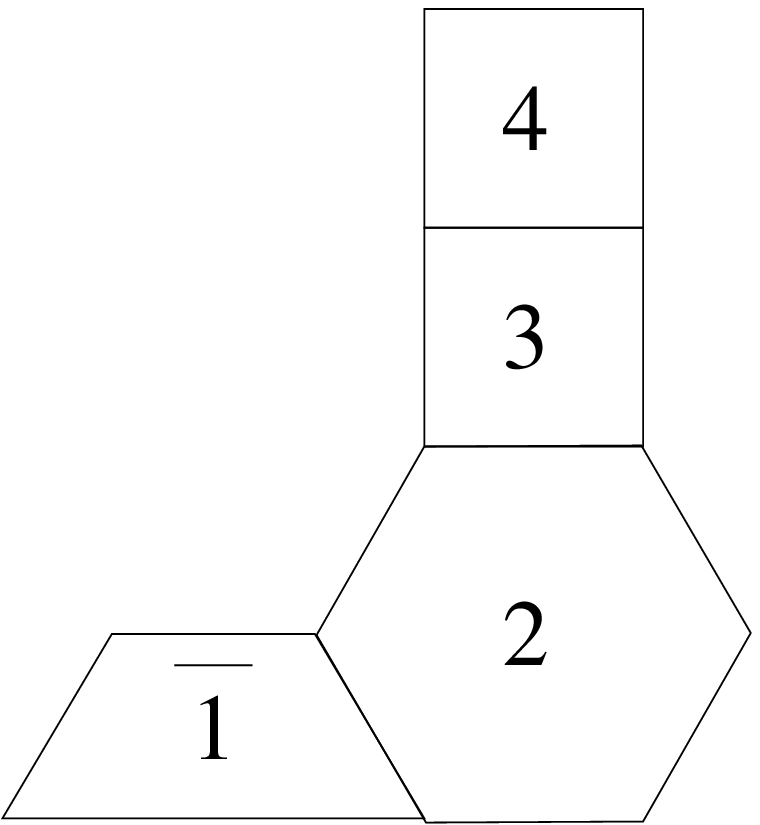} &~~~&~~~&~~~~ \\
~~~& \includegraphics[width = 0.8in , height = 0.7in]{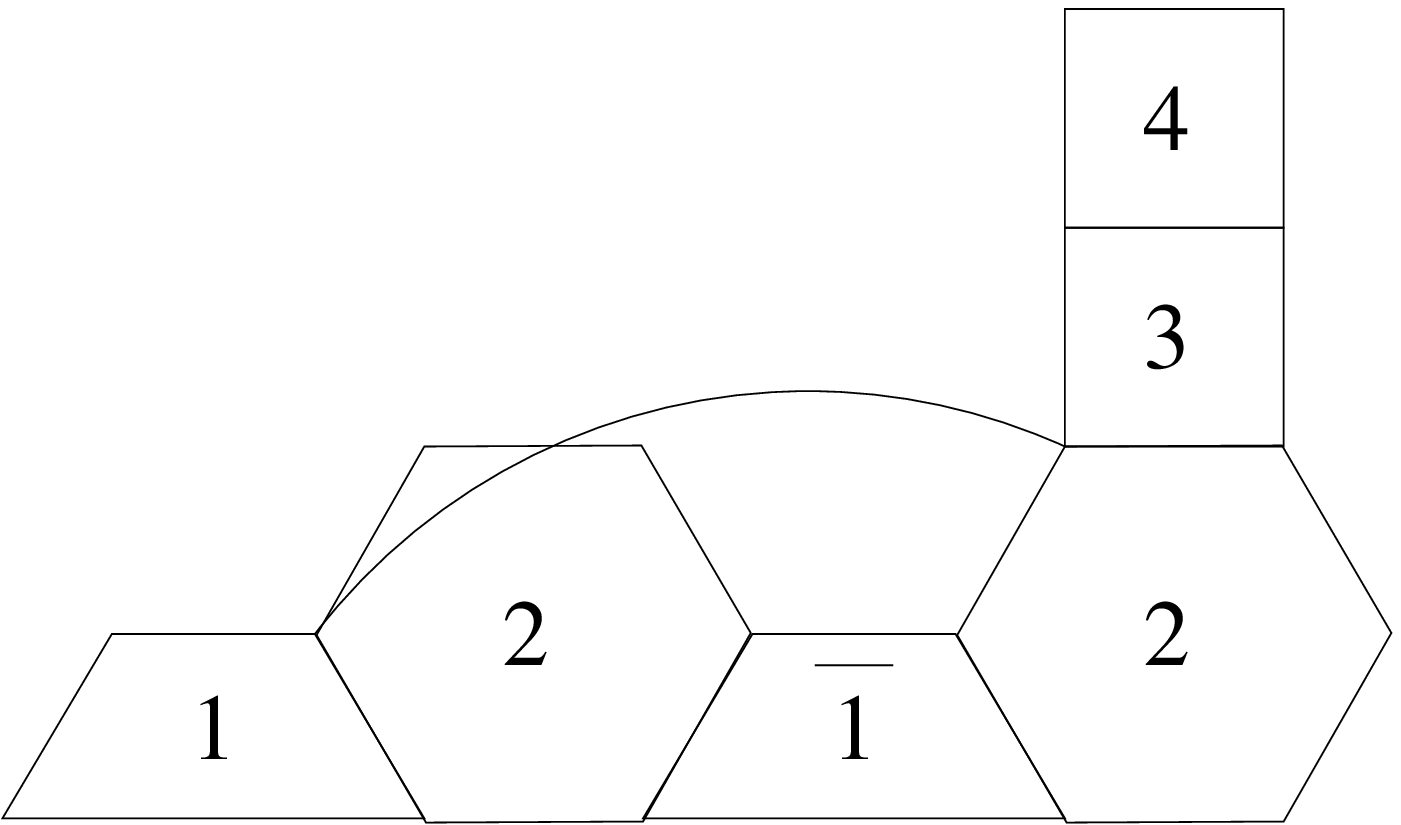} &~~~& \includegraphics[width = 0.25in , height = 0.6in]{BB23.eps}    \\
\includegraphics[width = 0.4in , height = 0.4in]{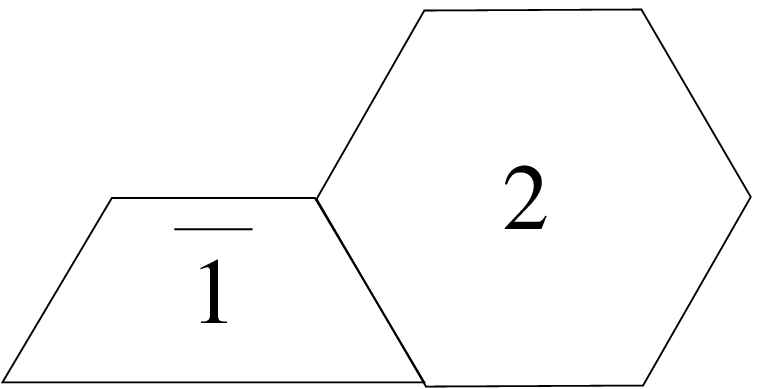}  &~~~&~~~&~~~~ \\ ~~~~& ~~~& \includegraphics[width = 0.25in , height = 0.7in]{BB24.eps} &~~~ \\
\includegraphics[width = 0.4in , height = 0.4in]{BB12.eps}  &~~~&~~~&~~~~ \\
~~~& \includegraphics[width = 0.25in , height = 0.4in]{BB2.eps} &~~~& \includegraphics[width = 0.2in , height = 0.2in]{BB4.eps}
\end{array}$

\section{$G_2$} The case of $G_2$ is the only cluster algebra of exceptional finite type for which we have been able to extend our graph theoretic interpretation.  We are able to do so since this case is analogous to $B_3$.  We use collection $\mathcal{T}_{G_2} = \{T_1,T_2\}$ with tile $T_1$ as in the $B_n$ case, and tile $T_2$ is again a hexagon, but now has all three nontrivial weights being value $x_1$.  There are six possible graphs that
correspond to the non-initial cluster variables. \vspace{1em}

\begin{center} 
$\begin{array}{ccc}
\includegraphics[width = 0.3in , height = 0.15in]{BB1.eps}  &~~~& \includegraphics[width = 0.9in , height = 0.4in]{BB121.eps} \\ \\
\includegraphics[width = 0.3in , height = 0.3in]{BB2.eps}  &~~~& \includegraphics[width = 0.9in , height = 0.75in]{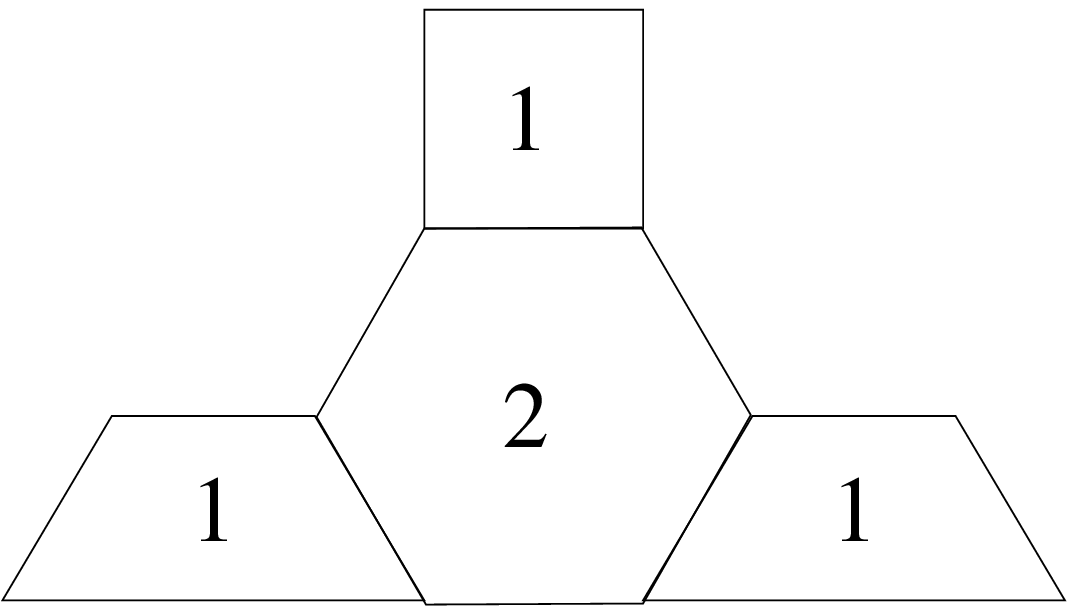} \\ \\
\includegraphics[width = 0.75in , height = 0.4in]{BB12.eps}  &~~~& \includegraphics[width = 0.9in , height = 0.75in]{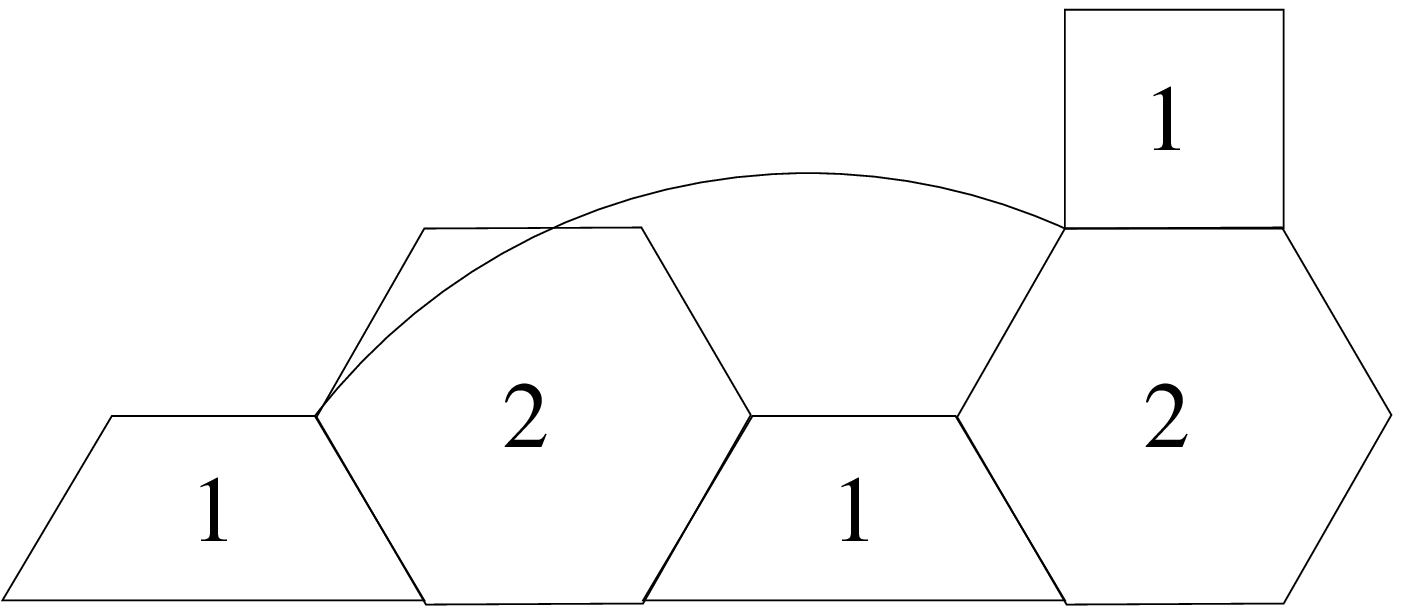}
\end{array}$
\\ Graphs for a cluster algebra of type
$G_2$. \end{center}  \vspace{1em}

\noindent $G_2$ has Dynkin diagram 
\hspace{0.8em}\includegraphics[width = 0.6in ,
height = 0.1in]{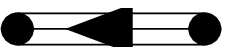}\hspace{0.8em}
and exchange matrix $\begin{bmatrix}
 0  & 1 \\
-3  &  0
\end{bmatrix}.$

\section{Future Directions}

Given the previous sections, coefficient-free cluster algebras of type $A_n$, $B_n$, $C_n$, $D_n$, or $G_2$ have a
combinatorial interpretation as a family of graphs such that the numerators of the cluster
variables enumerate the weighted number of matchings and the denominators encode the occurrences of
faces. Thus Theorem \ref{vargraph} is true in all of these cases. The next step would be to
extend Theorem \ref{vargraph} to include cluster algebras of type $E_6$, $E_7$, $E_8$, and $F_4$, and
thus have the result for all cluster algebras of finite type.

\begin{Rem}
Even though the Dynkin diagrams for the $E_n$'s are simply laced, fitting cluster algebras of these three types into patterns analogous to those of the $A_n$'s and $D_n$'s has been notoriously hard.  Such difficulties have rose elsewhere such as in the original proof of positivity in \cite{ClustII}, and also in recent models using $T$-paths on triangulated surfaces, for example in \cite{ClusIV} or \cite{ClustSurf} among other work.
\end{Rem}

Additionally, in the work of Schiffler and Carroll-Price for $A_n$, the cluster algebra considered is specifically the Ptolemy algebra, a cluster algebra \emph{with} coefficients.  In the $T$-paths model, the boundary of the polygon gives rise to $n+3$ additional coefficients which can be included in the exchange relations and cluster expansion formula.  Since the graphs we obtain in the above combinatorial interpretations are weighted so sparsely, perhaps a certain number of coefficients can be handled by the graph-model as well.

In \cite{MusPropp}, an analogous interpretation is given for rank $2$
cluster algebras of affine type and unpublished work \cite{Markoff,MusMark} done as a part of REACH, as described in \cite{MarkPropp}, gives a graph theoretic interpretation for a
totally cyclic rank $3$ cluster algebra.  This totally cyclic rank $3$ cluster algebra corresponds to a triangulated surface of genus one with exactly one puncture (i.e. interior marked point).  Such a cluster algebra has been studied geometrically including work of \cite{Thomas}.  Perhaps these graph theoretical interpretations could be extended to other cluster algebras thus providing proofs of Fomin and Zelevinsky's positivity conjecture for even further cases.

Lastly, we note that all the examples discussed above are families of \emph{planar} graphs associated to generators of cluster algebras.  When expanding our scope to include more complicated cluster algebras, is the category of planar graphs too restrictive?  More specifically, why did we need the extra arcs in the $B_n$, $D_n$ , $G_2$, and affine $A_1^{(2)}$ cases?  Perhaps it is an artifice of taking a higher dimensional object and projecting to two dimensions.

\vspace{2em}

\noindent {\bf Acknowledgments.}~~

The author would like to thank Andrei Zelevinsky for numerous helpful conversations including referring the author to \cite{YSys} where Fibonacci polynomials appear.  Discussions with Sergey Fomin, Jim Propp, Ralf Schiffler, and Hugh Thomas have also been very useful.  I wanted to especially thank Hugh Thomas and Andrei Zelevinsky for their comments on an earlier draft of this paper.

\end{document}